\documentclass[a4paper,11pt]{amsart}
\usepackage{amsmath,amsthm,amssymb}
\usepackage{graphicx}
\usepackage{subcaption}
\usepackage{float}\usepackage{slashed}
\usepackage[all]{xy}
\usepackage{cancel}
\usepackage{comment}
\newcommand{\mychoice}[2]{#1
}
\mychoice{
\newcommand{\plabel}[1]{ \label{#1}}
\newcommand{\gbibitem}[1]{ \bibitem{#1}}
\newcommand{\snewpage}{}

\specialcomment{commentx}{   }{   }\excludecomment{commentx}
}{
\usepackage{xcolor}
\newcommand{\plabel}[1]{ \label{#1}\rlap{\smash{${}^{^{[#1]}}$}}}
\newcommand{\gbibitem}[1]{ \bibitem{#1}\rlap{\smash{${}^{^{[#1]}}$}}}
\newcommand{\snewpage}{\newpage}

\newenvironment{commentx}{\color{magenta} }{\color{black} }
}

\DeclareMathOperator{\des}{des}
\DeclareMathOperator{\artanh}{artanh}

\DeclareMathOperator{\Tr}{Tr}
\DeclareMathOperator{\asc}{asc}
\DeclareMathOperator{\AC}{AC}
\DeclareMathOperator{\intD}{\mathring D}
\newcommand{\real}{\mathrm{real}}

\DeclareMathOperator{\Dbar}{D}
\DeclareMathOperator{\arsinh}{arsinh}
\DeclareMathOperator{\Rea}{Re}

\DeclareMathOperator{\Lie}{Lie}
\DeclareMathOperator{\Rexp}{exp_{R}}
\DeclareMathOperator{\Lexp}{exp_{L}}
\DeclareMathOperator{\ad}{ad}
\DeclareMathOperator{\BCH}{BCH}
\DeclareMathOperator{\spec}{sp}
\DeclareMathOperator{\Mip}{Mip}

\theoremstyle{definition}
\newtheorem{point}{}[section]

\newtheorem{method}[point]{Method}
\newtheorem{remark}[point]{Remark}
\newtheorem{example}[point]{Example}
\theoremstyle{plain}

\newtheorem{lemma}[point]{Lemma}
\newtheorem{cor}[point]{Corollary}
\newtheorem{theorem}[point]{Theorem}

\newcommand{\leaveout}[1]{}

\newcommand{\eqed}{
\pushQED{\qed}
\qedhere
\popQED
}
\newcommand{\eqedexer}{
\renewcommand{\qedsymbol}{$\diamondsuit$}
\pushQED{\qed}
\qedhere
\popQED
\renewcommand{\qedsymbol}{$\Box$}
}

\newcommand{\qedexer}{  \renewcommand{\qedsymbol}{$\diamondsuit$} \qed \renewcommand{\qedsymbol}{$\Box$}}
\newcommand{\qedremark}{  \renewcommand{\qedsymbol}{$\triangle$} \qed \renewcommand{\qedsymbol}{$\Box$}}

\newcommand{\proofremark}[1]{
\begin{proof}[Remark] #1
\renewcommand{\qedsymbol}{}
\end{proof}
}

\newcommand{\marginextend}[1]{ \addtolength{\oddsidemargin}{-#1}  \addtolength{\evensidemargin}{-#1}
  \addtolength{\textwidth}{#1}\addtolength{\textwidth}{#1}}
\newcommand{\updownextend}[1]{ \addtolength{\topmargin}{-#1}  \addtolength{\textheight}{#1}
\addtolength{\textheight}{#1}}
\marginextend{1cm}
\updownextend{0cm}
\allowdisplaybreaks[4]
\newcommand{\monom}{\mathbf M}
\newcommand{\dyh}{\mathrm d}
\newcommand{\loc}{\mathrm{loc}}
\title{Convergence estimates for the Magnus expansion III. Banach--Lie algebras}
\author{Gyula Lakos}
\address{Department of Geometry, Institute of Mathematics, E\"otv\"os University,
P\'azm\'any P\'eter s.~1/C,  Budapest, H--1117, Hungary}
\email{lakos@cs.elte.hu}
\begin{document}
\begin{abstract}
We  review and provide simplified proofs related to the Magnus expansion, and improve convergence estimates.
Observations and improvements concerning the Baker--Campbell--Hausdorff expansion are also made.

In this Part III, we consider the Banach--Lie algebraic setting.
We show how to improve the ``standard'' convergence bound $\boldsymbol\delta= 2.1737374\ldots$
using the customary generating function / ODE methods.
Then we discuss how to achieve better convergence bounds  using the resolvent method.
The emphasis here is on the variety of the methods, rather than pushing them to the extreme.
Nevertheless, regarding the cumulative convergence radii, we show how to establish
$2.427<\mathrm C_\infty^{\Lie}\leq 4$
for the Magnus expansion, and
$2.93<\mathrm C^{\Lie}_2 \leq 2\boldsymbol v_{\mathrm{Mi}}=5.4028570\ldots$
for the BCH expansion.
\end{abstract}
\maketitle

\section*{Introduction}
This paper is a direct continuation of Part I, \cite{L1}.
We assume general familiarity with the results and techniques presented there.
General sources for algebra, analysis, and combinatorics should also be taken from there.\\

\textbf{Introduction to the Banach--Lie algebraic setting.}
From the very beginning, Magnus \cite{M} (1954), it was realized that the Magnus expansion can be interpreted in Lie algebraic terms.
Indeed, it was conceived as the continuous generalization of the Baker--Campbell--Hausdorff expansion.
Apart form classical material on the BCH formula, exposed in Bonfiglioli, Fulci \cite{BF} (see also Achilles, Bonfiglioli \cite{AB}),
the related basic algebraic properties were clarified by Magnus \cite{M}, proving the Magnus recursion formula,
somewhat earlier by Dynkin \cite{Dyy},
and later by Goldberg \cite{G}, Solomon \cite{S},
Mielnik, Pleba\'{n}ski \cite{MP}, Helmstetter \cite{H}, cf. Reutenauer \cite{R}.

Convergence in the Lie algebraic setting can be studied using Lie norms, which are norms $\|\cdot\|$
such that $\|[X,Y]\|\leq\|X\|\cdot\|Y\|$ holds.
It is also natural to assume that the underlying space is a Banach space with respect to $\|\cdot\|$.

Regarding earlier works on convergence in this Banach--Lie setting,
the well-known lower estimate for the convergence radius of the Magnus expansion, in terms of the cumulative norm,  is the
Varadarajan--M\'erigot--Newman--So--Thompson number
\[\boldsymbol\delta=\int_{x=0}^{2\pi} \frac{\mathrm dx}{2+\frac x2-\frac x2\cot \frac x2 }= 2.1737374\ldots.\]
This number makes its appearance  in the setting of Baker--Campbell--Hausdorff formula,
due to the work of  Varadarajan \cite{Va} (1974), M\'erigot \cite{Me} (1974), Newman, So, Thompson \cite{NST}.
It was extended to the Magnus expansion setting by Blanes, Casas, Oteo, Ros \cite{BCOR1}, Moan \cite{Mo1} (1998).
It is the  bound generally cited in the literature, thus we call it as the `standard estimate'.
It is, however, but the trivial estimate which can be obtained from the Magnus recursion formulas.
Furthermore, it is very easy to improve, even if just by a little bit.
We explain this in Section \ref{sec:MagnusStandard}.
Regarding divergence, as a consequence of the Banach algebraic case, one
can immediately obtain counterexamples of cumulative norm $4$ (infinite dimensional case),
or just larger (finite dimensional case).

It was also quickly realized the standard estimate can be improved in the special case of the Baker--Campbell--Hausdorff expansion.
M\'erigot \cite{Me} (1974) proves that the BCH expansion with respect to $X,Y$ converges if
$\|X\|,\|Y\|<2\pi$ and
\[\|X\|<\int_{x=\|Y\|}^{2\pi} \frac{\mathrm dx}{2+\frac x2-\frac x2\cot \frac x2 }
\qquad\text{and}\qquad
\|Y\|<\int_{x=\|X\|}^{2\pi} \frac{\mathrm dx}{2+\frac x2-\frac x2\cot \frac x2 }
.\]
Unfortunately, \cite{Me} is not truly published, but its results are reported
in Michel \cite{Mi}, including a graphical representation;
and its content is exposed in Biagi \cite{Bia}; cf. also Biagi, Bonfiglioli \cite{BB}.
It is redeveloped by Day, So, Thompson \cite{DST}, Blanes, Casas \cite{BC}.
In particular, the BCH expansion converges if
$\|X\|,\|Y\|<\boldsymbol\delta_{\BCH}/2= 1.2357524\ldots$,
where $\boldsymbol\delta_{\BCH}/2$ is the solution of the equation
\[\boldsymbol\delta_{\BCH}/2=\int_{x=\boldsymbol\delta_{\BCH}/2}^{2\pi} \frac{\mathrm dx}{2+\frac x2-\frac x2\cot \frac x2 }.\]
(The numerical value was exhibited first by  Day, So, Thompson \cite{DST}.)
However, from the graphical representation, and its subsequent analysis,
it is also clear (cf.~Michel \cite{Mi}, Blanes, Casas \cite{BC}) that the BCH expansion converges for the cumulative norm
$\|X\|+\|Y\|<\boldsymbol\delta_{\BCH}= 2.4715048\ldots$.
On the other hand,  the Dynkin--Specht--Wever lemma implies
that if we have power series based \textit{convergence} estimates
(i.~e.~estimates of the convergence radius from below)
in the Banach algebraic setting, then they can essentially be transferred into the Banach-Lie-algebraic   setting
(cf. Dynkin \cite{Dy}, Strichartz \cite{St}).
From this viewpoint, according to Part I, the `standard estimate' presents an improvement relative to the Magnus case,
where the Banach algebraic cumulative convergence radius is $2$.
This is not so in the BCH case:
We know convergence in the plain Banach case for $|X|+|Y|< \mathrm C_2=2.89847930\ldots$;
thus, according to the Dynkin argument, we also know that
convergence holds in the Banach--Lie setting with $\|X\|+\|Y\|< \mathrm C_2$.
(A stronger result is possible in the unbalanced case, but in Part I, the analysis was not carried to those details.
Anyway, depending on the norm ration, what we obtain from the Banach algebraic case is converge for a cumulative norm
less than a value which is between $\mathrm C_2$ and $\pi$, the latter cumulative norm failing even in the case Hilbert space operators.)
Divergence (that is estimates of the convergence radius from above)
in the Banach--Lie algebraic setting was studied most eminently
by Michel \cite{Mi} although without numerical values given (but providing some graphical representation).
In general,  he finds that there are counterexamples to the convergence
of the BCH expansion regarding $X$, $Y$ such that even if arbitrary positive norm ratio $\|X\|:\|Y\|$
is prescribed, $\|X\|+\|Y\|<2\pi$ can be achieved.
 Most notably, there is a counterexample with $\|X\|=\|Y\|=\boldsymbol v_{\mathrm{Mi}}$ for which we will provide
the numerical value $\boldsymbol v_{\mathrm{Mi}}=2.70142851\ldots$.
For the reader's convenience, in the following Figure 1, the we redraw the ``classical'' convergence results
for the BCH expansion. Note, however, that the consequences which can be drawn from the Banach algebraic
case are represented in an imperfect form, as the unbalanced BCH problem with restricted norm ratios was not investigated  in \cite{L1}.
\begin{figure}[H]
   \begin{subfigure}[b]{5.0in}
    \centerline{ \includegraphics[width=3.0in]{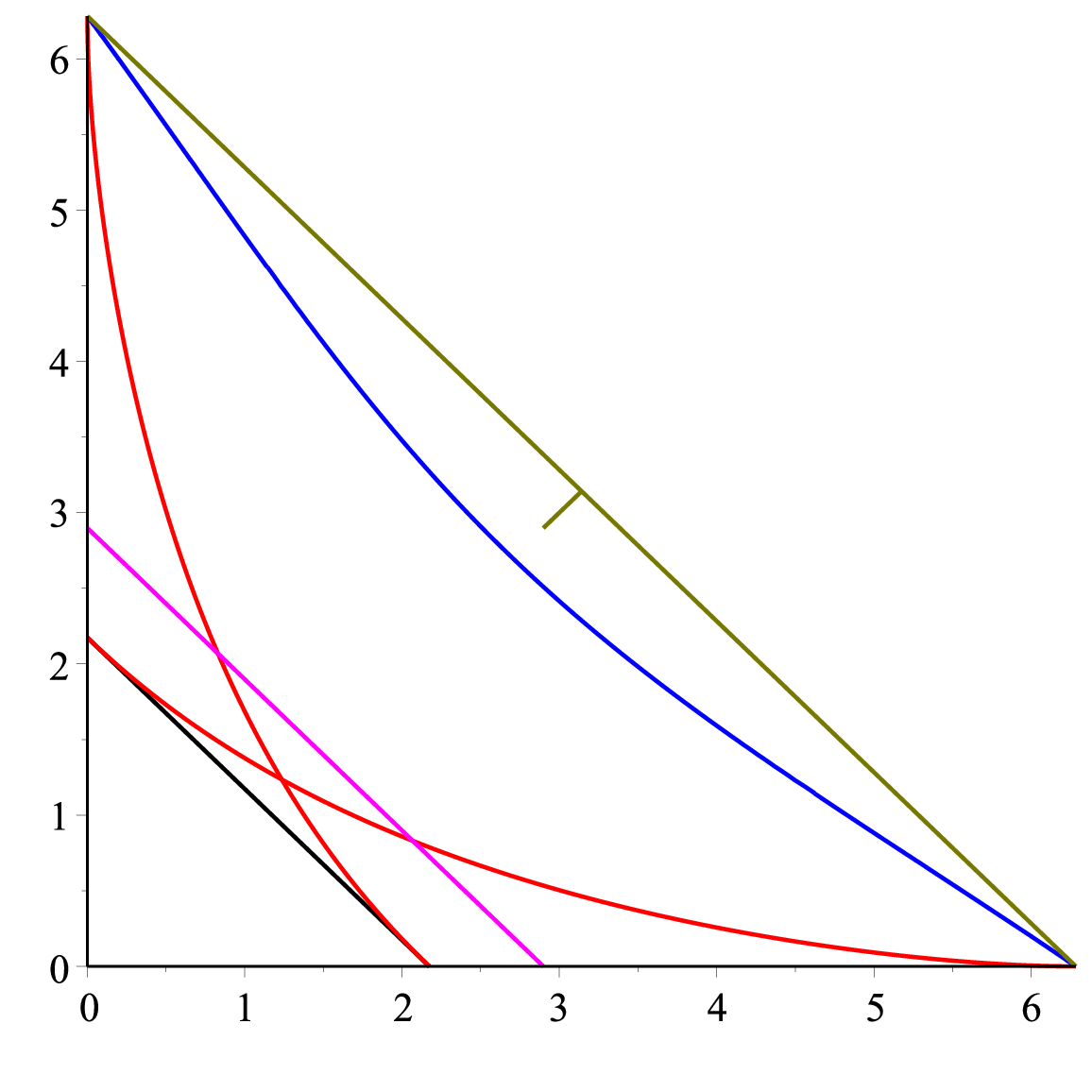} }
    \caption*{Fig.~\ref{fig:figC01}. Some classical convergence estimates for the BCH expansion in terms of
     the norms $\|X\|$ and $\|Y\|$. Convergence estimates:
     Black line: Varadarajan's result.
     Red lines: M\'erigot's result.
     Magenta line: A consequence of  the Banach algebraic result.
     Below those lines convergence is assured.
     Divergence estimates: Olive line: Results following from the Banach algebraic case.
     Blue line: Michel's results. At and beyond those lines divergence is possible.}
  \end{subfigure}
\phantomcaption
\plabel{fig:figC01}
\end{figure}~\\[-0.5cm]

\textbf{Outline of content.}
Our objective, in this Part III, is to improve the convergence / divergence
bounds regarding the Magnus / BCH expansions in the Banach--Lie setting.
Here the emphasis is on the variety of methods.
Unfortunately, all these methods are somewhat computational.
There was no intention to take any of them ``to the extreme''.
As a general phenomenon, all strict inequalities we give are quite easy to improve (a bit).
Nevertheless, regarding the cumulative convergence radii, we show how to establish
$2.427<\mathrm C_\infty^{\Lie}\leq 4$
for the Magnus expansion, and
$2.93<\mathrm C^{\Lie}_2 \leq 2\boldsymbol v_{\mathrm{Mi}}=5.4028570\ldots$
for the BCH expansion ($\boldsymbol v_{\mathrm{Mi}}$ being ``Michel's number'').

Section \ref{sec:MagnusBanachLie} discusses the Magnus commutators from Lie algebraic viewpoint; this is basically algebra.
Section \ref{sec:Lieconv} deals with principles of convergence in the Banach--Lie setting, and establishes the most basic results.
Section \ref{sec:bchupper} presents upper estimate for the convergence radius of the BCH expansion.
Section \ref{sec:MagnusStandard} discusses improvements of  the standard convergence bound $\boldsymbol\delta= 2.1737374\ldots$.
(These estimates will be superseded by the resolvent method later but recursive methods should not be abandoned altogether.)
Section \ref{sec:normext} is about how to consider the general Banach--Lie algebraic case as a Banach algebraic problem.
Section \ref{sec:resmagnus} gives lower estimates for the convergence radius of the Magnus expansion using the resolvent method.
Section \ref{sec:resbch} gives lower estimates in the case of the  BCH expansion using a simplified setting.
Section \ref{sec:numerics} compares the estimates above to the data from the appendices.
Appendix \ref{sec:someTheta} gives the exact Lie-$\ell^1$ norm of some Magnus commutators.
Appendix \ref{sec:someGamma} gives the exact Lie-$\ell^1$ norm of some BCH terms.

\textbf{Acknowledgements.}
The author would like thank Alp\'ar J\"uttner for some advices regarding computational issues.
\begin{commentx}
\tableofcontents
\end{commentx}
\snewpage
\section{The Magnus commutators as Lie polynomials}\plabel{sec:MagnusBanachLie}
Recall, the Magnus commutators $\mu_k(X_1,\ldots,X_k)$ are the 1-homogeneous noncommutative
polynomials in $X_1,\ldots,X_k$ defined by
\begin{equation}\mu_k(X_1,\ldots,X_k)=
\frac{\partial^k}{\partial t_1\cdot\ldots\cdot\partial t_k}\log(\exp(t_1X_1)\cdot\ldots\cdot\exp(t_kX_k))\Bigl|_{t_1=\ldots=t_k=0};
\plabel{eq:mudef}
\end{equation}
or, by
\begin{equation}\mu_k(X_1,\ldots,X_k)=\log(\exp(X_1)\cdot\ldots\cdot\exp(X_k))_{\text{the variables }
X_1\ldots X_k\text{ has multiplicity }1};\plabel{eq:sinmu}\end{equation}
or, by
\begin{equation}\log(\exp(X_1)\cdot\ldots\cdot\exp(X_k))=\sum_{\substack{\{i_1,\ldots, i_l\}\subset \{1,\ldots,k\}  \\ i_1<\ldots <i_l }}
\mu_l(X_{i_1},\ldots,X_{i_l}) +H(X_1,\ldots,X_k),\plabel{eq:bare}
\end{equation}
where $H(X_1,\ldots,X_k)$ collects higher multiplicities; or, by
\begin{equation}
X_1\cdot\ldots\cdot X_k=\sum_{\substack{I_1\dot\cup\ldots \dot\cup I_s=\{1,\ldots,k\}\\I_j=\{i_{j,1},\ldots,i_{j,l_j}\}\neq
\emptyset \\i_{j,1}<\ldots<i_{j,l_j}}}\frac1{s!}\cdot
\mu_{l_1}(X_{i_{1,1}},\ldots ,X_{i_{1,l_1}})\cdot \ldots \cdot\mu_{l_s}(X_{i_{s,1}},\ldots ,X_{i_{s,l_s}}).
\plabel{eq:ppod}
\end{equation}
Indeed, it easy to see that for \eqref{eq:mudef}, \eqref{eq:sinmu}, \eqref{eq:bare}, \eqref{eq:ppod}, any
line is quite synonymous to the next one.
Less trivially, according to the formula of Dynkin \cite{Dyy} or Mielnik, Pleba\'nski \cite{MP},
\begin{equation}\mu_k(X_1,\ldots,X_k)=
\sum_{\sigma\in\Sigma_k}(-1)^{\des(\sigma)}
\frac{\asc(\sigma)!\des(\sigma)!}{k!}X_{\sigma(1)}\cdot\ldots\cdot X_{\sigma(k)} ;
\plabel{eq:mpform}
\end{equation}
 where $\asc(\sigma)$ denotes the number  ascents,
and $\des(\sigma)$ denotes the number of its descents in the permutation $\sigma$.

Still in the setting of (Banach) algebras, let $\ad X$ denote the operator
$\ad X\,:\,\mathfrak A\rightarrow\mathfrak A,$ given by
$Y\mapsto [X,Y]$. Consider the meromorphic function
\begin{equation}\beta(x)=\frac x{\mathrm e^x-1}=\sum_{j=0}^\infty \beta_j x^j.\plabel{eq:schur}\notag\end{equation}
Note that this function has poles at $2\pi\mathrm i(\mathbb Z\setminus\{0\})$.
\begin{theorem}[F. Schur \cite{Sch1} (1890), \cite{Sch2} , Poincar\'e \cite{PH} (1899) ]
If $|X|<\pi$, or $\spec(X)\subset\{z\in\mathbb C\,:\,|z|<\pi\}$, or $\spec(X)\subset\{z\in\mathbb C\,:\,|\Rea z|<\pi\}$ , then
$|\ad X|_{\mathfrak A}<2\pi$, or  $\spec(\ad X)\subset\{z\in\mathbb C\,:\,| z|<2\pi\}$, or  $\spec(\ad X)\subset\{z\in\mathbb C\,:\,|\Rea z|<2\pi\}$, respectively.
In particular,
$\beta(\ad X):\mathfrak A\rightarrow \mathfrak A$ makes sense as an absolute convergent power series (first two cases)
 or as a homomorphic function of  $\ad X$.

In these cases, for $Y\in\mathfrak A$,
\begin{equation}
\frac{\mathrm d}{\mathrm dt}\log(\exp(tY)\exp(X))\Bigl|_{t=0}=\beta(\ad X)Y
\plabel{eq:schurL}\end{equation}
and
\begin{equation}
\frac{\mathrm d}{\mathrm dt}\log(\exp(X)\exp(tY))\Bigl|_{t=0}=\beta(-\ad X)Y
\plabel{eq:schurR}\end{equation}
hold;
with the usual $\log$ branch cut along the negative real axis.
\begin{proof}
This is discussed in Part I in detail.
\end{proof}
\end{theorem}

These formulae also extend, mutatis mutandis, for Banach--Lie groups.
However, at this point we are not interested in generalizations, but in the formal consequence:

\begin{cor}
If $X$ and $Y$ are formal noncommutative variables, then
\begin{equation}
\log(\exp(Y)\exp(X))_{\text{the multiplicity of $Y$ is $1$}}=\beta(\ad X)Y
\plabel{eq:schur2L}
\end{equation}
and
\begin{equation}
\log(\exp(X)\exp(Y))_{\text{the multiplicity of $Y$ is $1$}}=\beta(-\ad X)Y ;
\plabel{eq:schur2R}
\end{equation}
where $\beta(\ad X)$ is understood in the sense of formal power series.
\end{cor}

An immediate consequence is the Magnus recursion theorem.
\begin{theorem}[Magnus \cite{M}, 1954] \plabel{th:Magnusrec}
$\mu_k$ satisfies the recursions
\begin{multline}
\mu_k(X_1,\ldots,X_k)=\\
=\sum_{\substack{I_1\dot\cup\ldots \dot\cup I_s=\{2,\ldots,k\}\\I_j=\{i_{j,1},\ldots,i_{j,l_j}\}\neq
\emptyset \\i_{j,1}<\ldots<i_{j,l_j}}}\beta_s\cdot (\ad \mu_{l_1}(X_{i_{1,1}},\ldots ,X_{i_{1,l_1}}))\ldots (\ad \mu_{l_s}(X_{i_{s,1}},\ldots ,X_{i_{s,l_s}}))X_1
\plabel{eq:magL}
\end{multline}
and
\begin{multline}
\mu_k(X_1,\ldots,X_k)=\\
= \sum_{\substack{I_1\dot\cup\ldots \dot\cup I_s=\{1,\ldots,k-1\}\\I_j=\{i_{j,1},\ldots,i_{j,l_j}\}\neq
\emptyset \\i_{j,1}<\ldots<i_{j,l_j}}}(-1)^s\beta_s\cdot (\ad \mu_{l_1}(X_{i_{1,1}},\ldots ,X_{i_{1,l_1}}))\ldots (\ad \mu_{l_s}(X_{i_{s,1}},\ldots ,X_{i_{s,l_s}}))X_k.
\plabel{eq:magR}
\end{multline}

In particular, we find that $\mu_k(X_1,\ldots,X_k)$ is a commutator polynomial of its variables.
\proofremark{
This is combinatorially equivalent to the original version formulated  in a more ODE looking setting.
}
\begin{proof} Consider the first equation.
Let us apply the first formal Schur formula with $X=\log(\exp(X_2)\cdot\ldots\cdots \exp(X_{k}))$ and $Y=X_1$, and select the
terms where every variable $X_i$ has multiplicity $1$.
Considering \eqref{eq:bare} yields the formula immediately.
The second equation is analogous.
\end{proof}
\end{theorem}
\begin{remark}
The fact that $\mu_k$ is a commutator polynomial can also be established
using the PBW theorem and applying the Friedriechs trick $X_i\mapsto X_i\otimes 1+1\otimes X_i$,
cf. Reutenauer \cite{R} or Bonfiglioli, Fulci \cite{BF}.
\qedremark
\end{remark}

As $\mu_k(X_1,\ldots,X_k)$ is a commutator polynomial, we can associate a Lie polynomial
$\mu_k^{\Lie}(X_1,\ldots,X_k)$ to it.
According to the representability theorem of Magnus (which can be proven directly, but it is also a consequence of the
Poincar\'e--Birkhoff--Witt theorem over a field), the Lie polynomial itself is independent of which
commutator expression is used.
More generally, if $\mathfrak g$ is a Lie algebra over a field of characteristic $0$, then
we have natural maps
\[\xymatrix{\odot \mathfrak g\ar[r]_\iota&\ar@/_/[l]_\varpi\otimes \mathfrak g\ar[r]_U&\mathcal U\mathfrak g}.\]
Here $\iota(X_1\odot\ldots\odot X_n)=\frac1{n!}\sum_{\sigma\in\Sigma_n}X_{\sigma(1)}\otimes\ldots\otimes X_{\sigma(n)}$,
and $U$ is the factorization generated by $X_1\otimes X_2-X_2\otimes X_1=[X_1,X_2]$.
By the Poincar\'e--Birkhoff--Witt theorem, $U\circ \iota$ is a linear  isomorphism, and we can set $\varpi=(\iota\circ U)^{-1}\circ U$.
Let $\varpi_k$ denote the further projection of $\varpi$ to $\odot^k\mathfrak g$; this is the $k$-th canonical projection.

\begin{theorem}[Solomon \cite{S}, 1968] \plabel{th:solomon}
The first canonical projection is given by Magnus commutators:
\[
(U\circ \iota)\circ \varpi_1(X_1\otimes\ldots\otimes X_n)=\mu_n(X_1,\ldots,X_n)_{\,\mathrm{in}\,\mathcal U\mathfrak g},
\]
and
\[
\varpi_1(X_1\otimes\ldots\otimes X_n)=\mu_n^{\Lie}(X_1,\ldots,X_n).
\]
\begin{proof}(Cf. Helmstetter \cite{H}, 1989.)
Formula \eqref{eq:ppod} decomposes any product to  symmetric products of commutator expressions.
The part of  symmetric degree $1$ is exactly the Magnus commutator.
\end{proof}
\proofremark{
Solomon \cite{S} computes $(U\circ \iota)\circ\varpi_1$ directly to the RHS of \eqref{eq:mpform}.
Helmstetter \cite{H} understands the connection to the $\log\Pi\exp$-structure,
but does not mention the Magnus expansion as such.
Reutenauer \cite{R} has the full picture algebraically.
Equation  \eqref{eq:ppod} also shows how to express the higher canonical projections with Magnus commutators.
In fact, already Dynkin \cite{Dyy} (1949) has the Magnus commutator pretty much in our spirit,
but which he applies only to the BCH expansion.
}
\end{theorem}

\begin{cor}\plabel{cor:genMagnusrec} For $1\leq h_1,h_2\leq k$, $h_1+h_2\leq k$,
\begin{multline}
\mu_k(X_1,\ldots,X_k)
=\sum_{\substack{I_1\dot\cup\ldots \dot\cup I_s=\{h_1+1,\ldots,k-h_2\}\\I_j=\{i_{j,1},\ldots,i_{j,l_j}\}\neq
\emptyset \\i_{j,1}<\ldots<i_{j,l_j}}}\\
\mu_{s}(X_1,\ldots,X_{h_1},\mu_{l_1}(X_{i_{1,1}},\ldots ,X_{i_{1,l_1}}),\ldots, \mu_{l_s}(X_{i_{s,1}},\ldots ,X_{i_{s,l_s}}),
X_{k-h_2+1},\ldots,X_{k}
).\plabel{eq:gemma}
\end{multline}
\begin{proof}
This follows from applying the first canonical projection to
\begin{multline}
\log(\exp(X_1)\cdot\ldots\cdot\exp(X_{k}))= \log(\exp(X_1)\cdot\ldots\cdot\exp(X_{h_1})\cdot\\\exp (\log(\exp(X_{h_1+1})\cdot\ldots\cdot\exp(X_{k-h_1}))
\cdot\exp(X_{k-h_2+1})\cdot\ldots\cdot\exp(X_{k}) ).
\notag\qedhere
\end{multline}
\end{proof}
\begin{proof}[Remark]
\begin{commentx}

Equation \eqref{eq:gemma}, in case $h_1+h_2=1$,
relates to the classical Magnus recursion
and the polarized Schur identities
\[\sum_{\sigma\in\Sigma_{\{2,\ldots,n\}}}\mu(X_1,X_{\sigma(2)},\ldots,X_{\sigma(n)})=
\sum_{\sigma\in\Sigma_{\{2,\ldots,n\}}}\beta_{n-1}  (\ad X_{\sigma(2)}) \ldots (\ad X_{\sigma(n)})X_1,\]
\[\sum_{\sigma\in\Sigma_{\{1,\ldots,n-1\}}}\mu(X_{\sigma(1)},\ldots,X_{\sigma(n-1)},X_n)=
\sum_{\sigma\in\Sigma_{\{1,\ldots,n-1\}}}\beta_{n-1}  (\ad X_{\sigma(1)}) \ldots (\ad X_{\sigma(n-1)})X_n.\]

\end{commentx}
One can obtain many  identities in similar way; cf. Reutenauer \cite{R}.
\renewcommand{\qedsymbol}{$\triangle$}
\end{proof}
\renewcommand{\qedsymbol}{$\Box$}
\end{cor}

Using Theorem \ref{th:Magnusrec}, one can compute $\mu^{\Lie}_k$ effectively.
Another possibility is to use the Dynkin--Specht--Wever lemma in order to turn
\eqref{eq:mpform} into  an explicit Lie polynomial; as it was already done in Mielnik, Pleba\'nski, \cite{MP}.
Then, we obtain,
\begin{equation}\mu_k^{\Lie}(X_1,\ldots,X_k)=
\sum_{\sigma\in\Sigma_k}(-1)^{\des(\sigma)}
\frac{\asc(\sigma)!\des(\sigma)!}{k!k}(\ad X_{\sigma(1)})\cdot\ldots\cdot (\ad X_{\sigma(k-1)})  X_{\sigma(k)} .\notag
\end{equation}
(Here `$\ad$' was understood on the level of Lie algebras.)
There are several versions of this procedure, cf. \cite{L0}.

In fact, there is no obligation to use Schur's argument, or the Poincar\'e--Birkhoff--Witt lemma,
or the Dynkyn--Specht--Wever lemma in the way we used them.
It is possible to define the Magnus commutators directly as Lie polynomials,
using \eqref{eq:magL} and \eqref{eq:magR}, and then deduce the above mentioned classical statements (over a field of characteristic $0$)
using the Magnus commutators.
See \cite{L0} for this.

In what follows, we use the notation $\mu_k^{\Lie}$ only when it is particularly meant to be interpreted as a Lie polynomial.
Let us mention some algebraic properties of the Magnus commutators.
\begin{lemma}
\begin{equation}\mu_k(X_1,\ldots,X_k)=-\mu_k(-X_k,\ldots,-X_1).\plabel{eq:poref}\end{equation}
\begin{proof}
This follows from \eqref{eq:mpform} immediately:
Descents and ascents get inverted in permutations; which counts for $(-1)^{k-1}$ in the sign change.
(It also follows from \eqref{eq:mudef}.)
\end{proof}
\end{lemma}
For the sake of the next lemma, we introduce a formal trace $\Tr$ on the algebra noncommutative
polynomials of $X_1,\ldots,X_k$. This takes any monomial $X_{i_1}\ldots X_{i_s}$ into
a cyclic word $\Tr X_{i_1}\ldots X_{i_s}\equiv\langle X_{i_1}\ldots X_{i_s}\rangle$, and otherwise linear.

\begin{lemma}
\begin{equation}\Tr X_0\mu_k(X_1,\ldots,X_k)=\Tr \mu_k(X_0,\ldots,X_{k-1})X_k.\plabel{eq:pocyc}\end{equation}
\begin{proof}
This also follows from \eqref{eq:mpform}.
One can consider cyclic permutations $\hat\sigma\in\hat\Sigma_{k+1}$, with its ascent-descent patterns
but where we count not $k$ many ascents and descents
but $k+1$ many (due to cyclicity).
Then both sides
yield
\[\sum_{\hat\sigma\in\hat\Sigma_{k+1}}(-1)^{\des(\hat\sigma)-1}
\frac{(\asc(\hat\sigma)-1)!(\des(\hat\sigma)-1)!}{k!}
\langle X_{\hat\sigma(0)} \ldots X_{\hat\sigma(k)}\rangle; \]
except in the first case we imagine the cyclic permutations ordered to start with $0$,
and  in the second case we imagine the cyclic permutations ordered to end with $k$.
(It also follows from \eqref{eq:mudef}.)
\end{proof}
\end{lemma}
We can frame \eqref{eq:poref} and \eqref{eq:pocyc} as symmetries of $\mu_k^{\Lie}(X_1,\ldots,X_k)$ as follows.
Let $\mathfrak M_k$ be the vector space generated by $1$-homogeneous monomials of $X_1,\ldots,X_k$ of length $k$.
We define the operation $^\dag$ on $\mathfrak M_k$ by
\[(X_{i_1}\ldots X_{i_k})^\dag= (-1)^{k+1} X_{k+1-i_k}\ldots X_{k+1-i_1}.\]
Similarly, we define the operation $^\curvearrowright$ on $\mathfrak M_k$ by
\[(X_{i_1}\ldots X_{i_{s-1}}X_{i_s\equiv k}   X_{i_{s+1}}\ldots   X_{i_k})^ \curvearrowright=X_{i_{s+1}+1}\ldots
X_{i_k+1}X_1 X_{i_1+1}\ldots X_{i_{s-1}+1}.\]
Thus $\Tr X_0 M\mapsto \Tr X_0 M^\dag$ represents taking $X_i\mapsto X_{-i}$  (and a possible sign change) on $\Tr X_0\mathfrak M_k$,
and $\Tr X_0 M\mapsto \Tr X_0 M^\curvearrowright$ represents  taking $X_i\mapsto X_{i+1}$ on $\Tr X_0\mathfrak M_k$,
where the indices themselves count $\mod k+1$.
In particular, $^\dag$ and $^\curvearrowright$ define a dihedral action; which is  an action of $D_{2\cdot(k+1)}$.
It is easy to see that \eqref{eq:poref} and \eqref{eq:pocyc} express the invariance of $\mu_k(X_1,\ldots,X_k)$
under $^\dag$ and $^\curvearrowright$ respectively.
Now, we can also consider the vector space $\mathfrak M_k^{\Lie}$ generated by $1$-homogeneous Lie monomials of $X_1,\ldots,X_k$ of length $k$.

\begin{lemma} $^\dag$ and $^\curvearrowright$ descend to
$\mathfrak M_k^{\Lie}$. Moreover, they take Lie monomials into Lie monomials.

\proofremark{We  consider any Lie monomial times $-1$ also as a Lie monomial; but
multiplication by $-1$ can also be incorporated as a change of order.}
\begin{proof} We can consider the Lie monomials as commutator monomials.
The case of $^\dag$ is simple, because it involves only substitutions and a possible sign change.
Consider now $^\curvearrowright$.
Then any Lie monomial of $M$ of $X_1,\ldots,X_k$ can be written in form
\[M=[M_1,[M_2,[\ldots [M_s,X_{k}]\ldots] ,\] where the $M_p$ are other Lie monomials (of $X_1,\ldots,X_{k-1}$ altogether).
Then from the formal properties of trace,
\[\Tr X_0M=\Tr X_0[M_1,[M_2,[\ldots [M_s,X_{k}]\ldots]=\Tr[\ldots [X_{0},M_1]\ldots],M_{s-1}],M_s] X_{k} .\]
This indicates that $M^\curvearrowright$  is just $[\ldots [X_{0},M_1]\ldots],M_{s-1}],M_s] $ but the indices of the $X_i$
raised by $1$.
\end{proof}
\end{lemma}
\begin{cor}
$\mu_k^{\Lie}(X_1,\ldots,X_k)$ is invariant under the Lie-monomially (and dihedrally) acting symmetries $^\dag$ and $^\curvearrowright$.
\qed
\end{cor}
In fact, we can represent Lie monomials as rooted binary planar trees.
In that case $^\dag$ corresponds to reflection and $^\curvearrowright$ corresponds to re-rooting by $1$.
\begin{commentx}
More generally, one can abstract an algebra structure from $\Tr X_0M^{\Lie}(X_1,\ldots, X_k)$,
and one can approach the BCH and Magnus expansions through them.
\end{commentx}

There is a huge literature more or less related to the algebraic properties of $\mu_k$.
Reutenauer \cite{R} contains already a large amount of information in this direction;
and by now there is a considerable literature related to pre-Lie algebras, and other generalizations., cf.
Ebrahimi-Fard, Manchon \cite{EFM}; and several further articles.
\snewpage
\section{Principles of convergence}
\plabel{sec:Lieconv}
The convergence of Magnus expansion in the setting of  Banach--Lie-algebras
is customarily examined in terms of Banach--Lie norms.
In what follows, let $\mathfrak g$ be a Banach--Lie-algebra, i. e. Banach space endowed with a norm-compatible
Lie algebra structure $\|\cdot\|$ such that
\[\|[X,Y]\|\leq \|X\|\cdot \|Y\|\]
holds. This is not exactly compatible with the Banach-algebra settings.
If $\mathfrak A$ is a Banach algebra with norm $|\cdot|$, then it becomes a Banach--Lie algebra with
the norm
\[\|A\|=2|A|.\]
(Indeed, in this case $\|[A,B]\|=2|[A,B]|\leq 4|A|\cdot|B|=\|A\|\cdot\|B\|$.)

If $\phi$ is Banach--Lie algebra valued ordered measure of finite variation, we may consider
\begin{equation}
\mu_{\mathrm R}^{\Lie}(\phi)=\sum_{k=1}^\infty\mu_{k,\mathrm R}^{\Lie}(\phi),
\plabel{eq:MagnusRLie}
\end{equation}
where
\[\mu_{k,\mathrm R}^{\Lie}(\phi)=\int_{t_1\leq\ldots\leq t_k\in I}\mu_k^{\Lie}(\phi(t_1),\ldots,\phi(t_k));\]
etc., in analogy of the associative algebraic case.
In fact, the principal subject of the present paper is the converge properties of the sum \eqref{eq:MagnusRLie}.

The analogy of the Magnus expansion theorem \eqref{eq:MagnusRold}, giving the geometric interpretation of
$\mu_{\mathrm R}^{\Lie}(\phi)$ and the corresponding time-ordered integral should, obviously, be considered on some sort of Banach--Lie groups.
In this paper, we leave the full geometric treatment aside.
But, even in  geometric investigations, the question of convergence comes up naturally, thus this is a subject one should deal with.

Let $\monom_k^{\Lie}$ be the set of all Lie monomials of $X_1,\ldots,X_k$, where every variable is with multiplicity $1$.
Let
\[\Theta_k^{\Lie}:=\frac1{k!}\inf\left\{\sum_{\gamma\in\monom_k^{\Lie}} |\theta_\gamma|\,:
\,\mu^{\Lie}_k(X_1,\ldots,X_k)=\sum_{\gamma\in\monom_k^{\Lie}} \theta_\gamma\cdot \gamma(X_1,\ldots,X_k)
\,,\, \theta_\gamma\in\mathbb R \right\},
\]
that is $1/k!$ times the minimal sum of the absolute value of the coefficients of the presentations of $\mu^{\Lie}_k(X_1,\ldots,X_k)$.
We use the term `minimal presentation' where this sum is minimal; it is an easy compactness argument that such minimal
presentations exist.
We will not consider two Lie monomials different if they represent the same element in the free Lie algebra.
However, if $M$ is a Lie monomial, then we also consider $-M$ as a Lie monomial,
(although this makes difference only for $k=1$, and, in practice, not even then).
From that viewpoint, in minimal presentations, $\theta_\gamma\geq0$ can be assumed.
Furthermore, in minimal presentations only one of $M$ and $-M$ can occur with positive coefficients;
otherwise there would be a cancelation in the average, in contradiction to minimality.

More generally, we can consider the free Lie algebra $\mathrm F^{\Lie}[Y_\lambda:\lambda\in\Lambda]$, and it can be endowed
by the Lie-norm $\|\cdot\|_{\ell^1}$ defined in analogous manner.
(Here and in what follows we suppress the dependence on the base field $\mathbb K$, as, similarly to
the ordinary non-Lie $\ell^1$ case, the complex case will simply be an extension of the real case:
When we consider the coefficient minimization problem to real non-commutative polynomials,
there is no advantage in using complex coefficients as one can always restrict to the real part.)
Then  $\|\cdot\|_{\ell^1}$ decomposes in multidegrees, and there are minimal presentations in every multidegree.
This can be completed to a Banach--Lie norm on the Banach--Lie algebra $\mathrm F^{1,\Lie}[Y_\lambda:\lambda\in\Lambda]$
or to the locally convex  Banach--Lie algebra $\mathrm F^{1,\loc, \Lie}[Y_\lambda:\lambda\in\Lambda]$
in analogy to the ordinary spaces of noncommutative power series.
In this terminology,
$\Theta_k^{\Lie}=\frac{1}{k!}\|\mu^{\Lie}_k(Y_1,\ldots,Y_k)\|_{\ell^1}$.

Computing $\Theta_k^{\Lie}$ can be posed as a straightforward problem in rational linear programming.
(See the basic textbook Matou\v{s}ek, G\"{a}rtner \cite{MG}, and further references therein.)
The problem, however,  quickly grows intractable due to its size.
Nevertheless, a couple of terms can be computed.
Before going to the details, let us make some practical conventions.
We prefer lexicographically minimal presentations  in the order of variables.
E. g., we prefer $-[[X_1,X_3],X_2]$ to $[X_2,[X_1,X_3]]$.
For the sake of compactness, we will use the notation $X_{[[1,3],2]}\equiv [[X_1,X_3],X_2]$, etc.
Elements of the dual $(\mathfrak M_k)^*$ can be given as linear combinations $X_{j_1,\ldots,j_k}^*$,
such that $X_{j_1,\ldots,j_k}^*(X_{i_1}\ldots X_{i_k})=\delta_{i_1,j_1}\cdot\ldots\cdot\delta_{i_k,j_k}$.
The elements of the dual $(\mathfrak M_k^{\Lie})^*$ occur by restriction.
For that purpose, the generating elements $X_{1,j_2,\ldots,j_{k}}^*$ are sufficient.
Indeed, any element of $\mathfrak M_k$  can uniquely be written as a linear combination of
elements $[[\ldots [X_1,X_{i_{2}}]\ldots], X_{i_k}]$, and, on the other hand,
$X_{1,j_2,\ldots,j_{k}}^*( [[\ldots [X_1,X_{i_{2}}]\ldots], X_{i_k}] )=\delta_{i_2,j_2}\cdot\ldots\cdot\delta_{i_{k},j_{k}}$.

If one wants to demonstrate $\Theta_k^{\Lie}\geq\frac1{k!} c$, then it is sufficient to find a linear functional $\alpha\in(\mathfrak M_k^{\Lie})^*$
such that $-1\leq \alpha(M)\leq 1$ for any monomial $M\in\mathbf M_{k}^{\Lie}$ but $\alpha( \mu_k^{\Lie}(X_1,\ldots,X_k))\geq c$.
If one wants to demonstrate $\Theta_k^{\Lie}\leq\frac1{k!} c$, then it is sufficient merely to find a monomial representation
of $\mu^{\Lie}_k(X_1,\ldots,X_k)$ where the sum of absolute value of the coefficients of the monomials is at most $c$.

In more details, the convex span of the set $\mathbf M_{k}^{\Lie}$ of Lie monomials forms the unit ball $\mathcal I_k^{\Lie}$
of $\|\cdot\|_{\ell^1}$ in $\mathfrak M_{k}^{\Lie}$.
Its extremal points are exactly the elements of $\mathbf M_{k}^{\Lie}$.
(This is easy to show using an appropriate euclidean norm on $\mathfrak M_k$).
If $[0,\infty)\cdot\mu^{\Lie}_k(X_1,\ldots,X_k) $ intersects $\partial \mathcal I_k^{\Lie}$ in the interior of a
face of dimension $d_k$ and of number of vertices $v_k$, then the minimal presentations form an ``abstract'' polytope
of dimension $\Delta_k=v_k-d_k-1$.
If $\mu_k^{\Lie}(X_1,\ldots,X_k) $ is obtained from $\mathbf M_{k}^{\Lie}$ as convex combination in a quasi-canonical manner,
depending only the linear structure of  $\mathbf M_{k}^{\Lie}$, then its description must be invariant for the
symmetries $^\dag$ and $^\curvearrowright$, as they are linear automorphisms of $\mathcal I_k^{\Lie}$.
In particular, if the convex combination $\mu_k^{\Lie}(X_1,\ldots,X_k)$ is obtained from $\mathbf M_{k}^{\Lie}$
through a quasi-barycentric decomposition, then it is symmetric for $^\dag$ and $^\curvearrowright$.
In general, looking for symmetric decompositions reduces the size of the problem essentially by a factor of $k$,
although this puts us only one number ahead, as the size is about $k!$.
\begin{example}
\[\Theta_4^{\Lie}=\frac1{4!}\cdot\frac13. \]
Indeed, the linear functional
\begin{equation}\alpha_4={X}^*_{1234}-{X}^*_{1342}-{X}^*_{1423}-{X}^*_{1432}.
\plabel{eq:mu4dual}
\end{equation}
takes values on $\mathbf M_4^{\Lie}$ only from $\{-1,0,1\}$, yet it takes $\frac1{3}$ on $\mu_4(X_1,X_2,X_3,X_4)$;
proving $\Theta_4^{\Lie}\geq\frac1{4!}\cdot \frac1{3}$.
On the other hand, we have equality according to
\begin{equation}
\mu_4^{\Lie}(X_1,X_2,X_3,X_4)=\frac1{12}\Bigl(
X_{{[[1,3],[2,4]]}}+X_{{[[1,[2,3]],4]}}+X_{{[[1,2],[3,4]]}}+X_{{[1,[[2,3],4]]}}
\Bigr);
\plabel{eq:mu4}
\end{equation}
(the sum of the absolute values of the coefficients being $\frac13$),
which, therefore, yields a minimal presentation.
Taking a closer look, we find that the minimal presentations form a $\Delta_4=4$ dimensional simplex.
In more concrete terms, the minimal presentations are of shape
\begin{align}
\mu_4^{\Lie}(X_1,X_2,X_3,X_4)=\frac1{12}\Bigl(
&-X_{{[[[1,4],2],3]}}\lambda_1+X_{{[1,[2,[3,4]]]}} \left(\lambda_5+\lambda_1+\lambda_2 \right)\notag\\
&+X_{{[[1,[2,4]],3]}}\lambda_2+X_{{[[1,[2,3]],4]}} \left(\lambda_1+\lambda_2+\lambda_3 \right)\notag\\
&+X_{{[[1,3],[2,4]]}}\lambda_3+X_{{[[1,2],[3,4]]}} \left(\lambda_2+\lambda_3+\lambda_4 \right)\notag\\
&-X_{{[[[1,3],4],2]}}\lambda_4+X_{{[1,[[2,3],4]]}} \left(\lambda_3+\lambda_4+\lambda_5 \right)\notag\\
&-X_{{[[[1,4],3],2]}}\lambda_5+X_{{[[[1,2],3],4]}} \left(\lambda_4+\lambda_5+\lambda_1 \right)\notag
\Bigr),
\end{align}
where
\[\lambda_i\geq0,\qquad \lambda_1+\lambda_2+\lambda_3+\lambda_4+\lambda_5=1.\]

The presentation which is most economical and also symmetrical to $^\dag$, belongs to
$(\lambda_1,\lambda_2,\lambda_3,\lambda_4,\lambda_5)=(0,0,1,0,0)$.
It yields exactly \eqref{eq:mu4}.
However, there is only one minimal presentation which is symmetrical to $^\curvearrowright$ (and also symmetrical to $^\dag$),
which belongs to  $(\lambda_1,\lambda_2,\lambda_3,\lambda_4,\lambda_5)=(\frac15,\frac15,\frac15,\frac15,\frac15)$,
and yields
\begin{align}
\mu_4^{\Lie}(X_1,X_2,X_3,X_4)=\frac1{60}\Bigl(
&-X_{{[[[1,4],2],3]}}+3X_{{[1,[2,[3,4]]]}} \notag\\
&+X_{{[[1,[2,4]],3]}}+3X_{{[[1,[2,3]],4]}} \notag\\
&+X_{{[[1,3],[2,4]]}}+3X_{{[[1,2],[3,4]]}} \notag\\
&-X_{{[[[1,3],4],2]}}+3X_{{[1,[[2,3],4]]}} \notag\\
&-X_{{[[[1,4],3],2]}}+3X_{{[[[1,2],3],4]}} \notag
\Bigr).
\end{align}
This must also be any quasi barycentric decomposition;  notice, in particular, that all vertices of the critical face appear.
\qedexer\end{example}
From the viewpoint of linear programming, we have dealt with an optimization (minimalization) problem for
computing $4!\Theta_4^{\Lie}$;
the minimal presentation \eqref{eq:mu4} is called a primal solution of the LP problem;
and the good face-defining linear functional \eqref{eq:mu4dual}
is called a dual solution of the LP problem.
The machinery of linear programming can be applied in order to compute
other terms of $\Theta^{\Lie}(x)$:
\begin{cor}
\[\Theta_1^{\Lie}=1,\qquad \Theta_2^{\Lie}=\frac1{2!}\cdot\frac12,\qquad \Theta_3^{\Lie}=\frac1{3!}\cdot\frac13, \]
\[\Theta_4^{\Lie}=\frac1{4!}\cdot\frac13, \qquad \Theta_5^{\Lie}=\frac1{5!}\cdot\frac25,\qquad \Theta_6^{\Lie}=\frac1{6!}\cdot\frac{37}{60}, \]

\[\Theta_7^{\Lie}=\frac1{7!}\cdot\frac{1621}{1540},\qquad
\Theta_8^{\Lie}=\frac1{8!}\cdot
\frac{5242130799984621832318}{2419342933460499216625}.\]
The value of $\Theta_9^{\Lie}$ is also known.
\begin{proof}
Cf. Appendix \ref{sec:someTheta}.
\end{proof}
\end{cor}

We can give an estimate for the convergence of the Magnus expansion as follows.
We define the absolute Magnus--Lie characteristic $\Theta^{\Lie}$ by
\[\Theta^{\Lie}(x)=\sum_{k=1}^\infty\Theta_k^{\Lie}x^k.\]

Then
\[\left|\mu_{k,\mathrm R}^{\Lie}(\phi)\right|\leq\int_{t_1\leq\ldots\leq t_k\in I}\left|\mu_k^{\Lie}(\phi(t_1),\ldots,\phi(t_k))\right|\leq
\Theta_k^{\Lie}\cdot\left(\int|\phi| \right)^k,\]
and, consequently,
\begin{equation}
\sum_{k=1}^\infty\left|\mu_{k,\mathrm R}^{\Lie}(\phi)\right|\leq\sum_{k=1}^\infty\int_{t_1\leq\ldots\leq t_k\in I}\left|\mu_k^{\Lie}(\phi(t_1),\ldots,\phi(t_k))\right|\leq\Theta_\real^{\Lie}\left(\int|\phi|\right).
\plabel{eq:measumLie}
\end{equation}
Note that we have equalities in the case $\phi=c\cdot\mathrm Z^{1,\Lie}_{[0,1]}$, $c\in[0,\infty)$,
the ordered totally noncommutative continuous Lie mass of norm $c$.
Indeed, as $\mathrm F^{1}[Y_\lambda:\lambda\in\Lambda]$ generalizes to $\mathrm F^{1}([0,1])$,
$\mathrm F^{1,\Lie}[Y_\lambda:\lambda\in\Lambda]$ generalizes to $\mathrm F^{1,\Lie}([0,1])$
in the same manner; and we also have the tautological measure.
Thus, similarly to the associative algebraic case, the question of convergence
in the most general case reduces to the study of $\Theta_{\real}^{\Lie}(x)$.
We have convergence if $\Theta_{\real}^{\Lie}(\smallint |\phi|)<+\infty$,
and we can have divergence of $\Theta_{\real}^{\Lie}(\smallint |\phi|)=+\infty$ (using the tautological meausure).
We have to find the convergence radius of $\Theta^{\Lie}(x)$.

Regarding the convergence radius $\mathrm C_{\infty}^{\Lie}$ of $\Theta^{\Lie}$, we have the following easy estimates.
From the commutator expansion, it follows that
\[\Theta_n\leq2^{n-1}\Theta_n^{\Lie},\]
hence
\[2\,\Theta\left(\frac x2\right)\stackrel\forall\leq\Theta^{\Lie}(x),\]
consequently,
\[2\,\Theta_{\real}\left(\frac x2\right)\leq\Theta^{\Lie}_{\real}(x).\]
In particular, the convergence radius of $\Theta^{\Lie}$ is at most $4$.
On the other hand, the Dynkin--Specht--Wever lemma implies
\[\Theta_n^{\Lie}\leq\frac1n \Theta_n, \]
thus
\[\Theta^{\Lie}(x)\stackrel\forall\leq\int_{t=0}^1\frac{\Theta(xt)}t\,\mathrm dt.\]
Consequently,
\[\Theta^{\Lie}_{\real}(x)\leq\int_{t=0}^1\frac{\Theta_{\real}(xt)}t\,\mathrm dt=\int_{t=0}^x\frac{\Theta_{\real}(t)}t\,\mathrm dt.\]
In particular, the convergence radius of $\Theta^{\Lie}$ is at least $2$.
Thus $2\leq \mathrm C_{\infty}^{\Lie}\leq 4$.

Here the lower estimate is well known  not to be sharp (cf. the Introduction):
According to the standard estimate in the literature, we know that
\[\boldsymbol\delta=\int_{x=0}^{2\pi} \frac{\mathrm dx}{2+\frac x2-\frac x2\cot \frac x2 }= 2.1737374\ldots\]
provides a better lower estimate; $ \boldsymbol\delta\leq \mathrm C_{\infty}^{\Lie}\leq 4$.
In this paper, we will improve this lower bound further.

We can approach the convergence of Baker--Campbell--Hausdorff expansion in the same manner.
Indeed, $\Gamma^{\Lie}(x,y)$ can be defined analogously to $\Gamma(x,y)$, and
crude applications of the commutator expansion and the Dynkin--Specht--Wever lemma yield
\begin{equation}2\,\Gamma\left(\frac x2,\frac y2\right)\stackrel\forall\leq\Gamma^{\Lie}(x,y) \plabel{eq:gamup}\end{equation}
and
\begin{equation}\Gamma^{\Lie}(x,y)\stackrel\forall\leq \int_{t=0}^1\frac{\Gamma(xt,yt)}t\,\mathrm dt.\plabel{eq:gamdown}\end{equation}

We define the cumulative convergence radius $\mathrm C^{\Lie}_2$ as
\[\mathrm C^{\Lie}_2=\sup\left\{ x\in[0,+\infty)\,:\,\max_{0\leq x_1,x_2,x_1+x_2\leq x} \Gamma^{\Lie}(x_1,x_2) <+\infty  \right\};\]
and the diagonal convergence radius
\[\mathrm D^{\Lie}_2=\sup\left\{ x\in[0,+\infty)\,:\, \Gamma^{\Lie}\left(\frac x2,\frac x2\right) <+\infty  \right\},\]
i. e. as the convergence radius of the series $\Gamma^{\Lie}\left(\frac x2,\frac x2\right)$.
Obviously,
$\mathrm C^{\Lie}_2\leq\mathrm D^{\Lie}_2$.
Thus, according to the results in Part I and the previous argument, one has for the convergence radii of the
BCH expansion
\[\mathrm C_2\leq\mathrm C_2^{\Lie}\leq\mathrm D^{\Lie}_2\leq2\mathrm C_2,\]
where $\mathrm C_2=2.89847930\ldots$.
These upper and lower estimates are better than the ones one can read in the literature, but
we will (slightly) improve them later.

As we have indicated, we do not embark on a full geometrical treatment of the Magnus exponential
theorem
\begin{equation}
\sum_{k=1}^\infty\left\|\mu_{k,\mathrm R}(\phi)\right\|<+\infty\qquad\Rightarrow\qquad \exp_{\mathrm R}(\phi)=\exp \mu_{\mathrm R}(\phi).
\plabel{eq:MagnusRold}
\end{equation}
However, there are some weaker substitutes, which can be useful.
The simplest one is
\begin{theorem}\plabel{th:Magnuscor}
Suppose that $\phi$ is a $\mathfrak g$-valued measure, $Y\in\mathfrak g$.
If
\begin{equation}
\sum_{k=1}^\infty\left\|\mu_{\mathrm L[k]}^{\Lie}(\phi)\right\|<+\infty,
\plabel{eq:sersumcor}
\end{equation}
then
\begin{equation}
\exp(\mu_{\mathrm L}^{\Lie}\left(\phi)\right)Y=\Lexp(\ad\phi)Y.
\plabel{eq:magcor}
\end{equation}
\begin{proof} This a corollary of the Magnus exponential theorem applied to the adjoint representation.
(We have changed to `$\mathrm L$' formalism to make it more analytic looking.)
\end{proof}
\end{theorem}
A more sophisticated (and geometric) version is
\begin{theorem}\plabel{th:loccontract}
Suppose that $\phi_1,\phi_2,\phi_3$ are a $\mathfrak g$-valued measures.
If
\begin{equation}
\Theta^{\Lie}( \smallint\|\phi_1\|+\Theta^{\Lie}( \smallint\|\phi_2\|) + \smallint\|\phi_3\|) <+\infty ,
\end{equation}
then
\begin{equation}
\mu_{\mathrm R}(\phi_1\boldsymbol. \phi_2  \boldsymbol.\phi_2 )
=\mu_{\mathrm R}(\phi_1\boldsymbol. \mu_{\mathrm R}(\phi_2)\mathbf 1_{[0,1]}   \boldsymbol.\phi_3 ).
\end{equation}
\begin{proof} This is Corollary \ref{cor:genMagnusrec} integrated and contracted.
\end{proof}
\end{theorem}
Similar statements can be devised for
$\mu_{\mathrm R}(\phi_1\boldsymbol. \mu_{\mathrm R}(\phi_2)\mathbf 1_{[0,1]}   \boldsymbol.\phi_2
\boldsymbol. \mu_{\mathrm R}(\phi_4)\mathbf 1_{[0,1]}   \boldsymbol.\phi_5 )$ , etc.

Consider now the variant with $\phi=c\cdot\mathrm Z^{1,\loc,\Lie}_{[0,1)}$
over the locally convex Banach--Lie algebra $\mathrm F^{1,\loc,\Lie}([0,1))$.
Here convergence is not a problem (as we have a limit of nilpotent algebras).
Thus, for $c\geq0$, we obtain
\begin{equation}
\mu_{\mathrm R}(c\cdot\mathrm Z^{1,\loc,\Lie}_{[0,1)})
=\mu_{\mathrm R}\Bigl(\underbrace{\mu_{\mathrm R}(c\cdot\mathrm Z^{1,\loc,\Lie}_{[\omega_0,\omega_1)})\mathbf 1_{[\omega_0,\omega_1]}\boldsymbol.\ldots\boldsymbol.
\mu_{\mathrm R}(c\cdot\mathrm Z^{1,\loc,\Lie}_{[\omega_{s-1},\omega_s)})\mathbf 1_{[\omega_{s-1},\omega_s]}}_{ \mathrm{Red}_\omega(c\cdot\mathrm Z^{1,\loc,\Lie}_{[0,1)})} \Bigr).
\plabel{eq:somin}
\end{equation}
for any division $0=\omega_0<\omega_1<\ldots<\omega_s=1$.
However, for small $x$, $\Theta^{\Lie}(x)\sim x$, thus, for a sufficiently fine division,
$\|\mathrm{Red}_\omega(c\cdot\mathrm Z^{1,\loc,\Lie}_{[0,1)})\|_{\ell^1}=c+\varepsilon$, where $\varepsilon$ is arbitrarily small.
Then, due to equality \eqref{eq:somin} which holds in $\mathrm F^{1,\loc,\Lie}([0,1))$,
 if  $\Theta^{\Lie}(c)=+\infty$
(i. e. it yields the divergence of the Magnus series restricted to $\mathrm F^{1,\Lie}([0,1))$ ), then
$\phi=\mathrm{Red}_\omega(c\cdot\mathrm Z^{1,\loc,\Lie}_{[0,1)})$ gives an
example for divergence with $\smallint |\phi|=c+\varepsilon$ but $\phi$ being of multivariable BCH type.
This shows that we cannot expect a better convergence radius using smoother integrands instead of measures.
\snewpage
\section{Divergence estimate for the BCH expansion}
\plabel{sec:bchupper}
\begin{example}\plabel{ex:bchunwrap}
Let us consider the noncommutative polynomial algebra generated by $X,Y$
but factorized by the relations $X^2=0$, $Y^2=0$.
Consider its Lie subalgebra $\mathfrak{wsl}_2^{\mathrm{alg}}$
generated by $X$ and $Y$.
Linearly, it is generated by the elements
\begin{align}
X^{(2n+1)}&=(XY)^{n}X&&(n\geq0),\notag\\
Y^{(2m+1)}&=(YX)^{m}Y&&(m\geq0),\notag\\
XY^{(2k)}&=(XY)^{k}-(YX)^k &&(k\geq1).\notag
\end{align}

The commutation rules are given by
\begin{align}
[X^{(2n+1)},Y^{(2m+1)}]&=XY^{(2n+2m+2)},&[X^{(2n+1)},X^{(2\tilde n+1)}]&=0,\notag\\
[XY^{(2k)},X^{(2n+1)}]&=2X^{(2k+2n+1)},&[Y^{(2m+1)},Y^{(2\tilde m+1)}]&=0,\notag\\
[Y^{(2m+1)},XY^{(2k)}]&=2Y^{(2k+2m+1)},&[XY^{(2k)},XY^{(2\tilde k)}]&=0.\notag
\end{align}

We will define a norm $\|\cdot\|$ on it such that
$\| X^{(2n+1)}\|=2^{-n}$ $\| Y^{(2m+1)}\|=2^{-m}$,  $\| XY^{(2k)}\|=2^{-k+1}$
Indeed, we simply let
\begin{multline}
\left\|\sum_{n=0}^\infty b_nX^{(2n+1)}+\sum_{m=0}^\infty c_mY^{(2m+1)}+\sum_{k=1}^\infty d_mY^{(2k)}\right\|=\\=
\sum_{n=0}^\infty |b_n|2^{-n}+\sum_{m=0}^\infty |c_m|2^{(-m)}+\sum_{k=1}^\infty |d_m|2^{-k+1}.\notag
\end{multline}
Using the commutation rules, one can check that $\mathfrak{wsl}_2^{\mathrm{alg}}$ becomes a normed Lie algebra with $\|\cdot\|$;
which can be completed to  $\mathfrak{wsl}_2^{1}$ as usual.

Let us compute
$\BCH(tX,tY)=\log(\exp(tX)\exp(tY))$.
This can  done formally, because imposing  $X^2=0$, $Y^2=0$ is compatible with the grade filtration.
If $X,Y$ are formal variables, and $\mathcal R^{(\lambda)}(A)=\frac{A-1}{\lambda+(1-\lambda)A}$, then
\begin{align}
\mathcal R^{(\lambda)}((\exp X)(\exp Y))=\sum_{k=0}^\infty\biggl(
&(\lambda-1)^k\lambda^k\mathcal R^{(\lambda)}(\exp Y)\left(\mathcal R^{(\lambda)}(\exp X)\mathcal R^{(\lambda)}(\exp Y)\right)^{k}\notag\\
+&(\lambda-1)^k\lambda^k\mathcal R^{(\lambda)}(\exp X)\left(\mathcal R^{(\lambda)}(\exp Y)\mathcal R^{(\lambda)}(\exp X)\right)^{k}\notag\\
+&(\lambda-1)^{k}\lambda^{k+1}\left(\mathcal R^{(\lambda)}(\exp X)\mathcal R^{(\lambda)}(\exp Y)\right)^{k+1} \notag\\
+&(\lambda-1)^{k+1}\lambda^{k}\left(\mathcal R^{(\lambda)}(\exp Y)\mathcal R^{(\lambda)}(\exp X)\right)^{k+1} \qquad\,\,\biggr).\notag
 \end{align}
(Cf. Part I., \cite{L1}.) But in our case, $\mathcal R^{(\lambda)}(\exp X)=X$, $\mathcal R^{(\lambda)}(\exp Y)=Y$, thus
\begin{align*}
\mathcal R^{(\lambda)}((\exp X)(\exp Y))=\sum_{k=0}^\infty\biggl(
&(\lambda-1)^k\lambda^k(YX)^kY
+(\lambda-1)^k\lambda^k(XY)^kX\\
&+(\lambda-1)^{k}\lambda^{k+1}(XY)^{k+1}
+(\lambda-1)^{k+1}\lambda^{k}(YX)^{k+1} \biggr).
 \end{align*}
Integrated,
we find
\[\BCH(X,Y)=
\sum_{k=0}^\infty\frac{\left( -1 \right) ^{k}(k!)^2}{(2k+1)!}
\biggl(  X^{(2k+1)}+Y^{(2k+1)}+\frac12XY^{(2k+2)} \biggr),\]
formally. Scaled appropriately,
\[\BCH(tX,tY)=
\sum_{k=0}^\infty
\frac{\left( -1 \right) ^{k}(k!)^2}{(2k+1)!}
\biggl(  t^{2k+1}X^{(2k+1)}+t^{2k+1}Y^{(2k+1)}+t^{2k+2}\frac12XY^{(2k+2)} \biggr).\]

In particular, we obtain the estimate
\[\Gamma^{\Lie}(t,t)\stackrel\forall\geq \Xi(t)\equiv
\sum_{k=0}^\infty\frac{(k!)^2}{(2k+1)!}
\biggl(  t^{2k+1}2^{1-k}+t^{2k+2}2^{-k-1} \biggr).\]
Now
\[\Xi_\real(t)=
\begin{cases}
\dfrac{2\arcsin(2^{-3/2}t)}{\sqrt{8-t^2}}(4+t)&\text{if}\quad 0\leq t<2^{3/2}\\
+\infty&\text{if}\quad 2^{3/2}\leq t<+\infty.
\end{cases}\]
In particular, $\BCH(2\sqrt2 X, 2\sqrt2 Y)$ does not converge (nor $\log(\exp(2\sqrt2 X)\exp(2\sqrt2 Y))$ exists)
in $\mathfrak{wsl}_2^{1}$.
\qedexer
\end{example}
Before continuing, let us discuss the Banach--Lie Frobenius norm on $\mathfrak{sl}_2(\mathbb R)$.
For $\mathbf v=(a,b,c)$ let
\[S(\mathbf v)=\begin{bmatrix}a&\sqrt2b\\\sqrt2c&-a\end{bmatrix}.\]
We also set $\|S(\mathbf v)\|_{\mathrm F}=2|\mathbf v|$.
This is the Frobenius norm upscaled by $\sqrt2$:
One can see that
$\left\|
\left[\begin{smallmatrix}a&b\\c&-a\end{smallmatrix}\right]\right\|_{\mathrm F}=\sqrt2\sqrt{2a^2+b^2+c^2}$.
This yields a Banach--Lie algebra, because
\begin{multline}
\|[S(\mathbf v_1),S(\mathbf v_2) ]\|_{\mathrm F}
=\|2 S(\mathbf v_1\widetilde\times\mathbf v_2)\|_{\mathrm F}
=4|\mathbf v_1\widetilde\times\mathbf v_2 |=
4|\mathbf v_1\times\mathbf v_2 |\leq\\\leq 4|\mathbf v_1|\,|\mathbf v_2|
=\|S(\mathbf v_1)\|_{\mathrm F}\,\|S(\mathbf v_2)\|_{\mathrm F},\notag
\end{multline}
where $\langle a_1,b_1,c_1\rangle\widetilde\times\langle a_1,b_1,c_1\rangle=\langle b_1c_2-c_1b_2, a_1b_2-b_1a_2, c_1a_2-a_1c_2\rangle$.
(This Banach--Lie Frobenius norm was already considered by Michel \cite{Mi}.)

\begin{example}\plabel{ex:bchwrap}
The previous example can also be wrapped up as a counterexample on $\mathfrak{sl}_2$.
The idea is to take the homomorphism induced by
\[X\mapsto \hat X=\begin{bmatrix}0&\sqrt2/2\\&0\end{bmatrix},\qquad Y\mapsto \hat Y=\begin{bmatrix}0&\\\sqrt2/2&0\end{bmatrix}.\]
We will consider this with the Banach--Lie Frobenius norm.
Then $\|\hat X\|_{\mathrm F}=\|\hat Y\|_{\mathrm F}=1$, and
\begin{align}
\BCH(t\hat X,t\hat Y)=&
\sum_{k=0}^\infty\frac{\left( -1 \right) ^{k}(k!)^2}{(2k+1)!}
\biggl(t^{2k+1}2^{-k}\begin{bmatrix}0&\sqrt2/2\\\sqrt2/2&0\end{bmatrix}+t^{2k+2}2^{-k} \begin{bmatrix}1/4&0\\0&-1/4\end{bmatrix}\biggr)\notag\\
=&\frac{2t\arsinh(2^{-3/2}t)}{\sqrt{8+t^2}}
\begin{bmatrix}t&2\sqrt2\\2\sqrt2&-t\end{bmatrix}
=\AC\left(1+\frac14t^2\right)\frac t4
\begin{bmatrix}t&2\sqrt2\\2\sqrt2&-t\end{bmatrix}\notag
\end{align}
formally. (Cf. Part II, \cite{L2} for the notation `$\AC$'.)
Here we see that the BCH expansion diverges for $t=2\sqrt2$, although the actual logarithm analytically extends for $t\sim 2\sqrt2$
(but not for $t\sim 2\sqrt2\mathrm i$).
So, in a certain sense, this version is slightly less pronounced.

In any case, we see that the BCH expansion may fail to converge in a Banach--Lie algebra if
$\|X\|_{\mathrm F}=\|Y\|_{\mathrm F}=2\sqrt2=2.82842712\ldots$.
\qedexer
\end{example}

The known best counterexample in this very setting, however, stems from the work of Michel \cite{Mi}.
He considers the general unbalanced case of BCH expansions for divergence in  $\mathfrak{sl}_2(\mathbb R)$ with the
Banach--Lie Frobenius norm (what we have just also considered).
He obtains a  parametrization of some critical cases.
Although he does not give, exact numerical values, Michel \cite{Mi} contains a good graphical representation of the results.
However, in the balanced case, the critical case
(which is also the most favorable in terms of the cumulative norm) can also be obtained a relatively simple manner;
 from the ``rigid balanced''
counterexamples which were already featured in Part II, \cite{L2} but with respect to the operator norm in Hilbert spaces.
\begin{example}\plabel{ex:michelbalance}
Let $\delta\in(0,1)$, and let
\[
A_1^{[\delta]}=\frac{\arccos\delta}{\sqrt{1-\delta^2}}\begin{bmatrix}\delta&-1\\1&-\delta\end{bmatrix},
\qquad\text{and}\qquad
A_2^{[\delta]}=\frac{\arccos\delta}{\sqrt{1-\delta^2}}\begin{bmatrix}-\delta&-1\\1&\delta\end{bmatrix}.
\]
Then $\left(\exp A_1^{[\delta]}\right)\left(\exp A_2^{[\delta]}\right)=
\left[\begin{matrix} 2\delta&-1\\1&\end{matrix}\right]\left[ \begin{matrix} &-1\\1&2\delta\end{matrix}\right]=
\left[ \begin{matrix} -1&-4\delta\\&-1\end{matrix}\right]$,
which is not an exponential of any element of $\mathfrak{sl}_2(\mathbb R)$, therefore its Magnus expansion is completely divergent.
In the Frobenius Banach--Lie norm, however
\[\left\|A_1^{[\delta]}\right\|_{\mathrm F}=\left\|A_2^{[\delta]}\right\|_{\mathrm F}
=\frac{\arccos\delta}{\sqrt{1-\delta^2}}\cdot
2\sqrt{1+\delta^2}.\]
This is minimized at
$\delta_{\mathrm{Mi}}=0.55088194\ldots$ (the solution of $\frac{2\delta\arccos\delta}{(1+\delta^2)\sqrt{1-\delta^2}}=1$ for $\delta\in(0,1)$)
with value $\boldsymbol v_{\mathrm{Mi}}=
\dfrac{\arccos\delta_{\mathrm{Mi}}}{\sqrt{1-\delta_{\mathrm{Mi}}^2}}
\cdot2\sqrt{1+\delta_{\mathrm{Mi}}^2}=2.7014285\ldots$.
Despite not given by Michel \cite{Mi} explicitly, we call $\boldsymbol v_{\mathrm{Mi}}$ as Michel's number.
(Curiously, $\boldsymbol v_{\mathrm{Mi}}$ is the very same number as $\widetilde{\mathrm C}_2$, which comes up as a suboptimal
lower estimate for the convergence radius of the BCH expansion  in Part I, \cite{L1}.)
\qedexer
\end{example}
\plabel{cor:michel}
\begin{cor}
\[\Gamma^{\Lie}_\real(\boldsymbol v_{\mathrm{Mi}},\boldsymbol v_{\mathrm{Mi}})=+\infty.\]
In particular,
\[\mathrm C^{\Lie}_2\leq\mathrm D^{\Lie}_2\leq2\boldsymbol v_{\mathrm{Mi}}=5.4028570\ldots. \eqed\]
\end{cor}

\textbf{Michel's results.}
For the sake of completeness, we also recount Michel's results but in a slightly more user-friendly parametrization.
First, we claim that the function
\[\Mip: u\in(0,\pi)\mapsto \frac{1-\cos u}{1+\cos u}\sqrt{\frac{1+\frac{\sin u}u}{1-\frac{\sin u}u}}\in(0,+\infty)\]
gives a strictly monotone increasing bijection between $(0,\pi)$ and $(0,+\infty)$.
Therefore, the solution of the equation
\begin{equation}
\Mip(u_1)\Mip(u_2)=1
\plabel{eq:eMip}
\end{equation}
makes a strictly monotone decreasing bijection between $u_1\in(0,\pi)$ and $u_2\in(0,\pi)$.
Now, we set
\[X_1^{[u_1]}=\frac{u_1}{\sqrt[4]{1-\left(\frac{\sin u_1}{u_1}\right)^2}}
\begin{bmatrix}
\sqrt{\frac{1-\sqrt{1-\left(\frac{\sin u_1}{u_1}\right)^2}}{2}}&-\sqrt{\frac{1+\sqrt{1-\left(\frac{\sin u_1}{u_1}\right)^2}}{2}}
\\
\sqrt{\frac{1+\sqrt{1-\left(\frac{\sin u_1}{u_1}\right)^2}}{2}}&-\sqrt{\frac{1-\sqrt{1-\left(\frac{\sin u_1}{u_1}\right)^2}}{2}}
\end{bmatrix}
\]
and
\[X_2^{[u_2]}=\frac{u_2}{\sqrt[4]{1-\left(\frac{\sin u_2}{u_2}\right)^2}}
\begin{bmatrix}
-\sqrt{\frac{1-\sqrt{1-\left(\frac{\sin u_2}{u_2}\right)^2}}{2}}&-\sqrt{\frac{1+\sqrt{1-\left(\frac{\sin u_2}{u_2}\right)^2}}{2}}
\\
\sqrt{\frac{1+\sqrt{1-\left(\frac{\sin u_2}{u_2}\right)^2}}{2}}&\sqrt{\frac{1-\sqrt{1-\left(\frac{\sin u_2}{u_2}\right)^2}}{2}}
\end{bmatrix}.
\]
(Notice the sign differences between $X_1^{[u_1]}$ and $X_2^{[u_2]}$.)

Then, we claim, under \eqref{eq:eMip}, $\left(\exp X_1^{[u_1]}\right)\left(\exp X_2^{[u_2]}\right)$
has Jordan form  $\begin{bmatrix}-1&1\\&-1\end{bmatrix}$.
In particular, the BCH expansion of the pair $X_1^{[u_1]}, X_2^{[u_2]} $ is divergent.
Here
\[\|X_1^{[u_1]}\|_{\mathrm F}=\frac{2u_1}{\sqrt[4]{1-\left(\frac{\sin u_1}{u_1}\right)^2}}\qquad
\text{and}
\qquad
\|X_2^{[u_2]}\|_{\mathrm F}=\frac{2u_2}{\sqrt[4]{1-\left(\frac{\sin u_2}{u_2}\right)^2}}.\]
(Under $u_i\in(0,\pi)$, these norms are in $(0,2\pi)$.)

What is notable in this situation (under \eqref{eq:eMip}) that it is optimal with respect to the norms
and the non-exponential property of the exponential product.
More precisely, the following is true:

Suppose that $A_1,A_2\in\mathfrak{sl}_2(\mathbb R)$ such that
$\|A_1\|_{\mathrm F}\leq \|X_1^{[u_1]}\|_{\mathrm F}$ and $\|A_2\|_{\mathrm F}\leq \|X_2^{[u_2]}\|_{\mathrm F}$
but $\|A_1\|_{\mathrm F}+\|A_2\|_{\mathrm F}<\|X_1^{[u_1]}\|_{\mathrm F}+\|X_2^{[u_2]}\|_{\mathrm F}$.
Then $(\exp A_1)(\exp A_2)$ is an exponential of an element of $\mathfrak{sl}_2(\mathbb R)$.
(Furthermore, the example is unique up to orthogonal conjugation of the pair.)
\begin{proof}[Indication of proof]
The reader is advised to look up Michel \cite{Mi}, where the basic ideas can be followed easily,
(except there is a sudden change of variables $\theta^{(\prime)}\rightsquigarrow\frac\pi4-\theta^{(\prime)}$
just before the final formula).
We also made changes of variables for perceived simplifications; in this respect  we note $\sqrt{\Mip(\sqrt{\mu})}
=\left(\tan\frac{\sqrt\mu}{2}\right)\cdot\sqrt[4]{\frac{1+\frac{\sin\sqrt\mu}{\sqrt\mu}}{1-\frac{\sin\sqrt\mu}{\sqrt\mu}}}$.
\end{proof}
Here $u_1=u_2=u_{\mathrm{Mi}}=0.9749107\ldots$ yields
$\|X_1^{[u_1]}\|_{\mathrm F}=\|X_2^{[u_2]}\|_{\mathrm F}=\boldsymbol v_{\mathrm{Mi}}=2.7014285\ldots$.
That this $ \boldsymbol v_{\mathrm{Mi}}$ is the same as the one encountered before
(and $\delta_{\mathrm{Mi}}=
\frac{1-\sqrt{1-\left(\frac{\sin u_{\mathrm{Mi}}}{ u_{\mathrm{Mi}}}\right)^2}}{ \frac{\sin u_{\mathrm{Mi}}}{ u_{\mathrm{Mi}}} }$%
)
follows from the fact
in the balanced case the shape of the matrices in question is conveniently restricted, and the various optimization
approaches must yield the same result.

\begin{remark}
Returning to the case of the Magnus expansion, the Moan--Sch\"affer example
understood in $\mathfrak{sl}_2(\mathbb R)$
with the Banach--Lie Frobenius norm yields cumulative norm $\sqrt{2}\pi=4.4428829\ldots$
(cf.~Moan \cite{Mo1}, Sch\"affer \cite{Sch}, Part II, \cite{L2}).
In general, using other finite dimensional Banach--Lie algebras, one can get
counterexamples to the convergence of the Magnus expansion with
cumulative norm $4+\varepsilon$ (cf. Part IE, \cite{L11}).
\qedremark
\end{remark}
\snewpage
\section{The standard estimate and its improvements}\plabel{sec:MagnusStandard}
The Magnus recursion formulas imply
\[k!\Theta^{\Lie}_{k}\leq\sum_{l_1+\ldots+l_s=k-1} |\beta_s|\frac{(k-1)!}{l_1!
\cdot\ldots\cdot l_s!}(l_1!\Theta^{\Lie}_{l_1})\cdot\ldots\cdot(l_s!\Theta^{\Lie}_{l_s}) ;\]
which can be rewritten as
\[k \Theta^{\Lie}_{k}\leq\sum_{l_1+\ldots+l_s=k-1} |\beta_s| \Theta^{\Lie}_{l_1}\cdot\ldots\cdot \Theta^{\Lie}_{l_s}. \]
This can be turned into an estimate as follows.

Let us define the numbers $\psi_k$ $(k\in\mathbb N)$ by the recursion
\[\psi_0=0\]
and
\[k \psi_{k}=\sum_{l_1+\ldots+l_s=k-1} |\beta_s| \psi_{l_1}\cdot\ldots\cdot \psi_{l_s}. \]
We also consider the formal generating function
\[\tilde \beta(x):=\sum_{j=0}^\infty |\beta_j| x^j.\]

Then it is immediate that
\[\Theta_k^{\Lie}\leq\psi_k;\]
and the formal generator function
\[\psi(x)=\sum_{k=0}^\infty\psi_k x^k\]
satisfies the formal IVP
\begin{equation}\psi(0)=0\plabel{eq:stanstart}\end{equation}
\begin{equation}\psi'(x)=\tilde\beta(\psi(x)). \plabel{eq:standev}\end{equation}
\begin{theorem} (a) Interpreted as an analytic function around $0$,
\[\tilde \beta(x)=2+\frac x2-\frac x2\cot \frac x2;\]
and it is convergent in $\intD(0,\pi)$.

(b) The analytic version of IVP \eqref{eq:stanstart}--\eqref{eq:standev} has a solution around $0$ with
convergence radius
\[\boldsymbol\delta=\int_{y=0}^{2\pi} \frac{\mathrm dy}{\tilde \beta(y)}\approx 2.1737374\ldots.\]
As $x\nearrow\boldsymbol\delta$, we have $\psi(x)\nearrow2\pi$ and $\psi'(x)\nearrow+\infty$.
(So, in real function theoretic sense, $\psi(\boldsymbol\delta)=2\pi$ and $\psi'(\boldsymbol\delta)=+\infty$.)
Thus, $\boldsymbol\delta$ is the convergence radius of $\psi$.
\begin{proof}
(a) It follows from  $\frac{\beta(\mathrm ix)+\beta(-\mathrm ix)}2=\frac x2\cot \frac x2$
and the information on the signs of the Bernoulli numbers $B_j$ for $j\geq 2$.

(b) By Pringsheim's theorem it is sufficient to consider the development for $x\geq0$.
Then the standard method of separation of variables can be applied,
\[\frac{\mathrm d\psi}{\tilde \beta(\psi)}=\mathrm dx.\]
Thus the $\psi$ will be the inverse function of $y\mapsto \chi(y)=\int_{t=0}^y \frac{\mathrm dt}{\beta(t)}$,
as long as $\tilde \beta$ is positive, and develops a singularity when $\tilde \beta$ becomes infinite.
(We know that the solution $\psi$ is convex.)
This happens when  $t=2\pi$, which means that the range of the nonsingular
$\chi$ is $\left[0,\int_{t=0}^{2\pi} \frac{\mathrm dt}{\beta(t)}\right)$.
Thus, so is the (nonnegative) domain of the nonsingular $\psi$, and the behaviour around $x\sim \boldsymbol\delta$
follows from the inverse function picture.
\end{proof}
\end{theorem}

\begin{cor}[The standard estimate, cf. Introduction]
In real function theoretic sense,
\[\Theta^{\Lie}(x)\leq\psi(x),\]
showing, in particular, that the convergence radius of $\Theta^{\Lie}(x)$ is at least $\boldsymbol\delta$. \qed
\proofremark{
As there are several passes between the formal, analytic, and real function theoretic viewpoints in this section,
we will be more relaxed with the notation here.
}
\end{cor}

In the context of the Baker--Campbell--Hausdorff formula, a somewhat finer estimate is used, the
``$(1/4)$-commutative version''
 of M\'erigot \cite{Me} (see Biagi \cite{Bia})
Day, So, Thompson \cite{DST}, Blanes, Casas \cite{BC}, which we sketch
(but formulated in our more combinatorial view).
Let us define $\Theta\mathrm P_{k,l}^{\Lie}$, as the minimal possible sum of the absolute values of the coefficients of the
Lie-presentations of $\mu_{k+l}(X_1,\ldots,X_k, Y_1,\ldots,Y_l)$ but with the additional assumption that the variables $Y_i$
commute with each other. Then
$ \Theta\mathrm P_{k,l}^{\Lie}\leq \Theta_{k+l}^{\Lie}$.
Consider the formal generating function
\[\Theta\mathrm P^{\Lie}(x,y)=\sum_{k+l\geq 1}\Theta\mathrm P_{k,l}^{\Lie}x^ky^l. \]
Its coefficients can be estimated by the coefficients of the solution of the  formal IVP
\begin{equation}\grave{\psi}(0,y)=y\plabel{eq:sstanstart}\end{equation}
\begin{equation}\frac{\partial}{\partial x}\grave{\psi}(x,y)=\tilde\beta(\grave{\psi}(x,y)). \plabel{eq:sstandev}\end{equation}
Formally,
\[\grave{\psi}(x,y)=\psi(x+\psi^{-1}(y)).\]
By similar arguments as before,
$\grave{\psi}(x,y)$ will be finite in real analytical sense for $x,y\geq0$ if $x=0$, or $y\leq2\pi$ and
\[x\leq\int_{t=y}^{2\pi} \frac{\mathrm dt}{\tilde \beta(t)}.\]
Then we obtain convergence estimates for the BCH series through $\Gamma^{\Lie}(x,y)\leq \grave{\psi}(x,y)$ for $x,y\geq0$.
(Cf. the already cited sources).
We see, however, that the situation is not very well adapted to the BCH case, because in the form just given not even its $2\cdot(1/4)$-commutative features are utilized, meaning that the same argument also applies, more generally,
to the Magnus expansion with a commutative segment (in an appropriate form).
This makes it hard to compete with the estimate  $\mathrm C_2\leq \mathrm C_2^{\Lie}$
when it comes to the cumulative norm.
(Hence, we will refrain from discussing the improvements of the standard estimate in the BCH case,
although the arguments can be adapted to it.)

Let us return to the Magnus expansion.
By the earlier discussions, in terms of the convergence radius of $\Theta^{\Lie}$,
we have a gap between $\boldsymbol\delta$ and $4$.
Closing this gap likely requires some deeper insight.
However, an advantage of the algebraic formalism is that
it offers several ways to improve the standard estimate a bit.
We show some methods.
We are less interested in numerical constants but that they show greater convergence domains.

\begin{method}[Forced coefficients]
One can observe that
\[\psi(x)=x+\frac14\,{x}^{2}+{\frac {5}{72}}\,{x}^{3}+{\frac{11}{576}}\,{x}^{4}+{\frac {479}{86400}}\,{x}^{5}+{\frac {1769}{1036800}}\,{x}^{6}+\ldots\]
in contrast to
\[\Theta^{\Lie}(x)=x+\frac14\,{x}^{2}+\frac1{18}\,{x}^{3}+{\frac {1}{72}}{x}^{4}+{\frac {1}{
300}}{x}^{5}+{\frac {37}{43200}}\,{x}^{6}+\ldots.\]
Now,
\[\frac{\mathrm d\Theta^{\Lie}(x)}{\mathrm dx}-\tilde \beta(\Theta^{\Lie}(x))=
-\frac1{24}\,{x}^{2}-{\frac {1}{72}}{x}^{3}-{\frac {53}{8640}}\,{x}^{4}-{
\frac {11}{4320}}\,{x}^{5}+\ldots.\]
Thus, when we solve the IVP
\begin{equation}\hat\psi(0)=0\plabel{eq:hstanstart}\end{equation}
\begin{equation}\hat\psi'(x)=\tilde\beta(\hat\psi(x))\underbrace{-\frac1{24}\,{x}^{2}-{\frac {1}{72}}{x}^{3}-{\frac {53}{8640}}\,{x}^{4}-{
\frac {11}{4320}}\,{x}^{5}}_{-\Delta_6(x)}; \plabel{eq:hstandev}\end{equation}
we find that $\hat\psi(x)$ has the same coefficients as $\Theta^{\Lie}(x)$ up to order $6$, but after that
the majorizing property relative to $\Theta^{\Lie}(x)$ still holds; $\Theta_k^{\Lie}\leq\hat\psi_k\leq\psi_k$.

In order to demonstrate the larger convergence radius,
let us compare $\hat\psi$ to $\psi$ by a crude estimate.
For $x\in[0,\boldsymbol\delta]$, we see that
\[\hat\psi'(x)-\psi'(x)=\tilde\beta(\hat\psi(x))-\Delta_6(x)-\tilde\beta(\psi(x))\leq -\Delta_6(x). \]
Integrating in  $x$, we find
\[\hat\psi(x)-\psi(x)\leq -\int_{t=0}^{x} \Delta_6(t) \mathrm dt.\]
Thus, in particular,
\[\hat\psi(\boldsymbol\delta )\leq\psi(\boldsymbol\delta) -\int_{t=0}^{\boldsymbol\delta} \Delta_6(t) \mathrm dt.\]
In particular, $\hat\psi(\boldsymbol\delta)$ is strictly smaller than $\psi(\boldsymbol\delta)$.
At this point, even if we continue with the slower \eqref{eq:standev}, that would give an extra length
\[\hat L=\boldsymbol\delta- \chi\left(\psi(\boldsymbol\delta) -\int_{t=0}^{\boldsymbol\delta} \Delta_6(t) \mathrm dt\right)\]
for further development. (But we know that it gives even more.)
So, we know that the convergence radius of $\hat \phi$ is bigger than $\boldsymbol\delta+\hat L$.
In fact, we have the estimate
\[\hat\psi(x)\leq\begin{cases}
\psi(x) -\int_{t=0}^{x} \Delta_6(t) \mathrm dt&\text{if }x\in[0,\boldsymbol\delta]\\
\psi(x-\hat L)&\text{if }x\in[\boldsymbol\delta,\boldsymbol\delta+\hat L].
\end{cases}\]
Numerically $\hat L=0.0074001\ldots$, thus $\boldsymbol\delta+\hat L=2.1811375\ldots$ is obtained for
a larger convergence radius.

Our estimates above were really simplistic, though;
more precise numerical results show that the
convergence radius of $\hat\psi$ is around $2.2762\ldots$
\qedexer\end{method}

What hinders  the previous method is that we use the same naive recursion mechanism based on the Magnus recursion as originally.
We can achieve better results if we use recursion of higher order.
The next discussion will be essentially based on the identity
\[\sum_{j=1}^kP(X_1,\ldots,X_{j-1},[X_j,Y],\ldots,X_{j+1},\ldots,X_{k} )
=[P(X_1,\ldots,X_{k} ),Y],\]
where $P$ is a Lie polynomial.
This allows to reduce the size of some expressions.

In order to compactify our formulas, let us introduce some notation.
Instead of explaining it in advance, we show how the Magnus recursion can be expressed in this notation.

Expansion in the first variable ($Z_1$):
\[\mu=Z_1-\frac12[\mu,Z_1]+\sum_{j=1}^\infty \beta_{2j}[\underbrace{\mu,\ldots,\mu}_{2j},Z_1].\]

expansion in the last variable ($Z_n$):
\[\mu=Z_n+\frac12[\mu,Z_n]+\sum_{j=1}^\infty \beta_{2j}[\underbrace{\mu,\ldots,\mu}_{2j},Z_n].\]

What happens is that we consider Lie polynomials in $Z_1,\ldots,Z_n$, and
terms in the expressions are understood so that whenever we have  $\mu$'s
with some unspecified variables, then  the unnoted variables are distributed among
them  with no multiplicities, and ascendingly in every $\mu$.
The higher commutators should be resolved as
\[[X_1,\ldots,X_{k-1},X_k]=(\ad X_1)\ldots(\ad X_{k-1})X_k.\]
Clearly, one should be careful with this notation, but it has the advantage of being short.

We also use the notation
\[\tilde{\tilde\beta}(x)=\tilde\beta(x)-1-\frac12x.\]
\begin{method}[Magnus recursion of second order]

In the standard approach it did not matter if we used expansion by the first or last variable.
Here we take expansion by two variables, and, in order to gain  a little additional improvement,
we combine this with symmetrization.

If we expand in $Z_1$, and later in $Z_n$, then we find
\begin{align}
\mu=&Z_1\notag\\
 &+\frac12[Z_1,Z_n]+\frac14[Z_1,[\mu,Z_n]]+\sum_{k=1}^\infty \beta_{2k}\frac12[Z_1,[\underbrace{\mu,\ldots,\mu}_{2k},Z_n]]\notag\\
 &+\sum_{j=1}^\infty \beta_{2j}[\underbrace{\mu,\ldots,Z_n,\ldots,\mu}_{2j},Z_1]\notag\\
  &+\frac12\sum_{j=1}^\infty \beta_{2j}[[\underbrace{\mu,\ldots,\mu}_{2j},Z_1],Z_n]
  -\frac12\sum_{j=1}^\infty \beta_{2j}[\underbrace{\mu,\ldots,\mu}_{2j},[Z_1,Z_n]]\notag\\
 &+\sum_{k=1}^\infty\beta_{2k}\sum_{j=1}^\infty \beta_{2j}[[\underbrace{\mu,\ldots,\mu}_{2j},Z_1],[\underbrace{\mu,\ldots,\mu}_{2k-1},Z_n]]\notag\\
 &-\sum_{k=1}^\infty\beta_{2k}\sum_{j=1}^\infty \beta_{2j}[\underbrace{\mu,\ldots,\mu}_{2j},Z_1,\underbrace{\mu,\ldots,\mu}_{2k-1},Z_n].\notag
\end{align}
If we expand in $Z_n$, and later in $Z_1$, then we find
\begin{align}
\mu=&Z_n\notag\\
 &+\frac12[Z_1,Z_n]+\frac14[[Z_1,\mu],Z_n]+\sum_{k=1}^\infty \beta_{2k}\frac12[[\underbrace{\mu,\ldots,\mu}_{2k},Z_1],Z_n]\notag\\
 &+\sum_{j=1}^\infty \beta_{2j}[\underbrace{\mu,\ldots,Z_1,\ldots,\mu}_{2j},Z_n]\notag\\
  &+\frac12\sum_{j=1}^\infty \beta_{2j}[Z_1,[\underbrace{\mu,\ldots,\mu}_{2j},Z_n]]
  -\frac12\sum_{j=1}^\infty \beta_{2j}[\underbrace{\mu,\ldots,\mu}_{2j},[Z_1,Z_n]]\notag\\
 &+\sum_{k=1}^\infty\beta_{2k}\sum_{j=1}^\infty \beta_{2j}[[\underbrace{\mu,\ldots,\mu}_{2j},Z_n],[\underbrace{\mu,\ldots,\mu}_{2k-1},Z_1]]\notag\\
 &-\sum_{k=1}^\infty\beta_{2k}\sum_{j=1}^\infty \beta_{2j}[\underbrace{\mu,\ldots,\mu}_{2j},Z_n,\underbrace{\mu,\ldots,\mu}_{2k-1},Z_1].\notag
\end{align}
The averaged (symmetrized) expression is
\begin{align}
\mu=&Z+\frac12[Z_1,Z_n]\notag\\
 &+\frac18[Z_1,[\mu,Z_n]]+\frac18[[Z_1,\mu],Z_n]\notag\\
 &+\frac12\sum_{j=2}^\infty \beta_{2j}[\underbrace{\mu,\ldots,Z_n,\ldots,\mu}_{2j-1},\underbrace{\mu}_1,Z_1]+\frac12\sum_{j=2}^\infty \beta_{2j}[\underbrace{\mu,\ldots,Z_1,\ldots,\mu}_{2j-1},\underbrace{\mu}_1,Z_n]\notag\\
  &+\frac12\sum_{j=1}^\infty \beta_{2j}[[\underbrace{\mu,\ldots,\mu}_{2j},Z_1],Z_n]
  +\sum_{k=1}^\infty \beta_{2k}\frac12[Z_1,[\underbrace{\mu,\ldots,\mu}_{2k},Z_n]]\notag\\
 & -\frac12\sum_{j=1}^\infty \beta_{2j}[\underbrace{\mu,\ldots,\mu}_{2j},[Z_1,Z_n]]\notag\\
 &+\frac12\sum_{k=1}^\infty\beta_{2k}\sum_{j=1}^\infty \beta_{2j}[[\underbrace{\mu,\ldots,\mu}_{2j},Z_1],[\underbrace{\mu,\ldots,\mu}_{2k-1},Z_n]]\notag\\
 &+\frac12\sum_{k=1}^\infty\beta_{2k}\sum_{j=1}^\infty \beta_{2j}[[\underbrace{\mu,\ldots,\mu}_{2j},Z_n],[\underbrace{\mu,\ldots,\mu}_{2k-1},Z_1]]\notag\\
 &-\frac12\sum_{k=1}^\infty\beta_{2k}\sum_{j=1}^\infty \beta_{2j}[\underbrace{\mu,\ldots,\mu}_{2j},Z_1,\underbrace{\mu,\ldots,\mu}_{2k-1},Z_n]\notag\\
 &-\frac12\sum_{k=1}^\infty\beta_{2k}\sum_{j=1}^\infty \beta_{2j}[\underbrace{\mu,\ldots,\mu}_{2j},Z_n,\underbrace{\mu,\ldots,\mu}_{2k-1},Z_1].\notag
\end{align}

Using these, we can develop a majorizing series $\psi$ for $\Theta^{\Lie}$
by the recursion (formal IVP)
\begin{equation}\check\psi(0)=0\plabel{eq:masod1}\end{equation}
\begin{equation}\check\psi'(0)=1\plabel{eq:masod2}\end{equation}
\begin{equation}\check\psi''(x)=f(\psi(x)),\plabel{eq:masod3}\end{equation}
where
\begin{align}
f(x)=&\frac12+\frac14x+\tilde{\tilde\beta}(x)'-\frac1x\tilde{\tilde\beta}(x)\notag+\frac32\tilde{\tilde\beta}(x)
+\frac2x\tilde{\tilde\beta}(x)^2\notag\\
=&2+\frac x2+\frac1x-2\,\cot \left( \frac x2 \right)-\frac34\,x\cot \left( \frac x2 \right)
+\frac34\,x \left(\cot \left(\frac x2\right)  \right) ^{2}.\notag
\end{align}

The IVP \eqref{eq:masod1}--\eqref{eq:masod3} is one of the classically treatable ones, and  its leads to convergence radius
\[\boldsymbol\delta_2=\int_{u=0}^{2\pi}\frac{\mathrm du}{\sqrt{1+2\int_{t=0}^uf(t)\,\mathrm dt  }}\approx2.281\ldots\]

This improvement is still small, but better than in the previous case.
\qedexer\end{method}
In general, we are interested in the size of $\mu$.
However, in its estimation, the size of $\beta(\ad\mu)$ played a role,
which we estimated from the size of $\mu$ naively.
We can do better keeping a separate check on the size of $\beta(\ad\mu)$.
More precisely, we will keep a check on the size of  $\tilde{\tilde\beta}(\ad\mu)$.
\begin{method}[A simplest compartmentalization.]
Consider the equation
\[\mu=Z_1-\frac12[\mu,Z_1]+\sum_{j=1}^\infty \beta_{2j}[\underbrace{\mu,\ldots,\mu}_{2j},Z_1]\]
and it consequence
\begin{align}
 \sum_{k=1}^\infty\beta_{2k}[\underbrace{\mu,\ldots,\mu}_{2k},W]=&
\sum_{k=1}^\infty\beta_{2k}[\underbrace{\mu,\ldots,Z_1+\sum_{j=1}^\infty
\beta_{2j}[\underbrace{\mu,\ldots,\mu}_{2j},Z_1],\ldots\mu}_{2k},W]\plabel{line:har}\\
&-\frac12\sum_{k=1}^\infty\beta_{2k}[[\underbrace{\mu,\ldots,\mu}_{2k},W],Z_1]
+ \frac12\sum_{k=1}^\infty\beta_{2k}[\underbrace{\mu,\ldots,\mu}_{2k},[W,Z_1]].\notag
\end{align}
This leads to the majorizing sytem (formal IVP)
\begin{equation}\psi(0)=0,\qquad \tilde{\tilde\psi}(0)=0\end{equation}
\begin{equation}\psi'(x)=1+\frac12\psi(x)+\tilde{\tilde\psi}(x),\end{equation}
\begin{equation}\tilde{\tilde\psi}'(x)=\tilde{\tilde\beta}'(\psi(x))(1+\tilde{\tilde\psi}(x)) +\tilde{\tilde\psi}(x).\end{equation}
In this present form, this differential equation blows up around $x=2.204\ldots$,
which leads to a quite modest lower estimate for the convergence radius.
However, compartmentalization schemes like that, in general, allow to separate various algebraic patterns
in the Magnus expansion.
We do not pursue this direction in its full power now but we slightly improve this example.
\qedexer\end{method}
The disadvantage of the previous method is that in the RHS of line \eqref{line:har}
we still use exponential estimates.
The ideal thing would be keeping a check on the size $(\ad \mu)^k$ for every $k$, but this is
just too complicated for us to do here.
However, we will do this partially.
Let us use the notation
\[\beta^{(e)}(x)=\frac{x^2}{4\pi^2-x^2}=\sum_{j=1}^\infty\left(\frac{x}{2\pi}\right)^{2j}\]
\[\beta^{(o)}(x)=\frac{x^3}{2\pi(4\pi^2-x^2)}=\sum_{j=1}^\infty\left(\frac{x}{2\pi}\right)^{2j+1}\]
\[\mathring\beta(x)=\sum_{j=1}^\infty\left(\frac{x}{2\pi}\right)^{2j}2\sum_{N=2}^\infty\frac1{N^{2j}}\]
Then
\begin{equation}\tilde{\tilde\beta}(x)=2\beta^{(e)}(x)+ \mathring\beta(x).\plabel{eq:tilbetadec}\end{equation}
\begin{method}[A slightly more sophisticated compartmentalization.]
Here we keep track on the size of $\mu$, $\beta^{(e)}(\ad\mu)$,  $\beta^{(o)}(\ad\mu)$,  $\mathring\beta(\ad\mu)$.
The relevant equations are
\[\mu=Z_1-\frac12[\mu,Z_1]+\sum_{j=1}^\infty \beta_{2j}[\underbrace{\mu,\ldots,\mu}_{2j},Z_1]\]

\begin{multline}
 \sum_{k=1}^\infty\frac1{(2\pi)^{2k}}[\underbrace{\mu,\ldots,\mu}_{2k},W]=
\sum_{k=1}^\infty\frac1{(2\pi)^{2k}}[\underbrace{\mu,\ldots,Z_1+\sum_{j=1}^\infty
\beta_{2j}[\underbrace{\mu,\ldots,\mu}_{2j},Z_1],\ldots\mu}_{2k},W]\notag\\
-\frac12\sum_{k=1}^\infty\frac1{(2\pi)^{2k}}[[\underbrace{\mu,\ldots,\mu}_{2k},W],Z_1]
+ \frac12\sum_{k=1}^\infty\frac1{(2\pi)^{2k}}[\underbrace{\mu,\ldots,\mu}_{2k},[W,Z_1]],\notag
\end{multline}
\begin{multline}
 \sum_{k=1}^\infty\frac1{(2\pi)^{2k+1}}[\underbrace{\mu,\ldots,\mu}_{2k+1},W]=
\sum_{k=1}^\infty\frac1{(2\pi)^{2k+1}}[\underbrace{\mu,\ldots,Z_1+\sum_{j=1}^\infty
\beta_{2j}[\underbrace{\mu,\ldots,\mu}_{2j},Z_1],\ldots\mu}_{2k+1},W]\notag\\
-\frac12\sum_{k=1}^\infty\frac1{(2\pi)^{2k+1}}[[\underbrace{\mu,\ldots,\mu}_{2k+1},W],Z_1]
+ \frac12\sum_{k=1}^\infty\frac1{(2\pi)^{2k+1}}[\underbrace{\mu,\ldots,\mu}_{2k+1},[W,Z_1]],\notag
\end{multline}
\begin{multline}
 \sum_{k=1}^\infty\mathring\beta_{2k}[\underbrace{\mu,\ldots,\mu}_{2k},W]=
\sum_{k=1}^\infty\mathring\beta_{2k}[\underbrace{\mu,\ldots,Z_1+\sum_{j=1}^\infty
\beta_{2j}[\underbrace{\mu,\ldots,\mu}_{2j},Z_1],\ldots\mu}_{2k},W]\notag\\
-\frac12\sum_{k=1}^\infty\mathring\beta_{2k}[[\underbrace{\mu,\ldots,\mu}_{2k},W],Z_1]
+ \frac12\sum_{k=1}^\infty\mathring\beta_{2k}[\underbrace{\mu,\ldots,\mu}_{2k},[W,Z_1]].\notag
\end{multline}
This leads to the IVP
\begin{equation}\psi(0)=0,\quad \psi^{(e)}(0)=0,\quad \psi^{(o)}(0)=0,\quad \mathring\psi(0)=0,\quad  \end{equation}
\[\psi'(x)=1+\frac12\psi(x)+2 \psi^{(e)}(x)+\mathring\psi(x),\]
\[\psi^{(e)\prime}(x)=\frac1{2\pi}2\left(\frac{\psi(x)}{2\pi}+\psi^{(o)}(x)\right)(1+\psi^{(e)}(x) )
(1+2 \psi^{(e)}(x)+\mathring\psi(x))+\psi^{(e)}(x),\]
\[\psi^{(o)\prime}(x)=\frac1{2\pi}\left(2\psi^{(e)}(x)  +\psi^{(e)}(x)^2+\left(\frac{\psi(x)}{2\pi}+\psi^{(o)}(x)\right)^2 \right)
(1+2 \psi^{(e)}(x)+\mathring\psi(x))+\psi^{(o)}(x),\]
\[\mathring\psi'(x)=\mathring\beta'(\psi(x))(1+2 \psi^{(e)}(x)+\mathring\psi(x)) +\mathring\psi(x).\]
Numerical results show that this  IVP blows up at $x=2.297\ldots.$
This is our best lower bound for the convergence radius up to now.
\qedexer\end{method}
\begin{method}[A variant of the previous method]
Let
\[\beta^{(ee)}(x)=\sum_{j=1}^\infty\left(\frac{x}{2\cdot2\pi}\right)^{2j},\qquad
\beta^{(oo)}(x)=\sum_{j=1}^\infty\left(\frac{x}{2\cdot2\pi}\right)^{2j+1}.\]
According to this we have the analogue
\[\tilde{\tilde\beta}(x)=2\beta^{(e)}(x)+2\beta^{(ee)}(x)+ \ddot{\beta}(x)\]
of \eqref{eq:tilbetadec}.
The sizes of the expressions $\mu,\beta^{(e)}(\ad\mu),\ldots,\ddot\beta(\ad\mu)$ are described by the series
$\Theta^{\Lie}(x),$ $\Theta^{(e)}(x),\ldots,\ddot\Theta(x)$.
The appropriate recursion relations imply
\[\Theta(0)=0,\quad \Theta^{(e)}(0)=0,\quad \Theta^{(o)}(0)=0,
\quad \Theta^{(ee)}(0)=0,\quad \Theta^{(oo)}(0)=0,\quad \ddot\Theta(0)=0,  \]
\[\Theta^{\Lie\prime}(x)\leq\Delta\Theta(x)+\frac12\Theta^{\Lie}(x),\]
\[\Theta^{(e)\prime}(x)\leq\frac1{2\pi}2\left(\frac{\Theta^{\Lie}(x)}{2\pi}+\Theta^{(o)}(x)\right)(1+\Theta^{(e)}(x) )
\Delta\Theta(x)+\Theta^{(e)}(x),\]
\[\Theta^{(o)\prime}(x)\leq\frac1{2\pi}\left(2\Theta^{(e)}(x)  +\Theta^{(e)}(x)^2+
\left(\frac{\Theta^{\Lie}(x)}{2\pi}+\Theta^{(o)}(x)\right)^2 \right)\Delta\Theta(x)+\Theta^{(o)}(x),\]
\[\Theta^{(ee)\prime}(x)\leq\frac1{2\cdot2\pi}2\left(\frac{\Theta^{\Lie}(x)}{2\cdot2\pi}+\Theta^{(oo)}(x)\right)(1+\Theta^{(ee)}(x) )
\Delta\Theta(x)+\Theta^{(ee)}(x),\]
\[\Theta^{(oo)\prime}(x)\leq\frac1{2\cdot2\pi}\left(2\Theta^{(ee)}(x)  +\Theta^{(ee)}(x)^2+
\left(\frac{\Theta^{\Lie}(x)}{2\cdot2\pi}+\Theta^{(oo)}(x)\right)^2 \right)\Delta\Theta(x)+\Theta^{(oo)}(x),\]
\[\ddot\Theta'(x)\leq\ddot\beta'(\Theta^{\Lie}(x))\Delta\Theta(x) +\ddot\Theta(x).\]
where
\[\Delta\Theta(x)\equiv1+2 \Theta^{(e)}(x)+2 \Theta^{(ee)}(x)+\ddot\Theta(x)\]
However, from the series expansion we also know that
\[\ddot\Theta(x)\leq2\cdot2^2\sum_{N=3}^\infty\frac1{N^2}\cdot \Theta^{(ee)}(x).\]
Thus, we also have
 \begin{equation}\Theta(0)=0,\quad \Theta^{(e)}(0)=0,\quad \Theta^{(o)}(0)=0,
\quad \Theta^{(ee)}(0)=0,\quad \Theta^{(oo)}(0)=0, \end{equation}
\[\Theta^{\Lie\prime}(x)\leq\Delta\Theta(x)+\frac12\Theta^{\Lie}(x),\]
\[\Theta^{(e)\prime}(x)\leq\frac1{2\pi}2\left(\frac{\Theta^{\Lie}(x)}{2\pi}+\Theta^{(o)}(x)\right)(1+\Theta^{(e)}(x) )
\Delta\Theta(x)+\Theta^{(e)}(x),\]
\[\Theta^{(o)\prime}(x)\leq\frac1{2\pi}\left(2\Theta^{(e)}(x)  +\Theta^{(e)}(x)^2+
\left(\frac{\Theta^{\Lie}(x)}{2\pi}+\Theta^{(o)}(x)\right)^2 \right)\Delta\Theta(x)+\Theta^{(o)}(x),\]
\[\Theta^{(ee)\prime}(x)\leq\frac1{2\cdot2\pi}2\left(\frac{\Theta^{\Lie}(x)}{2\cdot2\pi}+\Theta^{(oo)}(x)\right)(1+\Theta^{(ee)}(x) )
\Delta\Theta(x)+\Theta^{(ee)}(x),\]
\[\Theta^{(oo)\prime}(x)\leq\frac1{2\cdot2\pi}\left(2\Theta^{(ee)}(x)  +\Theta^{(ee)}(x)^2+
\left(\frac{\Theta^{\Lie}(x)}{2\cdot2\pi}+\Theta^{(oo)}(x)\right)^2 \right)\Delta\Theta(x)+\Theta^{(oo)}(x),\]
where
\[\Delta\Theta(x)\equiv1+2 \Theta^{(e)}(x)+2\left(2^2\sum_{N=2}^\infty\frac1{N^2}\right) \Theta^{(ee)}(x).\]
We can draw a formal IVP for a majorizing series upon these inequalities, as before.
It turns out that this system blows up at $x=2.293\ldots$ which is not so good as in the case of the previous method.
Nevertheless, the system itself is polynomial, which offers some technical advantages.
\qedexer\end{method}
Now, the various methods can be combined and refined, leading to even better estimates.
Also note that using standard techniques, one can make the numerical results above completely robust.
However, as we will obtain stronger estimates using different methods later, we
do not make those numerical results more exact here.

Regarding the BCH expansion, the same convergence improvement methods apply.
As we said, we refrain from working out cases here, in favour of the resolvent method.
\snewpage
\section{The extension of the Banach-Lie norm in the formal case}\plabel{sec:normext}
We have already seen the Lie norm $\|\cdot\|_{\ell^1}$ on $\mathrm F^{\Lie}[Y_\lambda:\lambda\in\Lambda]$.
On the other hand, we can consider $\mathrm F^{\Lie}[Y_\lambda:\lambda\in\Lambda]\subset \mathrm F[Y_\lambda:\lambda\in\Lambda]$
through the commutator inclusion.
We will show that the Lie norm $\|\cdot\|_{\ell^1}$ extends to an associative algebraic norm on $\mathrm F[Y_\lambda:\lambda\in\Lambda]$.
\begin{lemma}\plabel{lem:adext}
Consider the Lie polynomial $P(Y_1,\ldots, Y_n)$, and a further, independent generator element $Z$. Then
\[\|P(Y_1,\ldots, Y_n)\|_{\ell^1}=\|[P(Y_1,\ldots, Y_n),Z]\|_{\ell^1}.\]
\begin{proof}

$\|P(Y_1,\ldots, Y_n)\|_{\ell^1}\geq\|[P(Y_1,\ldots, Y_n),Z]\|_{\ell^1}$ is rather obvious from the definition.
On the other hand, consider a minimal presentation of $[P(Y_1,\ldots, Y_n),Z]$.
The $P$ can be assumed to be homogeneous, and of shape
\[[P(Y_1,\ldots, Y_n),Z]=\sum_{j=1}^s \theta_j[M_{j,1},[\ldots,[M_{j,k_s},Z]\ldots],\]
where the $M_{j,k}$ are monomials of $Y_1,\ldots,Y_n$.
Due to the freeness, $[\cdot,Z]$ can be assumed to act as the degree derivation, thus
\[P(Y_1,\ldots, Y_n)\cdot (\deg P)= \sum_{j=1}^s \theta_j[M_{j,1},[\ldots,[M_{j,k_s}]\ldots]\cdot(\deg M_{j,k_s}).\]
So,
\[P(Y_1,\ldots, Y_n)= \sum_{j=1}^s \frac{\deg M_{j,k_s}}{\deg P}\theta_j[M_{j,1},[\ldots,[M_{j,k_s}]\ldots].\]
Then $\frac{\deg M_{j,k_s}}{\deg P}\leq1$ implies the inequality in the other direction.
(This is a variant of an argument used in the proof the Dynkin--Specht--Wever lemma.)
\end{proof}
\end{lemma}
The extended norm $\|Q(Y_1,\ldots,Y_k)\|_{\ell^1}'$ for $Q(Y_1,\ldots,Y_k)\in\mathrm F[Y_\lambda:\lambda\in\Lambda] $
can be defined in two ways.
The first one is to take
\begin{multline}\|Q(Y_1,\ldots,Y_k)\|_{\ell^1}'=\inf\Biggl\{\sum_{i=1}^s |\theta_j|\,:\,
Q(Y_1,\ldots,Y_k)=\sum_{i=1}^s \theta_j\,M_{i,1}\cdot \ldots \cdot M_{i,j_s},\\
\text{such that the $M_{i,j}$ are commutator monomials of the $Y_1,\ldots,Y_n$}\Biggr\}.
\notag\end{multline}
The second one is take
\[\|Q(Y_1,\ldots,Y_k)\|_{\ell^1}'=\|Q(\ad Y_1,\ldots,\ad Y_k) Z\|_{\ell^1}\quad(\text{$Z$ is an independent generator}).\]
The equivalence of the definitions and the extension property follow, taking the adjoint representation on $Z$, from  the previous lemma.

We remark that the following qualitative version of the DSW lemma holds:
\begin{cor}
Consider the Lie polynomial $P(Y_1,\ldots, Y_n)$.
Then any $\|\cdot\|_{\ell^1}'$-minimal presentation of $P$ in
$\mathrm F[Y_\lambda:\lambda\in\Lambda]$ (linear combination of products of commutator monomials with minimal sum of absolute values for the coefficients),
$0$ coefficient terms omitted,
is automatically of shape of an  $\|\cdot\|_{\ell^1}$-minimal presentation of $P$ in
 $\mathrm F^{\Lie}[Y_\lambda:\lambda\in\Lambda]$ (linear combination  of Lie monomials with minimal sum of absolute values for the coefficients).
\begin{proof}
This follows from the proof of Lemma \ref{lem:adext}, because, for nonzero coefficients, the equalities $\frac{\deg M_{j,k_s}}{\deg P}=1$ should hold.
\end{proof}
\end{cor}

Note that $\|\cdot\|'_{\ell^1}\leq|\cdot|_{\ell^1}$.
$\mathrm F[Y_\lambda:\lambda\in\Lambda]$ (and $\mathrm F^{1}[Y_\lambda:\lambda\in\Lambda]$)
will get completed to $\mathrm F^{1,(\Lie)}[Y_\lambda:\lambda\in\Lambda]$ (algebra) which naturally extends
$\mathrm F^{1,\Lie}[Y_\lambda:\lambda\in\Lambda]$ (Lie-algebra).
In what follows, we simply write $\|\cdot\|_{\ell^1}$ instead of $\|\cdot\|_{\ell^1}'$.

The tautological measure  construction also extends to this norm.
A possible interpretation is as follows.
If $X\in \mathrm F^1(I)$, then it can be decomposed as \[X=\sum_{k=0}^\infty X_k\]
according to homogeneous terms.
Then $X_k$ can be written as
\[X_k=\int_{(t_1,\ldots,t_k)\in[0,1]^k}h_k(t_1,\ldots,t_k) \,\mathrm Z^1_I(t_1)\cdot\ldots\cdot\mathrm Z^1_I(t_k)\]
where $h_k$ is a Lebesgue integrable function on $I^k$.
Then we define
\[\|X_k\|_{\ell^1}=\int_{t_1<\ldots< t_k\in I}\left\|\sum_{\sigma\in\Sigma_k}h_k(t_{\sigma(1)},\ldots,t_{\sigma(k)})
Y_{\sigma(1)}\cdot\ldots\cdot Y_{\sigma(k)}\right\|_{\ell^1}
\,\mathrm dt_1\ldots\mathrm dt_k.\]
More generally, we let
\[\|X\|_{\ell^1}=\sum_{k=0}^\infty \|X_k\|_{\ell^1}.\]
In general, $\|X\|_{\ell^1}\leq|X|_{\ell^1}$.
Actually, $\|\cdot\|_{\ell^1}$ is indeed a norm on $\mathrm F^1(I)$ not just a seminorm.
(Because $\|\cdot\|_{\ell^1}\leq|\cdot|_{\ell^1}$ are comparable in any finite degree.)
However,  $\mathrm F^1(I)$ needs to be completed in $\|\cdot\|_{\ell^1}$.
This yields $\mathrm F^{1,(\Lie)}(I)$ with a natural (weakly) contractive map
$\mathrm F^{1}(I)\rightarrow\mathrm F^{1,(\Lie)}(I)$, compatible to the algebra struture.
Here $\mathrm Z^{\Lie}_{[a,b)}$ is the image of $Z_{[a,b)}=\int_{t\in I}1_{[a,b)}\mathrm Z_I^{1}(t)$;
$\mathrm Z^1_I$ is the tautological (interval) measure sending $[a,b)\mapsto Z^{\Lie}_{[a,b)}$; in that way the image of $X_k$
in $\mathrm F^{1,(\Lie)}(I)$ is equal to
\[\int_{(t_1,\ldots,t_k)\in[0,1]^k}h_k(t_1,\ldots,t_k) \,\mathrm Z^{1,\Lie}_I(t_1)\cdot\ldots\cdot\mathrm Z^{1,\Lie}_I(t_k).\]
(Here we have identified the Banach--Lie algebraic $\mathrm Z^{1,\Lie}_I$ with the Banach algebraic  $\mathrm Z^{1,(\Lie)}_I$, etc.)

\begin{remark}
The advantages of the Banach algebraic extension are already implicit in some recursions methods of Section \ref{sec:MagnusStandard}.
\qedremark
\end{remark}
\snewpage

\section{The resolvent method and the Magnus expansion}\plabel{sec:resmagnus}

In this section we apply the resolvent method in order to get estimates to the
convergence radius of the Magnus expansion in the Banach--Lie setting.
Here we concentrate on obtaining actual estimates.
A parallel but more detailed and general discussion can be found in Part IA \cite{L15}.

\subsection{The spectral estimates}
\plabel{ss:specest}
~\\

Now, our objective is to estimate the convergence radius $\mathrm C^{\Lie}_\infty$ of
\[\Theta^{\Lie}(x)=\sum_{k=1}^\infty\Theta^{\Lie}_k,\]
where
\begin{align*}
\Theta_k^{\Lie}&=\left\|\mu_k(Y_1,\ldots,Y_k)\right\|_{\ell^1}\notag \\
&=\left\|\int_{t_1\leq\ldots\leq t_k\in[0,1]} \mu_k(\mathrm Z^{1,\Lie}_{[0,1]}(t_1)\ldots \mathrm Z^{1,\Lie}_{[0,1]}(t_k)) \right\|_{\ell^1}
\notag\\
&=\left\|\int_{\lambda=0}^1\int_{\mathbf t=(t_1,\ldots,t_k)\in[0,1]^k} \lambda^{\asc(\mathbf t)}
(\lambda-1)^{\des(\mathbf t)}\mathrm Z^{1,\Lie}_{[0,1]}(t_1)\ldots \mathrm Z^{1,\Lie}_{[0,1]}(t_k) \,\mathrm d\lambda\right\|_{\ell^1}.
\end{align*}

In the spirit of the resolvent method, we will examine
 the convergence radius of
\[\Theta^{(\lambda),\Lie}(x)=\sum_{k=1}^\infty\Theta^{(\lambda),\Lie}_k,\]
where
\begin{align*}
\Theta_k^{(\lambda),\Lie}&=\left\|\mu_k^{(\lambda)}(Y_1,\ldots,Y_k)\right\|_{\ell^1}\notag \\
&=\left\|\int_{t_1\leq\ldots\leq t_k\in[0,1]} \mu^{(\lambda)}_k(\mathrm Z^{1,\Lie}_{[0,1]}(t_1)\ldots \mathrm Z^{1,\Lie}_{[0,1]}(t_k)) \right\|_{\ell^1}
\notag\\
&=\left\|\int_{\mathbf t=(t_1,\ldots,t_k)\in[0,1]^k} \lambda^{\asc(\mathbf t)}
(\lambda-1)^{\des(\mathbf t)}\mathrm Z^{1,\Lie}_{[0,1]}(t_1)\ldots \mathrm Z^{1,\Lie}_{[0,1]}(t_k) \right\|_{\ell^1}.
\end{align*}

Let $1\leq q<k$ and $h=\lfloor(k-1)/p\rfloor$.
Applying the submultiplicative property of $\|\cdot\|_{\ell^1}$, we find
\begin{align}
\Theta_k^{(\lambda),\Lie}&\leq\int_{t_0,t_{p},\ldots, t_{ph}\in[0,1]}\plabel{eq:multdec}\\\notag
&\quad\| \mathrm Z^{1,\Lie}_{[0,1]}(t_0)\|_{\ell^1}\biggl\|\int_{\mathbf t_1=(t_1,\ldots,t_{p-1})\in[0,1]^k} \\\notag
&\qquad\lambda^{\asc(t_0,\mathbf t_1,t_p)}
(\lambda-1)^{\des(t_0,\mathbf t_1,t_p)}\mathrm Z^{1,\Lie}_{[0,1]}(t_1)\ldots \mathrm Z^{1,\Lie}_{[0,1]}(t_{p-1}) \biggr\|_{\ell^1}\\\notag
&\quad\ldots\\\notag
&\quad\| \mathrm Z^{1,\Lie}_{[0,1]}(t_{p(h-1)})\|_{\ell^1}\biggl\|\int_{\mathbf t_h=(t_{(h-1)p+1},\ldots,t_{hp-1})\in[0,1]^k}\\\notag
&\qquad\lambda^{\asc(t_{(h-1)p},\mathbf t_h,t_{hp})}
(\lambda-1)^{\des(t_{(h-1)p},\mathbf t_h,t_{hp})}\mathrm Z^{1,\Lie}_{[0,1]}(t_{h(p-1)+1})\ldots \mathrm Z^{1,\Lie}_{[0,1]}(t_{hp-1}) \biggr\|_{\ell^1}\\\notag
&\quad\| \mathrm Z^{1,\Lie}_{[0,1]}(t_{ph})\|_{\ell^1}\int_{t_{hp+1},\ldots, t_{k}\in[0,1]} \\\notag
&\qquad\lambda^{\asc(t_{hp},\ldots,t_{k-1})}(1-\lambda)^{\des(t_{hp},\ldots,t_{k-1})}
\| \mathrm Z^{1,\Lie}_{[0,1]}(t_{ph+1})\|_{\ell^1}\ldots \| \mathrm Z^{1,\Lie}_{[0,1]}(t_{k-1})\|_{\ell^1}.
\end{align}
Let us introduce the notation
\begin{multline}
K_{p-1}^{(\lambda),\Lie}(t_0,t_p)=\\=\biggl\|\int_{\mathbf t_1=(t_1,\ldots,t_{p-1})\in[0,1]^k}
\lambda^{\asc(t_0,\mathbf t_1,t_p)}
(\lambda-1)^{\des(t_0,\mathbf t_1,t_p)}\mathrm Z^{1,\Lie}_{[0,1]}(t_1)\ldots \mathrm Z^{1,\Lie}_{[0,1]}(t_{p-1}) \biggr\|_{\ell^1}.
\plabel{eq:kernel}
\end{multline}

With this notation \eqref{eq:multdec} yields
\begin{equation}
\Theta_k^{(\lambda),\Lie}\leq\int_{t_0,t_{p},\ldots, t_{ph}\in[0,1]}
K_{p-1}^{(\lambda),\Lie}(t_0,t_p)\ldots K_{p-1}^{(\lambda),\Lie}(t_{p(h-1)},t_{ph})\mathrm dt_{0}\ldots \mathrm dt_{ph}
\plabel{eq:multdec2}
\end{equation}
(if we apply the trivial $\leq1$ estimate in the last line of \eqref{eq:multdec}).

\begin{commentx}

\begin{lemma}
\plabel{lem:submult}
For $p-1,q-1\geq0$,
\[K_{p+q-1}^{(\lambda),\Lie}(t_0,t_{p+q})=
\int_{t_p=0}^1 K_{p-1}^{(\lambda),\Lie}(t_0,t_{p}) K_{q-1}^{(\lambda),\Lie}(t_p,t_{p+q})\]
\begin{proof}
Applying the submultiplicative property of $|\cdot|_{\mathrm F\Lie}$, we find
\begin{align*}
&K_{p+q-1}^{(\lambda),\Lie}(t_0,t_{p+q})=
\\&=\biggl|\int_{\mathbf t=(t_1,\ldots,t_{p+q-1})\in[0,1]^k}
\lambda^{\asc(t_0,\mathbf t,t_{p+q})}
(\lambda-1)^{\des(t_0,\mathbf t,t_{p+q})}\mathrm Z^{1}_{[0,1]}(t_1)\ldots \mathrm Z^{1}_{[0,1]}(t_{p+q-1}) \biggr|_{\mathrm F\Lie}.
\\&\leq\int_{t_p=0}^1
\\&\quad
\biggl|\int_{\mathbf t'=(t_1,\ldots,t_{p-1})\in[0,1]^k}
\lambda^{\asc(t_0,\mathbf t',t_{p})}
(\lambda-1)^{\des(t_0,\mathbf t',t_{p})}\mathrm Z^{1}_{[0,1]}(t_1)\ldots \mathrm Z^{1}_{[0,1]}(t_{p-1}) \biggr|_{\mathrm F\Lie}
\\&\quad
|\mathrm Z^{1}_{[0,1]}(t_{p})|_{\mathrm F\Lie}
\\&\quad
\biggl|\int_{\mathbf t''=(t_{p+1},\ldots,t_{p+q})\in[0,1]^k}
\lambda^{\asc(t_{p},\mathbf t'',t_{p+q})}
(\lambda-1)^{\des(t_p,\mathbf t'',t_{p+q})}\mathrm Z^{1}_{[0,1]}(t_{p+1})\ldots \mathrm Z^{1}_{[0,1]}(t_{p+q-1}) \biggr|_{\mathrm F\Lie}
\\&\quad \mathrm dt_p
\\&
=\int_{t_p=0}^1 K_{p-1}^{(\lambda),\Lie}(t_0,t_{p}) K_{q-1}^{(\lambda),\Lie}(t_p,t_{p+q}).
\end{align*}
\end{proof}
\end{lemma}

\begin{cor}
\plabel{cor:submult}
For $p-1,q-1\geq0$,
\[\|K_{p+q-1}^{(\lambda),\Lie}\|_{L^2}\leq \|K_{p-1}^{(\lambda),\Lie}\|_{L^2}\cdot\|K_{q-1}^{(\lambda),\Lie}\|_{L^2}.\]
\begin{proof}
Note that for operators of nonnegative kernel the norm can be tested only for nonnegative functions.
Then the statement follows from the previous lemma.
\end{proof}
\end{cor}
\end{commentx}

Here $K_{p-1}^{(\lambda),\Lie} $ is nonnegative, and a trivial estimate is $K_{p-1}^{(\lambda),\Lie} \leq 1$.
For $p-1\geq1$, the function $K_{p-1}^{(\lambda),\Lie}(t_0,t_p)$ is continuous.
We will naturally consider these $K_{p-1}^{(\lambda),\Lie} $ as nonnegative integral kernels.

In general, if $K\in L^2([0,1]^{2})$, then we define the associated integral operator $I_K$
such that for   $f\in L^2([0,1])$ it yields
\[(I_Kf)(s)=\int_{r=0}^1 K(s,r)f(r)\, \mathrm dr.\]
The norm of $I_K$ in $L^2$ sense will be denoted by $\|I_K\|_{L^2}$.
It is well-known that $\|I_K\|_{L^2}\leq |K|_{L^2}$.
(In Part IA, the standard theory is reviewed in greater detail.)

Then \eqref{eq:multdec2} reads as
\[
\Theta_k^{(\lambda),\Lie}\leq \left\langle 1_{[0,1]}\left(I_{K_{p-1}^{(\lambda),\Lie}}\right)^{\lfloor\frac {k-1}p\rfloor},   1_{[0,1]}\right\rangle.
\]
By a trivial estimate,

\begin{equation}
\Theta_k^{(\lambda),\Lie}\leq \left\|I_{K_{p-1}^{(\lambda),\Lie}}\right\|_{L^2}^{\lfloor\frac {k-1}p\rfloor}.
\plabel{eq:multdec3}
\end{equation}

Recall from Part I, that for $\lambda\in[0,1]$, we have already considered the expressions
\[{\mathrm C}_\infty^{(\lambda)}=\begin{cases}
2&\text{if }\lambda=\frac12,
\\
\dfrac{2\artanh (1-2\lambda)}{1-2\lambda}=\dfrac{\log\dfrac{1-\lambda}{\lambda}}{1-2\lambda}&\text{if }\lambda\in(0,1)
\setminus\{\frac12\},\\
+\infty&\text{if }\lambda\in\{0,1\};
\end{cases}\]
\[w^{(\lambda)}=1/{\mathrm C}_\infty^{(\lambda)};\]

Here ${\mathrm C}_\infty^{(\lambda)}$ was the convergence radius of $\Theta^{\lambda}(x)$.
(Note that $w^{(\lambda)}$ is a concave, nonnegative function on $\lambda\in[0,1]$, symmetric for $\lambda\mapsto1-\lambda$;  its maximum is $w(1/2)=1/2$.
In more concrete terms, $1/w(\lambda)$ is just the convergence radius of $G(\lambda x,(1-\lambda)x)$ in $x$, cf. Part I.)
In the same spirit,
let ${\mathrm C}_\infty^{(\lambda),\Lie}$ be the convergence radius of
$\Theta^{(\lambda),\Lie}(x)$, and let
$w^{(\lambda),\Lie}=1/{\mathrm C}_\infty^{(\lambda),\Lie}$.
Moreover, for $p-1\geq0$, we set
\[w^{(\lambda),\Lie}_{p-1}= \mathrm r\left(I_{K_{p-1}^{(\lambda),\Lie}}\right),\]
i.~e.~as the spectral radius of the given integral operator.
We set
${\mathrm C}_{\infty,p-1}^{(\lambda),\Lie}=1/w^{(\lambda),\Lie}_{p-1}$.

Previously we have defined $K_{p-1}^{(\lambda),\Lie}(t_0,t_p)$ with respect to the norm $\|\cdot\|_{\ell^1}$.
The same can also be done using $|\cdot|_{\ell^1}$, leading to the non-Lie version $K_{p-1}^{(\lambda)}(t_0,t_p)$.
The discussion is analogous but computationally much simpler, because $|\cdot|_{\ell^1}$ is just the sum of
the absolute values of the coefficients.
But then one can see that
\[I_{K_{p-1}^{(\lambda)}}=\left(I_{K_{0}^{(\lambda)}}\right)^p.\]
From that it is also easy to see that
\[\mathrm r\left(I_{K_{p-1}^{(\lambda)}}\right)=\mathrm r\left(I_{K_{0}^{(\lambda)}}\right)^p=\left(w^{(\lambda)}\right)^p.\]
(Indeed, as $K_0^{(\lambda)}$ is bounded and nonnegative, it is particularly easy to see that
in order to compute the spectral radius,
it is sufficient to consider $\langle1_{[0,1]},(I_{K_0^{(\lambda)}})^{k-2}1_{[0,1]}\rangle$, which reproduces $\Theta_{k}^{(\lambda)}$.)

\begin{theorem}
\plabel{th:pst}
For $p-1\geq0$,

\[w^{(\lambda),\Lie}\leq \sqrt[p]{w_{p-1}^{(\lambda),\Lie}}\leq w^{(\lambda)}.\]
Equivalently,
\[{\mathrm C}_\infty^{(\lambda)}\leq \sqrt[p]{{\mathrm C}_{\infty,{p-1}}^{(\lambda),\Lie}}\leq
 {\mathrm C}_\infty^{(\lambda),\Lie}.\]

\begin{proof}
It is sufficient to prove only the first set of inequalities.
The very first inequality can be proven by taking $\sqrt[k]{\cdot}$ of \eqref{eq:multdec3} and $\limsup$.
The second inequality follows from the fact that Lie case is dominated by plain algebraic case,
yet the in plain algebraic case the spectral properties are straightforward.
\end{proof}
\end{theorem}

It is easy to see that the kernels $K_{p-1}^{(\lambda)}$ are continuous in $\lambda$ (in $L^2$ sense).
Thus they yield a continuous family of compact operators, whose spectral radius is continuous.
(In general for a continuous family $\lambda\mapsto A_\lambda$ of bounded operators we case state only
$\liminf_{\lambda\rightarrow\lambda_0} \mathrm r(A_\lambda)\geq \mathrm r(A_{\lambda_0})$
i. e. lower semicontinuity.)
In particular, $w_{p-1}^{(\lambda)}$ is continuous function a $[0,1/2]$-valued function.
(The range of the majorizing $w^{\lambda}$ is $[0,\frac12]$.)
Consequently,
$\lambda\mapsto {\mathrm C}_{\infty,p-1}^{(\lambda),\Lie}$ is continuous as a $[2,+\infty]$-valued function.
We also have

\begin{lemma}\plabel{lem:Liecont}
For $\lambda_1,\lambda_2\in(0,1)$,
\[\left|\mathrm C^{(\lambda_1),\Lie}_\infty-\mathrm C^{(\lambda_2),\Lie}_\infty \right|\leq
\left| \log\frac{\lambda_1}{1-\lambda_1}  - \log\frac{\lambda_2}{1-\lambda_2} \right|\]
holds.
\begin{proof}
This follows from the formula
\begin{multline*}
\mathcal R^{(\lambda_2)}(\Rexp ( (t\cdot \mathrm Z^{\Lie}_{[0,1)})   )=\\=\frac{\lambda_1(1-\lambda_1)}{\lambda_2(1-\lambda_2)}
\mathcal R^{(\lambda_1)}\left(\Rexp \left( (t\cdot \mathrm Z^{\Lie}_{[0,1)})\boldsymbol.
\left(\log\frac{\lambda_1}{1-\lambda_1}  - \log\frac{\lambda_2}{1-\lambda_2}\right)
\mathbf 1_{[1,2)}  \right) \right)+\frac{\lambda_2-\lambda_1}{\lambda_2(1-\lambda_2)}.
\end{multline*}
This shows that one can obtain one resolvent from the other with an indicated tolerance in the cumulative norm.
\end{proof}
\end{lemma}
\begin{theorem}
\plabel{th:Liecont}
$\lambda\mapsto w^{(\lambda),\Lie}$ is continuous as a $[0,1/2]$-valued function;
$\lambda\mapsto {\mathrm C}_\infty^{(\lambda),\Lie}$ is continuous as a $[2,+\infty]$-valued function.
\begin{proof}
This is an immediate consequence of the previous lemma.
\end{proof}
\end{theorem}

Let
\[w^{(\log),\Lie}=\max_{\lambda\in[0,1]}  w^{(\lambda),\Lie},\]
and
\begin{equation}
{\mathrm C}_\infty^{(\log),\Lie}=\min_{\lambda\in[0,1]}  {\mathrm C}_\infty^{(\lambda),\Lie}.
\plabel{eq:Clogdef}
\end{equation}
Here $w^{(\log),\Lie} =1/{\mathrm C}_\infty^{(\log),\Lie}$ holds.

Recall that ${\mathrm C}_\infty^{\Lie}$ is the convergence radius of $\Theta^{\Lie}(x)$.
Let $w^{\Lie}=1/{\mathrm C}_\infty^{\Lie}$.

\begin{lemma}
\plabel{lem:mmm}
\[ {\mathrm C}_\infty^{(\log),\Lie}\leq {\mathrm C}_\infty^{\Lie}.\]
\begin{proof}
$(\lambda, t)\mapsto \lambda+(1-\lambda)\Rexp(t\cdot\mathrm Z^{\Lie}_{[0,1]}) $ is analytic and invertible
on  $[0,1]_\lambda\times \intD(0,[{\mathrm C}_\infty^{(\log),\Lie})$,
thus the resolvent expression is also analytic.
By the Cauchy formula
\[\frac{f^{(k)}(0)}{k!}=\frac1{2\pi\mathrm i}\int_{z\in \partial\Dbar(0,r)^\circlearrowleft}\frac{f(z)}{(z-0)^k}\,\mathrm dz,\]
we have some uniform estimates (independently from $\lambda$) for the
coefficients in $t$, which can be integrated in $\lambda$.
\end{proof}
\end{lemma}

\begin{commentx}
Note, that ${\mathrm C}_\infty^{(\log),\Lie}$ has more meaning than a simple numerical value set up by \eqref{eq:Clogdef}.
It is exactly threshold value which guarantees the existence of $\mu_{\mathrm R}(s\cdot Z_{[0,1)}^{\Lie})$ realized
as $\log(\Rexp(s\cdot Z_{[0,1)}^{\Lie}))$ (correctly, as analytical continuation shows).
Thus it is convergence radius of the Lie algebraic Magnus expansion in (the stronger) logarithmic sense.
\end{commentx}

Similarly as before, let
\[w^{(\log),\Lie}_{p-1}=\max_{\lambda\in[0,1]}  w^{(\lambda),\Lie}_{p-1},\]
and
\begin{equation}
{\mathrm C}_{\infty,p-1}^{(\log),\Lie}=\min_{\lambda\in[0,1]}  {\mathrm C}_{\infty,p-1}^{(\lambda),\Lie}.
\plabel{eq:Clogdefp}
\end{equation}
Here $w^{(\log),\Lie}_{p-1} =1/{\mathrm C}_{\infty,p-1}^{(\log),\Lie}$ holds.

\begin{theorem}
\plabel{th:Lieconv}
\[ w^{\Lie}\leq w^{(\log),\Lie}\leq \sqrt[p]{w^{(\log),\Lie}_{p-1}}\leq \frac12.\]

Or, equivalently,
\[2\leq \sqrt[p]{{\mathrm C}_{\infty,p-1}^{(\log),\Lie}}
\leq {\mathrm C}_\infty^{(\log),\Lie}\leq {\mathrm C}_\infty^{\Lie}.\]
\begin{proof}
This is just some of the previous information put together.
\end{proof}
\end{theorem}

Our general strategy is that if we obtain an upper estimate
$w^{(\log),\Lie}\leq C$ (via estimating $\sqrt[p]{w^{(\log),\Lie}_{p-1}}$), then it yields a lower estimate
$\frac1C\leq{\mathrm C}_\infty^{(\log),\Lie}\leq{\mathrm C}_\infty^{\Lie}$.

\snewpage
\subsection{The structure of the kernels}
\plabel{ss:kerstruct}
~\\

Let us take a closer look at $K_{p-1}^{(\lambda),\Lie}(t_0,t_p)$.
Assume that $t_0<t_p$.
In \eqref{eq:kernel}, the integrand is best to be decomposed according to the distribution of
$\{t_1,\ldots,t_{p-1}\}$ relative to $t_0,t_p$.
Here we mean imagine $a$ to be the number of indices smaller than $t_0$ and $t_p$;
$b$ be the number of indices between $t_0$ and $t_p$;
$c$ be the number of indices greater than $t_0$ and $t_p$.
For $a+b+c=p-1$, let
\[p_{a,b,c}(t_0,t_p)=\frac{(p-1)!}{a!b!c!}t_0^a(t_p-t_0)^b(1-t_p)^c;\]
and
\[\mu_{a,b,c}^{(\lambda)}(X_1,\ldots,X_{p-1})=\sum_{\sigma\in\Sigma_{p-1}}
\lambda^{\asc(a+\frac12,\sigma,p-\frac12-c)}(\lambda-1)^{\des(a+\frac12,\sigma,p-\frac12-c)}
 X_{\sigma(1)}\ldots X_{\sigma(p-1)};\]
and
\[\Theta_{a,b,c}^{(\lambda),\Lie}=\frac1{(p-1)!}\left\| \mu_{a,b,c}^{(\lambda)}(Y_1,\ldots,Y_{p-1})\right\|_{\ell^1}.\]
Then
\begin{equation}
K^{(\lambda),\Lie}_{p-1}(t_0,t_p)=\sum_{a+b+c=p-1}p_{a,b,c}(t_0,t_p)\Theta_{a,b,c}^{(\lambda),\Lie}.
\plabel{eq:Kdeco}
\end{equation}
Here $p_{a,b,c}(t_0,t_p)$ refers to the probability of the configuration, and $\Theta_{a,b,c}^{\Lie}$
is the contribution of the corresponding noncommutative term.

There is a similar analysis for $t_0>t_p$.
Let
\[\tilde \mu_{a,b,c}^{(\lambda)}(X_1,\ldots,X_{p-1})=\sum_{\sigma\in\Sigma_{p-1}}
\lambda^{\asc(p-\frac12-c,\sigma,a+\frac12)}(\lambda-1)^{\des(p-\frac12-c,\sigma,a+\frac12)}
 X_{\sigma(1)}\ldots X_{\sigma(p-1)};\]
and
\[\tilde\Theta_{a,b,c}^{(\lambda),\Lie}=\frac1{(p-1)!}\left\| \tilde\mu_{a,b,c}^{(\lambda)}(Y_1,\ldots,Y_{p-1})\right\|_{\ell^1}.\]
Then
\begin{equation}
K^{(\lambda),\Lie}_{p-1}(t_0,t_p)=\sum_{a+b+c=p-1}p_{a,b,c}(t_p,t_0)\tilde\Theta_{a,b,c}^{(\lambda),\Lie}.
\plabel{eq:KdecoAnt}
\end{equation}
Similarly to \eqref{eq:poref},
\[\tilde \mu_{a,b,c}^{(\lambda)}(X_1,\ldots,X_{p-1})=-\mu_{c,b,a}^{(1-\lambda)}(-X_{p-1},\ldots,-X_{1})\]
holds in general; it implies
\[\tilde \Theta_{a,b,c}^{(\lambda),\Lie}=\Theta_{c,b,a}^{(1-\lambda),\Lie}.\]
Also,
\[p_{a,b,c}(t_0,t_p)=p_{c,b,a}(1-t_p,1-t_0)\]
holds. Thus
\begin{equation}
K^{(\lambda),\Lie}_{p-1}(t_0,t_p)=\sum_{a+b+c=p-1}p_{c,b,a}(1-t_0,1-t_p)\Theta_{c,b,a}^{(1-\lambda),\Lie}.
\plabel{eq:KdecoAntis}
\end{equation}
Therefore,
\[K^{(\lambda),\Lie}_{p-1}(t_0,t_p)=K^{(1-\lambda),\Lie}_{p-1}(1-t_0 ,1-t_p)\]
holds generally.
(Remark: From the forthcoming discussion
\[K^{(\lambda),\Lie}_{p-1}(t_0,t_p)=K^{(1-\lambda),\Lie}_{p-1}(t_p ,t_0)\]
will also be clear.)

Restricted to $t_0<t_p$ or $t_0>t_p$  , $K^{(\lambda),\Lie}_{p-1}(t_0,t_p)$ is a sum of products of terms polynomial in $t_0,t_p$
and terms piecewise polynomial in $\lambda$.
Indeed, $\left\| \mu_{a,b,c}^{(\lambda)}(Y_1,\ldots,Y_{p-1})\right\|_{\ell^1}$ is the maximum of the linear
functions defining the hyperfaces of the unit ball of $\|\cdot\|_{\ell^1}$ evaluated on $\mu_{a,b,c}^{(\lambda)}(Y_1,\ldots,Y_{p-1})$.

Now, one can greatly simplify \eqref{eq:Kdeco} and \eqref{eq:KdecoAnt}/\eqref{eq:KdecoAntis}.
The idea is to use
\begin{lemma}\plabel{lem:csere} For $a+b+c+1=p-1$,
\[\mu_{a+1,b,c}^{(\lambda)}(X_1,\ldots,X_{p-1})=\mu_{a,b,c+1}^{(\lambda)}(X_2,\ldots,X_{p-1},X_1).\]
\begin{proof}
If we rename the lowest position to the highest position, then it also yields one descent and one ascent, while
the descent/ascent relations between other indices remain the same.
\end{proof}
\end{lemma}
\begin{cor}
$\Theta_{a,b,c}^{(\lambda)}$ depends only on $\lambda$, $a+c$, and $b$.
\qed
\end{cor}
This suggests that we can do this presentation with only two indices.
Let us define
\[p_{a,b}(t)=\frac{(p-1)!}{a!b!}(1-t)^at^b;\]
and
\[\mu_{a,b}^{(\lambda)}(X_1,\ldots,X_{p-1})=\sum_{\sigma\in\Sigma_{p-1}}
\lambda^{\asc(a+\frac12,\sigma)}(\lambda-1)^{\des(a+\frac12,\sigma)}
 X_{\sigma(1)}\ldots X_{\sigma(p-1)}\]
making $\mu_{a,b,0}^{(\lambda)}(X_1,\ldots,X_{p-1}) =\lambda\cdot\mu_{a,b}^{(\lambda)}(X_1,\ldots,X_{p-1}) $; and
\begin{equation}
\Theta_{a,b}^{(\lambda),\Lie}=\frac1{(p-1)!}\left\| \mu_{a,b}^{(\lambda)}(Y_1,\ldots,Y_{p-1})\right\|_{\ell^1}.
\plabel{eq:ThetaLieDef}
\end{equation}
making $\Theta^{(\lambda),\Lie}_{a,b,c}=\lambda \Theta^{(\lambda),\Lie}_{c+a,b}$ and
$\tilde\Theta^{(\lambda),\Lie}_{a,b,c}=(1-\lambda) \Theta^{(1-\lambda),\Lie}_{c+a,b}$.

\begin{theorem}
For $t_0<t_p$,
\begin{equation}
K^{(\lambda),\Lie}_{p-1}(t_0,t_p)=\lambda\cdot\sum_{a+b=p-1}p_{a,b}(t_p-t_0)\Theta_{a,b}^{(\lambda),\Lie}.
\plabel{eq:Kcoup}
\end{equation}
For $t_0>t_p$,
\begin{equation}
K^{(\lambda),\Lie}_{p-1}(t_0,t_p)=(1-\lambda)\cdot\sum_{a+b=p-1}p_{a,b}(t_0-t_p)\Theta_{a,b}^{(1-\lambda),\Lie}.
\plabel{eq:Kdecoup}
\end{equation}
\begin{proof}
This is just \eqref{eq:Kdeco} and \eqref{eq:KdecoAnt}/\eqref{eq:KdecoAntis} transcribed.
\end{proof}
\end{theorem}

\begin{cor}
For a fixed p,
$K^{(\lambda),\Lie}_{p-1}(t_0,t_p)$ depends only on $\lambda$ and $t_p-t_0$.
\qed
\end{cor}
Thus the notation $K^{(\lambda),\Lie}_{p-1}(t_0,t_p)\equiv K^{(\lambda),\Lie}_{p-1}(t_p-t_0)$ is reasonable.
Furthermore, \[K^{(\lambda),\Lie}_{p-1}(t)=K^{(1-\lambda),\Lie}_{p-1}(-t) .\]

\begin{lemma}\plabel{lem:ThetaLT}

\[\mu_{a,b}^{(\lambda)}(X_1,\ldots,X_{p-1})
=\mu_{b,a}^{(1-\lambda)}(-X_{p-1},\ldots,-X_{1})
.
\plabel{eq:door4}\]
\begin{proof}
This a consequence of the previous identities combined.
\end{proof}
\end{lemma}
\begin{cor}\plabel{cor:ThetaLT}
\[\Theta_{a,b}^{(1-\lambda),\Lie}=\Theta_{b,a}^{(\lambda),\Lie}.
\eqed\]
\end{cor}
It is easy to see that
\[p_{a,b}(1-t)=p_{b,a}(t).\]
Thus, another consequence of Lemma \ref{lem:ThetaLT} is that for $t_0>t_p$,
\begin{equation}
K^{(\lambda),\Lie}_{p-1}(t_0,t_p)=(1-\lambda)\cdot\sum_{a+b=p-1}p_{b,a}(1+t_p-t_0)\Theta_{b,a}^{(\lambda),\Lie}.
\plabel{eq:Kdecoupvar}
\end{equation}

\begin{commentx}
Let us define the reduced kernel
\begin{equation}
\widetilde K_{p-1}^{(\lambda),\Lie}(t)=
\begin{cases}
\sum_{a+b=p-1}p_{a,b}(t)\Theta_{a,b}^{(\lambda),\Lie}&\text{if}\quad t\in(0,1]\\
\sum_{a+b=p-1}p_{a,b}(t+1)\Theta_{a,b}^{(\lambda),\Lie}&\text{if}\quad t\in[-1,0),
\end{cases}
\plabel{eq:Kred}
\end{equation}
for $\lambda\in[0,1]$, $t\in[-1,1]$.
Then, for $t\in[-1,1]\setminus\{0\}$,
\[K_{p-1}^{(\lambda),\Lie}(t)=\lambda^{\asc(0,t)}(1-\lambda)^{\des(0,t)}   \widetilde K_{p-1}^{(\lambda),\Lie}(t).\]
\end{commentx}

\snewpage
\subsection{Estimates in low orders }
\plabel{ss:estlow}
~\\

Recall, by Lemma \ref{lem:ThetaLT}, we have to compute $\Theta_{a,p-1-a}^{(\lambda),\Lie}$ only for $0\leq a\leq\lfloor\frac{p-1}2\rfloor$.

\begin{example}
For $p-1=1$,
\[\Theta_{0,1}^{(\lambda),\Lie}=\lambda\quad\lambda\in[0,1].\]
\qedexer
\end{example}
\begin{example}
For $p-1=2$,
\[\Theta_{0,2}^{(\lambda),\Lie}=\begin{cases}
-{\lambda}^{2}+\lambda&\text{if }\lambda\in\left[0,\frac12\right]\\
{\lambda}^{2}&\text{if }\lambda\in\left[\frac12,1\right].
\end{cases}\]
\[\Theta_{1,1}^{(\lambda),\Lie}=2\lambda(1-\lambda)\quad\text{if }\lambda\in[0,1].\]
\qedexer
\end{example}
\begin{example}
For $p-1=3$,
\[\Theta_{0,3}^{(\lambda),\Lie}=
\begin{cases}
-2\,{\lambda}^{3}+\lambda&\text{if }\lambda\in\left[0,\frac13\right]\\
-5\,{\lambda}^{3}+4\,{\lambda}^{2}&\text{if }\lambda\in\left[\frac13,\frac12\right]\\
-5\,{\lambda}^{3}+6\,{\lambda}^{2}-\lambda&\text{if }\lambda\in\left[\frac12,\frac23\right]\\
-2\,{\lambda}^{3}+4\,{\lambda}^{2}-\lambda&\text{if }\lambda\in\left[\frac23,1\right];
\end{cases}
\]
\[\Theta_{1,2}^{(\lambda),\Lie}=
\begin{cases}
2\,{\lambda}^{3}-4\,{\lambda}^{2}+2\,\lambda&\text{if }\lambda\in\left[0,\frac13\right]\\
-4\,{\lambda}^{3}+4\,{\lambda}^{2}&\text{if }\lambda\in\left[\frac13,1\right].
\end{cases}
\]
\qedexer
\end{example}
\begin{example}
For $p-1=4$,
\[\Theta_{0,4}^{(\lambda),\Lie}=
\begin{cases}
12\,{\lambda}^{4}-19\,{\lambda}^{3}+5\,{\lambda}^{2}+\lambda&\text{if }\lambda\in\left[0,\frac{5-\sqrt{17}}4\right]\\
14\,{\lambda}^{4}-26\,{\lambda}^{3}+11\,{\lambda}^{2}&\text{if }\lambda\in\left[\frac{5-\sqrt{17}}4,\frac{\sqrt3-1}2\right]\\
12\,{\lambda}^{4}-26\,{\lambda}^{3}+14\,{\lambda}^{2}-\lambda&\text{if }\lambda\in\left[\frac{\sqrt3-1}2,\frac12\right]\\
-12\,{\lambda}^{4}+10\,{\lambda}^{3}+2\,{\lambda}^{2}-\lambda&\text{if }\lambda\in\left[\frac12,\frac{3-\sqrt3}2\right]\\
-14\,{\lambda}^{4}+16\,{\lambda}^{3}-{\lambda}^{2}-\lambda&\text{if }\lambda\in\left[\frac{3-\sqrt3}2,\frac{\sqrt{17}-1}4\right]\\
-12\,{\lambda}^{4}+17\,{\lambda}^{3}-3\,{\lambda}^{2}-\lambda&\text{if }\lambda\in\left[\frac{\sqrt{17}-1}4,1\right];
\end{cases}
\]
\[\Theta_{1,3}^{(\lambda),\Lie}=
\begin{cases}
12\,{\lambda}^{4}-15\,{\lambda}^{3}+{\lambda}^{2}+2\,\lambda&\text{if }\lambda\in\left[0,\frac{3-\sqrt3}6\right]\\
18\,{\lambda}^{4}-27\,{\lambda}^{3}+8\,{\lambda}^{2}+\lambda&\text{if }\lambda\in\left[\frac{3-\sqrt3}6,\frac14\right]\\
22\,{\lambda}^{4}-36\,{\lambda}^{3}+14\,{\lambda}^{2}&\text{if }\lambda\in\left[\frac14,\frac{\sqrt3-1}2\right]\\
18\,{\lambda}^{4}-36\,{\lambda}^{3}+20\,{\lambda}^{2}-2\,\lambda&\text{if }\lambda\in\left[\frac{\sqrt3-1}2,\frac12\right]\\
2\,{\lambda}^{4}-12\,{\lambda}^{3}+12\,{\lambda}^{2}-2\,\lambda&\text{if }\lambda\in\left[\frac12,\frac23\right]\\
-{\lambda}^{4}-7\,{\lambda}^{3}+10\,{\lambda}^{2}-2\,\lambda&\text{if }\lambda\in\left[\frac23,\frac34\right]\\
-5\,{\lambda}^{4}+7\,{\lambda}^{2}-2\,\lambda&\text{if }\lambda\in\left[\frac34,1\right];
\end{cases}
\]

\[\Theta_{2,2}^{(\lambda),\Lie}=
\begin{cases}
3\,{\lambda}^{4}+2\,{\lambda}^{3}-9\,{\lambda}^{2}+4\,\lambda&\text{if }\lambda\in\left[0,\frac29\right]\\
\frac{15}2\,{\lambda}^{4}-8\,{\lambda}^{3}-\frac52\,{\lambda}^{2}+3\,\lambda&\text{if }\lambda\in\left[\frac29,\frac4{17}\right]\\
16\,{\lambda}^{4}-27\,{\lambda}^{3}+10\,{\lambda}^{2}+\lambda&\text{if }\lambda\in\left[\frac4{17},\frac14\right]\\
20\,{\lambda}^{4}-36\,{\lambda}^{3}+16\,{\lambda}^{2}&\text{if }\lambda\in\left[\frac14,\frac{\sqrt5}5\right]\\
{\frac {35}{2}}\,{\lambda}^{4}-{\frac {67}{2}}\,{\lambda}^{3}+{\frac {33}{2}}\,{\lambda}^{2}-\frac12\,\lambda&\text{if }\lambda\in\left[\frac{\sqrt5}5,\frac{6-\sqrt{21}}3\right]\\
18\,{\lambda}^{4}-36\,{\lambda}^{3}+{\frac {58}{3}}\,{\lambda}^{2}- \frac43\,\lambda&\text{if }\lambda\in\left[\frac{6-\sqrt{21}}3,\frac{\sqrt{21}-3}3\right]\\
{\frac {35}{2}}\,{\lambda}^{4}-{\frac {73}{2}}\,{\lambda}^{3}+21\,{\lambda}^{2}-2\,\lambda&\text{if }\lambda\in\left[\frac{\sqrt{21}-3}3,\frac{5-\sqrt5}5\right]\\
20\,{\lambda}^{4}-44\,{\lambda}^{3}+28\,{\lambda}^{2}-4\,\lambda&\text{if }\lambda\in\left[\frac{5-\sqrt5}5,\frac34\right]\\
16\,{\lambda}^{4}-37\,{\lambda}^{3}+25\,{\lambda}^{2}-4\,\lambda&\text{if }\lambda\in\left[\frac34,\frac{13}{17}\right]\\
\frac{15}2\,{\lambda}^{4}-22\,{\lambda}^{3}+{\frac {{37}}{2}\,\lambda}^{2}-4\,\lambda&\text{if }\lambda\in\left[\frac{13}{17},\frac79\right]\\
3\,{\lambda}^{4}-14\,{\lambda}^{3}+15\,{\lambda}^{2}-4\,\lambda&\text{if }\lambda\in\left[\frac79,1\right].
\end{cases}
\]
\qedexer
\end{example}

\begin{remark}
\plabel{rem:algor}
Finding $\Theta_{a,b}^{(\lambda),\Lie}$ (as a maximum of finitely many polynomials) can be algorithmized very easily.
However, in practice, it is a computationally very intensive problem.
Fortunately, we do not have to find the exact expressions for good estimates, as we will see.
\qedremark
\end{remark}
Having $\Theta_{a,b}^{(\lambda),\Lie}$ with $p-1=a+b$, we can write down $K_{p-1}^{(\lambda),\Lie}$ explicitly.
The next question is how to locate $\mathrm r\left(I_{ K_{p-1}^{(\lambda),\Lie}}\right)$ effectively.
Here it is a major advantage that we deal with operators of nonnegative kernels.
In general, rectangularly based non-negative step-functions on $[0,1]^2$ realize
classical Perron--Frobenius theory, and these step-functions are dense among nonnegative kernels
in $L^2$ sense.
Consequently, there are several statements extending from classical Perron--Frobenius
theory to nonnegative operators without much trouble.
As the kernels $K_{p-1}^{(\lambda),\Lie}$ are continuous for $p-1\geq1$ anyway;
instead of referring to several articles, we will quote only Anselone, Lee \cite{AL}
 in order to get an impression about the matters.
(A more detailed discussion can be found in Part IA \cite{L15}.)
We only quote the seemingly innocent
\begin{theorem}
\plabel{th:average}
(Averaging principle, special case.)
Assume that $K\geq0$.
For $n\in\mathbb N$,
\[n\mapsto\left[\mathrm{ess\,inf\,}\frac{(I_K)^{n+1}1_{[0,1]}}{(I_K)^n1_{[0,1]}},
\mathrm{ess\,sup\,}\frac{(I_K)^{n+1}1_{[0,1]}}{(I_K)^n1_{[0,1]}}\right]
\]
yields a sequence of encapsulated intervals (all) containing ${\mathrm r}(I_K)$.

(Here $\frac00=$``undecided''; if the quotient is $\frac00$ almost everywhere, i.~e.~if $(I_K)^n1_{[0,1]}=0$ is reached, then we set the interval to be $[0,0]$.)
\begin{proof}
See Anselone, Lee \cite{AL} (Theorem 7.1) for the (critical) containment of spectral radius in the continuous case;
or see the discussion in Part IA \cite{L15}.
\end{proof}
\end{theorem}

Applying the previous theorem, we can localize $\sqrt[p]{w^{(\lambda),\Lie}_{p-1}}=\sqrt[p]{\mathrm r\left(I_{ K_{p-1}^{(\lambda),\Lie}}\right)}$
in order to obtain an upper estimate $w^{(\log),\Lie}\leq C$.
Then, according to Theorem \ref{th:Lieconv},  that yields a lower estimate
$\frac1C\leq{\mathrm C}_\infty^{(\log),\Lie}\leq{\mathrm C}_\infty^{\Lie}$.
\begin{example}

(a) In case of $p-1=2$, it yields
\[\mathrm C^{\Lie}_\infty>2.250\quad.\]

(b) In case of $p-1=3$, it yields
\[\mathrm C^{\Lie}_\infty>2.282\quad.\]

(c) In case of $p-1=4$, it yields
\[\mathrm C^{\Lie}_\infty>2.364\quad.\eqedexer\]
\end{example}

\begin{commentx}
In general, one does not need  the precise function $\Theta_{a,b}^{(\lambda),\Lie}$
in order to have good estimates.
We need only good upper estimates for them.
For a large set $\Lambda$ of relevant $\lambda$, using linear programming, we find (hopefully) near-optimal
presentations for computing \eqref{eq:ThetaLieDef}.
One can round up the coefficients to rational numbers in order to find approximate presentations for $\lambda\in\Lambda$.
Then, using $\|\cdot\|_{\ell^1}\leq |\cdot|_{\ell^1}$, we can give upper estimates for $\Theta_{a,b}^{(\lambda),\Lie}$, $\lambda\in\Lambda$.
Again using  $\|\cdot\|_{\ell^1}\leq |\cdot|_{\ell^1}$, the modulus of continuity in $\lambda$ can be estimated well, leading to
upper estimates for  $\Theta_{a,b}^{(\lambda),\Lie}$, $\lambda\in[0,1]$.
In the final step, when we estimate the spectral radii, we have to find only greatest real eigenvalue according to Frobenius--Perron theory
(cf. Gantmacher \cite{Ga}).

\end{commentx}

For $p-1\geq 5$ the $\Theta^{(\lambda),\Lie}_{a,b}$ (with $a,b=p-1$) become cumbersome to write down in general.
All is not lost, however.
For any concrete $\lambda\in[0,1]$, the value $\Theta^{(\lambda),\Lie}_{a,b}$ can be still computed.
Even less, it is sufficient to get an upper estimate for $\Theta^{(\lambda),\Lie}_{a,b}$
from getting any concrete presentation of $\mu_{a,b}^{(\lambda)}(Y_1,\ldots,Y_{a+b})$
as a linear combination  products of Lie monomials.
(Thus we do not have to do the exact linear programming problem, it is
sufficient to obtain a merely acceptable approximate optimization and correct it to a non-optimal but
exact presentation by plain monomials.)
Then we get an upper estimate for $K_{p-1}^{(\lambda),\Lie}$.
As the spectral radius is monotone for nonnegative kernels, this leads to
an upper estimate for $w_{p-1}^{(\lambda),\Lie}$, hence for $w^{(\lambda,\Lie)}$.
Then we have a lower estimate for $\mathrm C^{(\lambda),\Lie}_\infty$.
Although this is for only one concrete $\lambda$, it still has global consequences by Lemma \ref{lem:Liecont}.
Moreover we do not have particulary concerned for $\lambda\sim0$ or $\lambda\sim1$ as we know
$\mathrm C^{(\lambda),\Lie}_\infty\geq \mathrm C^{(\lambda)}_\infty$ anyway.
Ultimately by a well-adapted choice of testing $\lambda$ values, we can obtain a good estimate for $\mathrm C^{(\log),\Lie}_\infty$,
hence for $\mathrm C^{\Lie}_\infty$.

\begin{example} In the above manner:

(a) In case of $p-1=5$, it yields
\[\mathrm C^{\Lie}_\infty>2.396\quad.\]

(b) In case of $p-1=6$, it yields
\[\mathrm C^{\Lie}_\infty>2.427\quad.\eqedexer\]
\end{example}
\begin{remark}
(a)
Naively, one would expect to get around $\mathrm C^{\Lie}_\infty>2.5$ in this way, having infinite resources;
which is below of our naive expectations, cf. Section \ref{sec:numerics}.

(b) The resolvent method as presented here is certainly not able to show $\mathrm C^{(\log),\Lie}_\infty>\pi$.
Indeed, replacing $\mathrm Z^{1,\Lie}_{[0,1)}$ by the Lebesgue measure $\mathbf 1_{[0,1)}$, we find
$\mathrm C^{(\lambda),\Lie}_\infty\leq \mathrm C^{(\lambda),\boldsymbol\varepsilon}$
(in the ternology of Part I).
In particular, $\mathrm C^{(1/2),\Lie}_\infty\leq \pi$.
\qedremark
\end{remark}
\snewpage

\section{The resolvent method and the BCH expansion}\plabel{sec:resbch}
The methods of the previous section also apply to the BCH expansion.
In that case, due to the heterogeneity of the domain of kernels, it is better to write the kernels as $2\times2$ matrices.
However, we will not follow that path here.
Instead, we concentrate on giving a direct argument showing that the Banach algebraic estimate can be improved.

As a reminder, for $x_1,x_2\geq0$,
\[\Gamma^{\Lie}(x_1,x_2)=\left\| \BCH(x_1Y_1,x_2Y_2) \right\|_{\ell^1},\]
where $Y_1,Y_2$ are appropriate formal variables.
We are looking for $x_1,x_2$ such that $\Gamma^{\Lie}(x_1,x_2)<+\infty$.
For $\lambda\in[0,1]$, we also define
\[\Upsilon^{(\lambda)}(x_1Y_1,x_2Y_2)=\lambda(1-\lambda)\mathcal R^{(\lambda)}(\exp x_1Y_1)\mathcal R^{(\lambda)}(\exp x_2Y_2).\]

\begin{theorem}
Suppose that $0\leq x_1,x_2<\pi$. If for some  $n\geq1$,
\[\sup_{\lambda\in[0,1]}\sqrt[n]{\left\|\Upsilon^{(\lambda)}(x_1Y_1,x_2Y_2)^n\right\|_{\ell^1}}<1, \]
then
\[\Gamma^{\Lie}(x_1,x_2)<+\infty.\]

In particular, if for the $\|\cdot\|_{\ell^1}$-spectral radius
\[\sup_{\lambda\in[0,1]}\mathrm r_{\|\cdot\|_{\ell^1}}\left(\Upsilon^{(\lambda)}(x_1Y_1,x_2Y_2)^n\right)<1, \]
then the conclusion applies.
\begin{proof} According to  Part I,  (formally)
\begin{align}
\BCH(x_1Y_1,x_2Y_2)=&\int_{\lambda=0}^1\mathcal R^{(\lambda)}((\exp x_1Y_1)(\exp x_2Y_2))\,\mathrm d\lambda\notag\\
=&\int_{\lambda=0}^1(1-\Upsilon^{(\lambda)}(x_1Y_1,x_2Y_2)   )^{-1}\mathcal R^{(\lambda)}(\exp x_1Y_1)\notag\\
&+\mathcal R^{(\lambda)}(\exp x_2Y_2)(1-\Upsilon^{(\lambda)}(x_1Y_1,x_2Y_2)  )^{-1}\notag\\
&+\lambda \mathcal R^{(\lambda)}(\exp x_1Y_1)\mathcal R^{(\lambda)}(\exp x_2Y_2)(1-\Upsilon^{(\lambda)}(x_1Y_1,x_2Y_2)  )^{-1} \notag\\
&+(\lambda-1) \mathcal R^{(\lambda)}(\exp x_2Y_2)(1-\Upsilon^{(\lambda)}(x_1Y_1,x_2Y_2)  )^{-1}\mathcal R^{(\lambda)}(\exp x_1Y_1) \notag\\
&\mathrm d\lambda\notag,
\end{align}
completely well-defined in every $(Y_1,Y_2)$-grade.
According to our assumptions, the terms $\mathcal R^{(\lambda)}((\exp x_1Y_1)$ and $\mathcal R^{(\lambda)}((\exp x_2Y_2)$ are bounded in $\|\cdot\|_{\ell^1}$, and
so are the relevant Neumann  series.
\end{proof}
\end{theorem}
The statement also applies to the  the case of  $|\cdot|_{\ell^1}$ (cf. Part I),
except in that case there is no difference between the spectral radius and the norm of $\Upsilon^{(\lambda)}(x_1Y_1,x_2Y_2)$.
So, in Part I only the $|\cdot|_{\ell^1}$ norm was used.
We have demonstrated in Part I that on the domain $0\leq x_1+x_2\leq\mathrm  C_2=2.89847930\ldots$, $\lambda\in[0,1]$ the inequality
\[\left|\Upsilon^{(\lambda)}(x_1Y_1,x_2Y_2)\right|_{\ell^1}=\sqrt[2]{\left|\Upsilon^{(\lambda)}(x_1Y_1,x_2Y_2)^2\right|_{\ell^1}}\leq1, \]
holds; and in case of equality $x_1=x_2=\frac12 \mathrm  C_2$ and $0.35865<\min(\lambda,1-\lambda)<0.35866$.

Let us now compare $\left\|\Upsilon^{(\lambda)}(x_1Y_1,x_2Y_2)^2\right\|_{\ell^1}$ and $\left|\Upsilon^{(\lambda)}(x_1Y_1,x_2Y_2)^2\right|_{\ell^1}$.
The first one is less or equal than the second one, actually formally degree-wise (in $x_1$ and $x_2$ jointly).
Let us consider the the coefficients of $x_1^2x_2^4$, which are coming from $\Upsilon^{(\lambda)}(Y_1,Y_2)$.
After some computation, one finds the degree component to be
\begin{multline}\left(\Upsilon^{(\lambda)}(Y_1,Y_2)^2\right)_{\deg_{(Y_1,Y_2)}=(2,4) }=\lambda^2(1-\lambda)^2\\
\underbrace{\left(
\left( {\lambda}^{2}-\lambda+\frac14 \right) Y_{122122}
+\left( {\lambda}^{2}-\lambda+\frac16 \right) (Y_{121222}+Y_{122212})
\right)}_H.\plabel{eq:san}
\end{multline}
(Here  $Y_{121222}\equiv Y_1Y_2Y_1Y_2Y_2Y_2$, etc.)

Thus, the coefficient of $x_1^2x_2^4$ in $\left|\Upsilon^{(\lambda)}(Y_1,Y_2)\right|_{\ell^1}$
is (the sum of the absolute value of the coefficients of the monomials)
\[\upsilon_{2,4}^{2,(\lambda)}=\begin{cases}
\lambda^2(1-\lambda)^2\left(-3\lambda(1-\lambda)+{\frac {7}{12}}\right)
\\\qquad\qquad\qquad\text{if }\min(\lambda,1-\lambda)\in\left[0,\frac12-\frac{\sqrt3}6\right],
\text{ i. e. if }\lambda(1-\lambda)\in \left[0,\frac16 \right],\\
\lambda^2(1-\lambda)^2\left(\lambda(1-\lambda)-\frac1{12}\right)
\\\qquad\qquad\qquad\text{if }\min(\lambda,1-\lambda)\in\left[\frac12-\frac{\sqrt3}6,\frac12\right],\text{ i. e. if }
\lambda(1-\lambda)\in \left[\frac16,\frac14 \right].
\end{cases}\]

On the other hand, the coefficient of $x_1^2x_2^4$ in $\left\|\Upsilon^{(\lambda)}(Y_1,Y_2)\right\|_{\ell^1}$
is
\[\upsilon_{2,4}^{2,(\lambda),\Lie}=\begin{cases}
\lambda^2(1-\lambda)^2\left(-3\lambda(1-\lambda)+{\frac {7}{12}}\right)
\\\qquad\qquad\text{if }\min(\lambda,1-\lambda)\in\left[0,\frac12-\frac{\sqrt3}6\right],
\text{ i. e. if } \lambda(1-\lambda)\in \left[0,\frac16 \right],\\
\lambda^2(1-\lambda)^2\left(  -2\,\lambda\, \left( 1-\lambda \right) +{\frac {5}{12}} \right)
\\\qquad\qquad\text{if }\min(\lambda,1-\lambda)\in\left[\frac12-\frac{\sqrt3}6,\frac12-\frac{\sqrt2}6 \right],\text{ i. e. if }\lambda(1-\lambda)\in \left[\frac16,\frac7{36} \right],
\\\lambda^2(1-\lambda)^2\left( -\frac12\,\lambda\, \left( 1-\lambda\right) +\frac18 \right)
\\\qquad\qquad\text{if }\min(\lambda,1-\lambda)\in\left[\frac12-\frac{\sqrt2}6,\frac12-\frac{\sqrt5}{10}\right],\text{ i. e. if }  \lambda(1-\lambda)\in \left[\frac7{36},\frac15 \right],
\\\lambda^2(1-\lambda)^2\left({\frac {17 }{6}}\,\lambda\, \left( 1-\lambda \right)-{\frac {13}{24}}\right)
\\\qquad\qquad\text{if }\min(\lambda,1-\lambda)\in\left[ \frac12-\frac{\sqrt5}{10} ,\frac12\right],\text{ i. e. if }
\lambda(1-\lambda)\in \left[\frac15,\frac14 \right].
\end{cases}\]
(Cf.
\begin{align*}
H=&\left(-\lambda(1-\lambda)+\frac14 \right)Y_{122122}
+\left(-\lambda(1-\lambda)+\frac16 \right) (Y_{121222}+Y_{122212})\\
H=&\left( -3\,\lambda\, \left( 1-\lambda \right) +{\frac {7}{12}}
 \right) {Y}_{{122122}}+ \left( \lambda\,
 \left( 1-\lambda \right) -\frac16 \right) {Y}_{{12[[2,[1,2]]2}}\\
H=&\left( \frac32\,\lambda\, \left( 1-\lambda \right) -{\frac {7}{24}}
 \right)  \left( -{Y}_{{122221}}-{Y}_{{112222}}-{Y}_{{1[2,[2,[2,[1,2]]]]}}
 \right)\\
&+ \left( -5\,\lambda\, \left( 1-\lambda \right) +1 \right) {Y}_{{12[[2,[1,2]]2}}\\
H=&\left( -\frac16\,\lambda\, \left( 1-\lambda \right) +\frac1{24} \right)  \left(
-{Y}_{{122221}}-{Y}_{{112222}}-{Y}_{{1[2,[2,[2,[1,2]]]]}} \right)\\
&+ \left( \frac53\,
\lambda\, \left( 1-\lambda \right) -\frac13 \right)  \left( -{Y}_{{122212}}-{Y}_{{121222}}
 \right)\\
\end{align*}
for the minimal presentations of the critical part of \eqref{eq:san} in the four cases, respectively.
Minimality is attested by the linear functionals
\[Y^*_{112222}+Y^*_{121222}+Y^*_{122122}+Y^*_{122212}\]
\[-Y^*_{112222}+Y^*_{122122}+Y^*_{122212}\]
\[-Y^*_{112222}+\frac12Y^*_{122122}-Y^*_{122221}+\frac18Y^*_{212122}-\frac14Y^*_{212221}
\]
\[-Y^*_{112222}-Y^*_{121222}-\frac56\,Y^*_{122122}-Y^*_{122212}-
Y^*_{122221}-\frac16Y^*_{212212}-\frac14Y^*_{212221}
\]
respectively.)

Now,  near the critical place, $x_1=x_2=\frac12 \mathrm  C_2$ and $0.35865<\min(\lambda,1-\lambda)<0.35866$, apparently there is a
difference between
(the coefficients of $x_1^2x_2^4$ in)
$\left\|\Upsilon^{(\lambda)}(x_1Y_1,x_2Y_2)^2\right\|_{\ell^1}$ and $\left|\Upsilon^{(\lambda)}(x_1Y_1,x_2Y_2)^2\right|_{\ell^1}$.
So, the former one  will not reach up to $1$ as the latter one does.
Yet, both expressions are continuous on a larger domain.
Consequently, we know that the convergence domain will extend.

\begin{example}
If we apply the theorem above with $n=2$, correction gains counted up to expansion degree $9$, then we find
\[\mathrm C_{2}^{\Lie}>2.93\qquad.\eqedexer\]
\end{example}
Increasing $n$ is not cheap, however.
Estimates for higher $n$ are better to use a slightly finer accounting than simply correcting up to a few expansion degrees, but
 like, corrections to powers of corrected series, etc.

In unbalanced case of BCH expansion (with prescribed norm ratios), the result above in terms of cumulative norm is not optimal.
In that case one can obtain convergence results by (improved versions of) M\'erigot's method, or
 one can obtain result convergence and divergence results by considering, again, modifications and consequences of the
 Banach algebraic case.
As our interest is more the cumulative norm, we will not work out any cases here.
\snewpage
\section{Numerics and heuristics }\plabel{sec:numerics}
We have seen that for the cumulative convergence radius of the Magnus expansion $2.427<\mathrm C_\infty^{\Lie}\leq 4$.
So, one still wonders about the exact value of $\mathrm C_\infty^{\Lie}$, that is about the convergence radius of the positive power series $\Theta^{\Lie}(x)$ .
We know the coefficients $\Theta^{\Lie}_n$ for  $n\leq9$; cf.  Appendix \ref{sec:someTheta}.
The value of a few coefficients is, of course, not conclusive about the convergence of a power series in any way.
(Although the exact values can be used in arguments improving convergence estimates; and some recursive methods can be tested against them.)
That is said, one may come up with a naive guess.
In the case of Magnus expansion, we see that $\dfrac1{4-(\Theta^{\Lie}_n)^{1/n}}$ is roughly linear in $n$,
which suggests that  $\mathrm C_\infty^{\Lie}$ may be around $3.85\pm0.15$, that is rather close to upper bound $4$, may be even equal to it.
This is, in some way, surprising because the improvement in the lower estimates earlier was relatively slow.
If $\mathrm C_\infty^{\Lie}<4$, then we can say that the Banach--Lie algebraic setting is slightly less natural than the Banach algebraic one.
On the other hand, if $\mathrm C_\infty^{\Lie}=4$, then we should conclude that the Banach--Lie algebraic setting is more natural and fundamental than
the Banach algebraic one.

Regarding the Baker--Campbell--Hausdorff expansion, the corresponding estimates are $2.93<\mathrm C^{\Lie}_2 \leq \mathrm D^{\Lie}_2
\leq 2\boldsymbol v_{\mathrm{Mi}}=5.4028570\ldots$.
These are about the convergence properties of the power series $\Gamma^{\Lie}(x,y)$.
Again, we know the actual coefficients   $\Gamma^{\Lie}_{n,m}$ for $n+m\leq16$, cf. Appendix \ref{sec:someGamma}.
Taking a look at them, firstly they suggest that $\mathrm C^{\Lie}_2 =\mathrm D^{\Lie}_2$.
Predicting the common value is certainly risky, but one would guess $5.15\pm0.15$, which is also close to the upper limit,
and, with a stretch, possibly equal to it.
\snewpage
\appendix
\section{Magnus commutators and $\Theta^{\Lie}_n$ of low degree}\plabel{sec:someTheta}

\begin{example}
\[\Theta_1^{\Lie}=1,\qquad \Theta_2^{\Lie}=\frac1{2!}\cdot\frac12,\qquad \Theta_3^{\Lie}=\frac1{3!}\cdot\frac13. \]
In fact,
\[\alpha_1=X^*_1,\qquad \alpha_2=X^*_{12},\qquad \alpha_3=X^*_{123} \]
attest the  minimality of the presentations
\[\mu_1^{\Lie}(X_1)=X_1,\quad \mu_2^{\Lie}(X_1,X_2)=\frac12X_{[1,2]},\quad \mu_3^{\Lie}(X_1,X_2,X_3)=\frac16X_{[[1,2],3]}+\frac16X_{[1,[2,3]]},\]
respectively. These are also unique.
\qedexer\end{example}
\begin{example}
\[\Theta_4^{\Lie}=\frac1{4!}\cdot\frac13. \]
The minimal presentations form a $\Delta_4=4$ dimensional polytope.
A nicest one is
\[\mu_4^{\Lie}(X_1,X_2,X_3,X_4)=\frac1{12}\Bigl(
X_{{[[1,3],[2,4]]}}+X_{{[[1,[2,3]],4]}}+X_{{[[1,2],[3,4]]}}+X_{{[1,[[2,3],4]]}}
\Bigr).\]

Indeed, minimality is attested by the linear functional
\[\alpha_4={X}^*_{1234}-{X}^*_{1342}-{X}^*_{1423}-{X}^*_{1432}.\]

The dihedrally invariant presentations form a $\Delta_4'=0$ dimensional polytope.
Then, necessarily,
\[\mu_4^{\Lie}(X_1,X_2,X_3,X_4)=\frac1{12}\Bigl(
X_{{[[1,3],[2,4]]}}^{\dyh5}
+3\cdot
X_{{[[1,2],[3,4]]}}^{\dyh5}
\Bigr).\]
Here we used the following convention: If the Lie monomial $M$ has orbit $\{M_1,\ldots,M_p\}$ under the dihedral action,
the $M^{\dyh p}$ denotes $\frac1p M_1+\ldots+\frac 1p M_p$.

[There are $10$ critical vertices spanning a $(4-1)!-1$ dimensional hyperface.
Thus, formally, the minimal presentations form a $10-(4-1)!=4$ dimensional abstract polytope.
The dihedrally invariant (dual) space is $2$ dimensional.
In particular, there is a minimal presentation with at most $2$ dihedrally symmetrized Lie monomials
(in fact, there is only one such).]
\qedexer
\end{example}
\begin{example}
\[\Theta_5^{\Lie}=\frac1{5!}\cdot\frac25. \]
The minimal presentations form a nontrivial $\Delta_5=32$ dimensional polytope.
A sort of most economical presentation is given by
\begin{align}
\mu_5^{\Lie}(X_1,&X_2,X_3,X_4,X_5)=\plabel{eq:penta}\\
=\frac1{120}\Bigl(&+4\,X_{{[[1,2],[3,[4,5]]]}}+4\,X_{{[[[1,2],[3,4]],5]}}+4\,X_{{[[[1,2],3],[4,5]]}}+4\,X_{{[1,[[2,3],[4,5]]]}}\notag\\
&+4\,X_{{[[1,[2,[3,5]]],4]}}-4\,X_{{[[[[1,3],4],5],2]}}+4\,X_{{[[1,[[2,3],4]],5]}}+4\,X_{{[1,[[2,[3,4]],5]]}}\notag\\
&+2\,X_{{[[1,3],[[2,4],5]]}}+2\,X_{{[[1,4],[[2,5],3]]}}-2\,X_{{[[[1,4],3],[2,5]]}}+2\,X_{{[[1,[2,4]],[3,5]]}}\notag\\
&+2\,X_{{[[[[1,5],4],3],2]}}-2\,X_{{[[[[1,5],2],3],4]}}+2\,X_{{[[[1,2],[3,5]],4]}}-2\,X_{{[[[1,3],[4,5]],2]}}\Bigl).\notag
\end{align}
(This is shorter than the formula given in Prato, Lamberti \cite{PL}.)

Minimality is attested by the linear functional
\begin{align}\alpha_5=
&+{X}^*_{12345}-{X}^*_{12453}-{X}^*_{12534}-{X}^*_{12543}-{X}^*_{13425}-{X}^*_{13452}\notag
\\&-{X}^*_{14235}-{X}^*_{14253}-{X}^*_{14325}-{X}^*_{14523}-{X}^*_{15234}+{X}^*_{15432}
\notag\end{align}
(taking values from $\{-1,0,1\}$ on Lie monomials).

The dihedrally invariant presentations form a $\Delta_5'=3$ dimensional polytope.
One such presentation, which uses all the vertices of the critical face, is
\begin{align}
\mu_5^{\Lie}(X_1,&X_2,X_3,X_4,X_5)=\plabel{eq:pentasymm}\\
=\frac1{50}\Bigl(&
+X_{{[1,[[2,[3,5]],4]]}}^{{\dyh12}}
-X_{{[1,[[[2,5],3],4]]}}^{{\dyh12}}
+X_{{[1,[2,[3,[4,5]]]]}}^{{\dyh6}}
+X_{{[1,[[2,4],[3,5]]]}}^{{\dyh6}}
\notag\\&
+X_{{[[1,4],[[2,5],3]]}}^{{\dyh6}}
+3\,X_{{[1,[[2,3],[4,5]]]}}^{{\dyh2}}
+6\,X_{{[1,[2,[[3,4],5]]]}}^{{\dyh6}}
+6\,X_{{[[1,3],[[2,4],5]]}}^{{\dyh6}}
\notag\Bigl).\end{align}

More generally, the  dihedrally invariant presentations are of form
\begin{align*}
\mu_5^{\Lie}(X_1,&X_2,X_3,X_4,X_5)=\notag\\
=\frac1{60}\Bigl(&
+6\lambda_0X_{{[1,[[2,[3,5]],4]]}}^{{\dyh12}}\\&
- \left( 3\lambda_1+6\lambda_4 \right) X_{{[1,[[[2,5],3],4]]}}^{{\dyh12}}
\notag\\&
+ \left( 3\lambda_3+6\lambda_4 \right) X_{{[1,[2,[3,[4,5]]]]}}^{{\dyh6}}\\&
+ \left( 3\lambda_1+2\lambda_2 \right) X_{{[1,[[2,4],[3,5]]]}}^{{\dyh6}}
\notag\\&
+ \left( 2\lambda_2+3\lambda_3 \right) X_{{[[1,4],[[2,5],3]]}}^{{\dyh6}}\\&
+ \left( 6\lambda_0+3\lambda_1+4\lambda_2+3\lambda_3 \right) X_{{[1,[[2,3],[4,5]]]}}^{{\dyh2}}
\notag\\&
+ \left( 6\lambda_0+9\lambda_1+8\lambda_2+6\lambda_3+6\lambda_4 \right) X_{{[1,[2,[[3,4],5]]]}}^{{\dyh6}}
\notag\\&
+ \left( 6\lambda_0+6\lambda_1+8\lambda_2+9\lambda_3+6\lambda_4 \right) X_{{[[1,3],[[2,4],5]]}}^{{\dyh6}}
\notag\Bigl),\end{align*}
where $\lambda_0+\lambda_1+\lambda_2+\lambda_3+\lambda_4=1$.
They form an irregular quadrangle based pyramid.
Here \eqref{eq:pentasymm} belongs to $(\lambda_0,\lambda_1,\lambda_2,\lambda_3,\lambda_4)=(\frac15,\frac25,0,\frac25,0),(\frac15,0,\frac35,0,\frac15) $, etc.;
i. e. it lies on the segment connecting the tip and the intersection of the diagonals of the base.
The most economical symmetric presentations belong to  $(\lambda_0,\lambda_1,\lambda_2,\lambda_3,\lambda_4)=(1,0,0,0,0)$ and $(0,0,0,0,1)$.
The dihedral symmetrization of \eqref{eq:penta} belongs to $(\lambda_0,\lambda_1,\lambda_2,\lambda_3,\lambda_4)=(0,0,1,0,0)$.

[There are $56$ critical vertices spanning a $(5-1)!-1$ dimensional hyperface.
Thus, formally, the minimal presentations form a $56-(5-1)!=32$ dimensional abstract polytope.
The dihedrally invariant (dual) space is $5$ dimensional.
In particular, there is a minimal presentation with at most $5$ dihedrally symmetrized Lie monomials
(but one can also do it with $4$, cf. the discussion above).]
\qedexer
\end{example}
\begin{example}
\[\Theta_6^{\Lie}=\frac1{6!}\cdot\frac{37}{60}. \]
The minimal presentations form a nontrivial $\Delta_5=370$ dimensional polytope.
A possible presentation is
\begin{align}
\mu_6^{\Lie}(X_1,&X_2,X_3,X_4,X_5,X_6)=\notag\\
=\frac1{240}\Bigl(&
         +4\,X_{{[[1,[3,5]],[[2,4],6]]}}+4\,X_{{[[1,[4,5]],[[2,3],6]]}}-4\,X_{{[[[1,4],5],[2,[3,6]]]}}+4\,X_{{[[1,[2,3]],[[4,5],6]]}}
\notag\\&+4\,X_{{[[[1,2],3],[4,[5,6]]]}}+4\,X_{{[[1,[2,4]],[[3,5],6]]}}+4\,X_{{[[1,[2,5]],[[3,4],6]]}}+4\,X_{{[[1,[3,4]],[[2,5],6]]}}
\notag\\&+4\,X_{{[1,[[2,[3,[4,5]]],6]]}}+4\,X_{{[1,[[[2,[3,4]],5],6]]}}-4\,X_{{[1,[[[[2,5],3],4],6]]}}+4\,X_{{[[1,3],[[2,[4,5]],6]]}}
\notag\\&+4\,X_{{[[1,[2,[[3,4],5]]],6]}}+4\,X_{{[[1,[[[2,3],4],5]],6]}}-4\,X_{{[[1,[[[2,5],4],3]],6]}}+4\,X_{{[[1,[[2,3],5]],[4,6]]}}
\notag\\&+4\,X_{{[[1,4],[[2,[3,5]],6]]}}+4\,X_{{[[1,5],[[2,[3,4]],6]]}}+4\,X_{{[[1,[[2,4],5]],[3,6]]}}+4\,X_{{[[1,[[3,4],5]],[2,6]]}}
\notag\\&-2\,X_{{[[[[1,3],[4,5]],6],2]}}-2\,X_{{[[[[[1,3],4],5],6],2]}}+2\,X_{{[[[[[1,3],6],5],4],2]}}+4\,X_{{[[[[[1,4],5],6],3],2]}}
\notag\\&+2\,X_{{[[[[[1,2],6],5],4],3]}}-2\,X_{{[[[[1,2],6],[4,5]],3]}}-2\,X_{{[[[[[1,2],4],5],6],3]}}-4\,X_{{[[[1,[4,[5,6]]],2],3]}}
\notag\\&+2\,X_{{[[1,[2,[3,[5,6]]]],4]}}-2\,X_{{[[[[1,[5,6]],2],3],4]}}-2\,X_{{[[[1,[5,6]],[2,3]],4]}}-4\,X_{{[[[[[1,2],3],6],5],4]}}
\notag\\&+2\,X_{{[[1,[2,[3,[4,6]]]],5]}}+2\,X_{{[[1,[[2,3],[4,6]]],5]}}-2\,X_{{[[[[1,[4,6]],2],3],5]}}+4\,X_{{[[[1,[2,[3,6]]],4],5]}}
\notag\\&+2\,X_{{[[1,2],[[3,[4,6]],5]]}}+2\,X_{{[[1,2],[3,[[4,5],6]]]}}+2\,X_{{[[1,2],[[3,4],[5,6]]]}}+4\,X_{{[[1,2],[[3,[4,5]],6]]}}
\notag\\&+2\,X_{{[[[1,2],[3,4]],[5,6]]}}+2\,X_{{[[[1,[2,3]],4],[5,6]]}}-2\,X_{{[[[[1,3],4],2],[5,6]]}}+4\,X_{{[[1,[[2,3],4]],[5,6]]}}\plabel{eq:formensol}
\\\notag&+2\,X_{{[[1,3],[[[2,4],5],6]]}}-2\,X_{{[[1,3],[[[2,6],4],5]]}}+2\,X_{{[[1,[2,[3,5]]],[4,6]]}}-2\,X_{{[[[[1,5],3],2],[4,6]]}}
\Bigl);
\end{align}
but many similarly looking ones exist. (For example,
\[+2\,X_{{[[1,[2,[3,4]]],[5,6]]}}+2\,X_{{[[1,[[2,3],4]],[5,6]]}}-2\,X_{{[[[[1,4],2],3],[5,6]]}}+4\,X_{{[[[1,[2,3]],4],[5,6]]}}\]
can replace line \eqref{eq:formensol}, but many other possibilities exist.)

Minimality is attested by
\begin{align}\alpha_6=
&+{X}^*_{123456}-{X}^*_{123564}-{X}^*_{123645}-{X}^*_{123654}-{X}^*_{124536}-{X}^*_{124563}\notag\\
&-{X}^*_{125346}-{X}^*_{125364}-{X}^*_{125436}-{X}^*_{125634}-{X}^*_{126345}+{X}^*_{126543}\notag\\
&-{X}^*_{134256}-{X}^*_{134526}-{X}^*_{134562}-{X}^*_{134625}+{X}^*_{135264}+{X}^*_{135642}\notag\\
&-{X}^*_{136245}-{X}^*_{136254}+{X}^*_{136452}+{X}^*_{136542}-{X}^*_{142356}-{X}^*_{142536}\notag\\
&+{X}^*_{142635}+{X}^*_{142653}-{X}^*_{143256}-{X}^*_{143625}-{X}^*_{145236}+{X}^*_{145362}\notag\\
&+{X}^*_{145632}+{X}^*_{146253}+{X}^*_{146352}+{X}^*_{146523}+{X}^*_{146532}-{X}^*_{152346}\notag\\
&+{X}^*_{152463}+{X}^*_{152643}+{X}^*_{153264}+{X}^*_{153462}+{X}^*_{153624}+{X}^*_{153642}\notag\\
&+{X}^*_{154263}+{X}^*_{154326}+{X}^*_{154362}+{X}^*_{154623}+{X}^*_{154632}+{X}^*_{156243}\notag\\
&+{X}^*_{156324}-{X}^*_{162345}+{X}^*_{162453}+{X}^*_{162534}+{X}^*_{162543}+{X}^*_{163425}\notag\\
&+{X}^*_{163452}+{X}^*_{163524}+{X}^*_{163542}+{X}^*_{164235}+{X}^*_{164253}+{X}^*_{164325}\notag\\
&+{X}^*_{164352}+{X}^*_{165234}+{X}^*_{165243}+{X}^*_{165324}-{X}^*_{165432}\notag
\end{align}
(taking values from $\{-1,0,1\}$ on Lie monomials).

The dihedrally invariant presentations form a $\Delta_6'=26$ dimensional polytope.
One such presentation is
\begin{align}
\mu_6^{\Lie}(X_1,X_2,&X_3,X_4,X_5,X_6)=\notag\\
=\frac1{15120}\Bigl(&
+28\,X_{{[1,[[2,4],[3,[5,6]]]]}}^{{\dyh14}}
+28\,X_{{[1,[[2,[3,[5,6]]],4]]}}^{{\dyh14}}
-28\,X_{{[1,[[2,[[3,6],5]],4]]}}^{{\dyh14}}
\notag\\&
-28\,X_{{[1,[[[2,[3,6]],5],4]]}}^{{\dyh14}}
-28\,X_{{[1,[[[2,[5,6]],3],4]]}}^{{\dyh7}}
+28\,X_{{[[1,4],[2,[[3,6],5]]]}}^{{\dyh14}}
\notag\\&
+28\,X_{{[[1,5],[[2,4],[3,6]]]}}^{{\dyh14}}
-28\,X_{{[[1,5],[[[2,4],6],3]]}}^{{\dyh14}}
+28\,X_{{[[[1,4],[2,6]],[3,5]]}}^{{\dyh14}}
\notag\\&
+54\,X_{{[1,[2,[3,[4,[5,6]]]]]}}^{{\dyh7}}
+54\,X_{{[1,[2,[[3,5],[4,6]]]]}}^{{\dyh14}}
-54\,X_{{[1,[2,[[[3,5],6],4]]]}}^{{\dyh14}}
\notag\\&
+54\,X_{{[1,[[2,5],[[3,4],6]]]}}^{{\dyh14}}
+54\,X_{{[1,[[2,5],[[3,6],4]]]}}^{{\dyh14}}
+54\,X_{{[1,[[2,[3,6]],[4,5]]]}}^{{\dyh14}}
\notag\\&
+54\,X_{{[1,[[2,[[3,4],6]],5]]}}^{{\dyh7}}
+54\,X_{{[[1,3],[[2,5],[4,6]]]}}^{{\dyh14}}
-54\,X_{{[[1,3],[[[2,5],6],4]]}}^{{\dyh14}}
\notag\\&
-54\,X_{{[[1,4],[[[2,6],3],5]]}}^{{\dyh7}}
-54\,X_{{[[1,5],[[[2,6],4],3]]}}^{{\dyh14}}
+54\,X_{{[[1,[3,5]],[2,[4,6]]]}}^{{\dyh14}}
\notag\\&
+54\,X_{{[[[1,3],[2,[4,6]]],5]}}^{{\dyh14}}
+54\,X_{{[[[1,5],3],[[2,6],4]]}}^{{\dyh7}}
+198\,X_{{[[1,4],[[2,[3,6]],5]]}}^{{\dyh14}}
\notag\\&
+227\,X_{{[1,[2,[3,[[4,5],6]]]]}}^{{\dyh14}}
+227\,X_{{[1,[2,[[3,[4,6]],5]]]}}^{{\dyh14}}
-227\,X_{{[1,[2,[[[3,6],4],5]]]}}^{{\dyh14}}
\notag\\&
+227\,X_{{[[1,3],[[2,[4,6]],5]]}}^{{\dyh14}}
-227\,X_{{[[1,3],[[[2,6],4],5]]}}^{{\dyh14}}
+227\,X_{{[[[1,3],5],[[2,6],4]]}}^{{\dyh14}}
\notag\\&
-251\,X_{{[[1,3],[[[2,6],5],4]]}}^{{\dyh14}}
-297\,X_{{[1,[2,[[[3,6],5],4]]]}}^{{\dyh14}}
+337\,X_{{[[1,3],[[[2,4],5],6]]}}^{{\dyh7}}
\notag\\&
+396\,X_{{[[1,[2,[4,6]]],[3,5]]}}^{{\dyh7}}
+833\,X_{{[1,[2,[[3,[4,5]],6]]]}}^{{\dyh7}}
+847\,X_{{[1,[[2,4],[[3,5],6]]]}}^{{\dyh14}}
\notag\\&
-1072\,X_{{[1,[[[[2,6],3],4],5]]}}^{{\dyh14}}
+1317\,X_{{[1,[[2,[3,[4,6]]],5]]}}^{{\dyh14}}
+1406\,X_{{[1,[2,[[3,4],[5,6]]]]}}^{{\dyh14}}
\Bigl);
\notag
\end{align}
where every vertex of the critical face appears (up to dihedral action).
As such, this is a longest symmetric presentation.
Other ones, like
\begin{align}
\mu_6^{\Lie}(X_1,&X_2,X_3,X_4,X_5,X_6)=\notag\\
=\frac1{60}\Bigl(&
-X_{{[1,[[[2,[5,6]],3],4]]}}^{{\dyh7}}
+X_{{[[1,4],[[2,[3,6]],5]]}}^{{\dyh14}}
+2\,X_{{[1,[[2,4],[[3,5],6]]]}}^{{\dyh14}}
\notag\\&
-2\,X_{{[[1,3],[[[2,6],5],4]]}}^{{\dyh14}}
+2\,X_{{[[1,[2,[4,6]]],[3,5]]}}^{{\dyh7}}
+3\,X_{{[[1,3],[[2,[4,6]],5]]}}^{{\dyh14}}
\notag\\&
+4\,X_{{[1,[2,[3,[4,[5,6]]]]]}}^{{\dyh7}}
+4\,X_{{[1,[[2,[3,[4,6]]],5]]}}^{{\dyh14}}
-5\,X_{{[1,[[[[2,6],3],4],5]]}}^{{\dyh14}}
\notag\\&
+6\,X_{{[1,[2,[[3,[4,5]],6]]]}}^{{\dyh7}}
-7\,X_{{[1,[2,[[[3,6],4],5]]]}}^{{\dyh14}}
\Bigl),
\notag
\end{align}
can be much simpler.

[There are $490$ critical vertices spanning a $(6-1)!-1$ dimensional hyperface.
Thus, formally, the minimal presentations form a $490-(6-1)!=370$ dimensional abstract polytope.
The dihedrally invariant (dual) space is $13$ dimensional.
In particular, there is a minimal presentation with at most $13$ dihedrally symmetrized Lie monomials
(cf. the formula above with $11$  dihedrally symmetrized Lie monomials).]
\qedexer
\end{example}
In the previous examples the critical faces of $\mathcal I_k^{\Lie}$, in the directions $\mu_k^{\Lie}$, were hyperfaces.
Thus the linear functionals $\alpha_k$ which defined them (took $1$ on the critical face),
were unique (restricted to $\mathfrak M_k^{\Lie}$).
In particular, as the dihedral action is a linear automorphism of   $\mathcal I_k^{\Lie}$,
the $\alpha_k$, restricted to $\mathfrak M_k^{\Lie}$, were invariant for the dihedral action.
The $\alpha_k$ also had the integrality property that they took only the values ${-1,0,1}$ on $\mathbf M^{\Lie}_k$.
The following example shows that this is more complicated in general.

\begin{example}
\[\Theta_7^{\Lie}=\frac1{7!}\cdot\frac{1621}{1540}.\]

Indeed, consider the linear function\\

$\alpha_7=\dfrac{1}{176}\Bigl(
115 {X}^*_{1234567}-61 {X}^*_{1234576}-61 {X}^*_{1234657}-176 {X}^*_{1234675}\\
-155 {X}^*_{1234756}-176 {X}^*_{1234765}-61 {X}^*_{1235467}-61 {X}^*_{1235476}-176 {X}^*_{1235647}-176 {X}^*_{1235674}\\
-96 {X}^*_{1235746}-35 {X}^*_{1235764}-176 {X}^*_{1236457}-115 {X}^*_{1236475}-176 {X}^*_{1236547}-14 {X}^*_{1236574}\\
-115 {X}^*_{1236745}+35 {X}^*_{1236754}-176 {X}^*_{1237456}-35 {X}^*_{1237465}-26 {X}^*_{1237546}+35 {X}^*_{1237564}\\
+26 {X}^*_{1237645}+176 {X}^*_{1237654}-61 {X}^*_{1243567}-61 {X}^*_{1243576}-61 {X}^*_{1243657}+21 {X}^*_{1243756}\\
-176 {X}^*_{1245367}-26 {X}^*_{1245376}-176 {X}^*_{1245637}-149 {X}^*_{1245673}-149 {X}^*_{1245736}+27{X}^*_{1245763}\\
-129 {X}^*_{1246357}+38 {X}^*_{1246375}-68 {X}^*_{1246537}+27 {X}^*_{1246573}-35 {X}^*_{1246735}+141 {X}^*_{1246753}\\
-129 {X}^*_{1247356}-47 {X}^*_{1247365}+21 {X}^*_{1247536}+162 {X}^*_{1247563}+26 {X}^*_{1247635}+176 {X}^*_{1247653}\\
-176 {X}^*_{1253467}-26 {X}^*_{1253476}-136 {X}^*_{1253647}+75 {X}^*_{1253746}+115 {X}^*_{1253764}-176 {X}^*_{1254367}\\
-26 {X}^*_{1254376}-26 {X}^*_{1254637}+27 {X}^*_{1254673}-97 {X}^*_{1254736}+53 {X}^*_{1254763}-94 {X}^*_{1256347}\\
+82 {X}^*_{1256374}+14 {X}^*_{1256437}+135 {X}^*_{1256473}+27 {X}^*_{1256734}+149 {X}^*_{1256743}-22 {X}^*_{1257346}\\
+112 {X}^*_{1257364}+18 {X}^*_{1257436}+147 {X}^*_{1257463}+142 {X}^*_{1257634}+114 {X}^*_{1257643}-176 {X}^*_{1263457}\\
-9 {X}^*_{1263475}-26 {X}^*_{1263547}+150 {X}^*_{1263574}+52 {X}^*_{1263745}+176 {X}^*_{1263754}-68 {X}^*_{1264357}\\
+99 {X}^*_{1264375}-28 {X}^*_{1264537}+114 {X}^*_{1264573}+129 {X}^*_{1264735}+147 {X}^*_{1264753}+14{X}^*_{1265347}\\
+176 {X}^*_{1265374}+122 {X}^*_{1265437}+149 {X}^*_{1265473}+121 {X}^*_{1265734}+135 {X}^*_{1265743}+124 {X}^*_{1267354}\\
+150 {X}^*_{1267435}+53 {X}^*_{1267453}+53 {X}^*_{1267534}+27 {X}^*_{1267543}-176 {X}^*_{1273456}+141 {X}^*_{1273564}\\
+141 {X}^*_{1273645}+176 {X}^*_{1273654}+68 {X}^*_{1274365}+150 {X}^*_{1274536}+176 {X}^*_{1274563}+141 {X}^*_{1274635}\\
+162 {X}^*_{1274653}+115 {X}^*_{1275346}+168 {X}^*_{1275364}+176 {X}^*_{1275436}+141 {X}^*_{1275463}+95 {X}^*_{1275634}\\
+27 {X}^*_{1275643}+176 {X}^*_{1276345}+150 {X}^*_{1276354}+176 {X}^*_{1276435}+27 {X}^*_{1276453}+27 {X}^*_{1276534}\\
-149 {X}^*_{1276543}-61 {X}^*_{1324567}-61 {X}^*_{1324576}-61 {X}^*_{1324657}+21 {X}^*_{1324756}-61 {X}^*_{1325467}\\
-61{X}^*_{1325476}-26 {X}^*_{1325647}-35 {X}^*_{1325746}-26 {X}^*_{1326457}+47 {X}^*_{1326475}-26 {X}^*_{1326547}\\
+21 {X}^*_{1326574}+47 {X}^*_{1326745}+21 {X}^*_{1326754}+35 {X}^*_{1327546}+26 {X}^*_{1327645}-155 {X}^*_{1342567}\\
+21 {X}^*_{1342576}+21 {X}^*_{1342657}+176 {X}^*_{1342675}+115 {X}^*_{1342756}+176 {X}^*_{1342765}-176 {X}^*_{1345267}\\
-176 {X}^*_{1345627}-164 {X}^*_{1345672}-164 {X}^*_{1345726}+12 {X}^*_{1345762}-129 {X}^*_{1346257}+47 {X}^*_{1346275}\\
+12 {X}^*_{1346572}-23 {X}^*_{1346725}+153 {X}^*_{1346752}-61 {X}^*_{1347256}+106 {X}^*_{1347265}+12 {X}^*_{1347526}\\
+165 {X}^*_{1347562}+118 {X}^*_{1347625}+165 {X}^*_{1347652}-96 {X}^*_{1352467}-35 {X}^*_{1352476}+75 {X}^*_{1352647}\\
+176 {X}^*_{1352674}+89 {X}^*_{1352746}+96 {X}^*_{1352764}-26 {X}^*_{1354267}+35 {X}^*_{1354276}+12 {X}^*_{1354672}\\
+7 {X}^*_{1354726}+33 {X}^*_{1354762}-22 {X}^*_{1356247}+154 {X}^*_{1356274}+115 {X}^*_{1356427}+165 {X}^*_{1356472}\\
+28 {X}^*_{1356724}+176 {X}^*_{1356742}-45 {X}^*_{1357246}+86 {X}^*_{1357264}+131 {X}^*_{1357426}+132 {X}^*_{1357462}\\
+42 {X}^*_{1357624}+109 {X}^*_{1357642}-149{X}^*_{1362457}-32 {X}^*_{1362475}-97 {X}^*_{1362547}+47 {X}^*_{1362574}\\
+144 {X}^*_{1362745}+47 {X}^*_{1362754}+21 {X}^*_{1364257}+56 {X}^*_{1364275}+150 {X}^*_{1364527}+153 {X}^*_{1364572}\\
+71 {X}^*_{1364725}+129{X}^*_{1364752}+18 {X}^*_{1365247}+107 {X}^*_{1365274}+176 {X}^*_{1365427}+165 {X}^*_{1365472}\\
+84 {X}^*_{1365724}+95 {X}^*_{1365742}+89 {X}^*_{1367245}+107 {X}^*_{1367254}+159 {X}^*_{1367425}+147 {X}^*_{1367452}\\
+84{X}^*_{1367524}-11 {X}^*_{1367542}-164 {X}^*_{1372456}-23 {X}^*_{1372465}+7 {X}^*_{1372546}+96 {X}^*_{1372564}\\
+101 {X}^*_{1372645}+176 {X}^*_{1372654}+12 {X}^*_{1374256}+65 {X}^*_{1374265}+98 {X}^*_{1374526}+176 {X}^*_{1374562}\\
+165 {X}^*_{1374625}+95 {X}^*_{1374652}+131 {X}^*_{1375246}+86 {X}^*_{1375264}+154 {X}^*_{1375426}+109 {X}^*_{1375462}\\
+42 {X}^*_{1375624}-67 {X}^*_{1375642}+176 {X}^*_{1376245}+154 {X}^*_{1376254}+165 {X}^*_{1376425}-11 {X}^*_{1376452}\\
+28 {X}^*_{1376524}-148 {X}^*_{1376542}-176 {X}^*_{1423567}+176 {X}^*_{1423675}+176 {X}^*_{1423756}+176 {X}^*_{1423765}\\
-115 {X}^*_{1425367}+47 {X}^*_{1425376}-9 {X}^*_{1425637}+27 {X}^*_{1425673}-32 {X}^*_{1425736}+74 {X}^*_{1425763}\\
+38 {X}^*_{1426357}+144 {X}^*_{1426375}+99 {X}^*_{1426537}+53 {X}^*_{1426573}+83 {X}^*_{1426735}+94 {X}^*_{1426753}\\
+47 {X}^*_{1427356}+56 {X}^*_{1427536}+115 {X}^*_{1427563}+3 {X}^*_{1427635}-176 {X}^*_{1432567}+176 {X}^*_{1432675}\\
+176 {X}^*_{1432756}+176 {X}^*_{1432765}-35 {X}^*_{1435267}+12 {X}^*_{1435672}-23 {X}^*_{1435726}+12 {X}^*_{1435762}\\
-47 {X}^*_{1436257}+68 {X}^*_{1436527}+12 {X}^*_{1436572}+3 {X}^*_{1436725}+3 {X}^*_{1436752}+106 {X}^*_{1437256}\\
+144 {X}^*_{1437265}+65 {X}^*_{1437526}+50 {X}^*_{1437562}+83 {X}^*_{1437625}-11{X}^*_{1437652}-115 {X}^*_{1452367}\\
+47 {X}^*_{1452376}+52 {X}^*_{1452637}+176 {X}^*_{1452673}+144 {X}^*_{1452736}+162 {X}^*_{1452763}+26 {X}^*_{1453267}\\
+26 {X}^*_{1453276}+141 {X}^*_{1453627}+176 {X}^*_{1453672}+101{X}^*_{1453726}+40 {X}^*_{1453762}+176 {X}^*_{1456273}\\
+176 {X}^*_{1456327}+176 {X}^*_{1456372}+27 {X}^*_{1456723}+176 {X}^*_{1456732}+89 {X}^*_{1457236}+135 {X}^*_{1457263}\\
+176 {X}^*_{1457326}+156 {X}^*_{1457362}+53{X}^*_{1457623}+26 {X}^*_{1457632}-35 {X}^*_{1462357}+83 {X}^*_{1462375}\\
+129 {X}^*_{1462537}+135 {X}^*_{1462573}+100 {X}^*_{1462735}+88 {X}^*_{1462753}+26 {X}^*_{1463257}+3 {X}^*_{1463275}\\
+141 {X}^*_{1463527}+82 {X}^*_{1463572}+49 {X}^*_{1463725}-12 {X}^*_{1463752}+150 {X}^*_{1465237}+162 {X}^*_{1465273}\\
+176 {X}^*_{1465327}+70 {X}^*_{1465372}+74 {X}^*_{1465723}+47 {X}^*_{1465732}+118 {X}^*_{1467235}+176 {X}^*_{1467253}\\
+132 {X}^*_{1467325}+71 {X}^*_{1467352}+115 {X}^*_{1467523}-47 {X}^*_{1467532}-23 {X}^*_{1472356}+3 {X}^*_{1472365}\\
+71 {X}^*_{1472536}+176 {X}^*_{1472563}+49{X}^*_{1472635}+61 {X}^*_{1472653}+118 {X}^*_{1473256}+83 {X}^*_{1473265}\\
+165 {X}^*_{1473526}+176 {X}^*_{1473562}+100 {X}^*_{1473625}+12 {X}^*_{1473652}+159 {X}^*_{1475236}+88 {X}^*_{1475263}\\
+165 {X}^*_{1475326}+70 {X}^*_{1475362}+94 {X}^*_{1475623}-47 {X}^*_{1475632}+132 {X}^*_{1476235}+61{X}^*_{1476253}\\
+118 {X}^*_{1476325}-58 {X}^*_{1476352}-176 {X}^*_{1476532}-176 {X}^*_{1523467}+176 {X}^*_{1523674}+176 {X}^*_{1523746}\\
+176 {X}^*_{1523764}-14 {X}^*_{1524367}+21 {X}^*_{1524376}+150 {X}^*_{1524637}+176 {X}^*_{1524673}+47 {X}^*_{1524736}\\
+35 {X}^*_{1524763}+82 {X}^*_{1526347}+176 {X}^*_{1526374}+176 {X}^*_{1526437}+154 {X}^*_{1526473}+132 {X}^*_{1526734}\\
+105 {X}^*_{1526743}+154 {X}^*_{1527346}+132 {X}^*_{1527364}+107 {X}^*_{1527436}+101 {X}^*_{1527463}+80 {X}^*_{1527634}\\
-62 {X}^*_{1527643}-35 {X}^*_{1532467}+115 {X}^*_{1532647}+176 {X}^*_{1532674}+96 {X}^*_{1532746}+70 {X}^*_{1532764}\\
+35 {X}^*_{1534267}+141 {X}^*_{1534627}+176 {X}^*_{1534672}+96 {X}^*_{1534726}+26 {X}^*_{1534762}+112 {X}^*_{1536247}\\
+132 {X}^*_{1536274}+168 {X}^*_{1536427}+153 {X}^*_{1536472}+106 {X}^*_{1536724}+78 {X}^*_{1536742}+86 {X}^*_{1537246}\\
+85 {X}^*_{1537264}+86 {X}^*_{1537426}+{X}^*_{1537462}-3 {X}^*_{1537624}-45 {X}^*_{1537642}+35 {X}^*_{1542367}\\
+21 {X}^*_{1542376}+176 {X}^*_{1542637}+176 {X}^*_{1542673}+47 {X}^*_{1542736}+12 {X}^*_{1542763}+176 {X}^*_{1543267}\\
+176 {X}^*_{1543627}+176 {X}^*_{1543672}+176 {X}^*_{1543726}+124 {X}^*_{1546237}+124 {X}^*_{1546273}+150 {X}^*_{1546327}\\
+35 {X}^*_{1546372}+36 {X}^*_{1546723}+9 {X}^*_{1546732}+107 {X}^*_{1547236}+6 {X}^*_{1547263}+154 {X}^*_{1547326}\\
+22{X}^*_{1547362}-46 {X}^*_{1547623}-141 {X}^*_{1547632}+27 {X}^*_{1562347}+132 {X}^*_{1562374}+121 {X}^*_{1562437}\\
+143 {X}^*_{1562473}+176 {X}^*_{1562734}+94 {X}^*_{1562743}+142 {X}^*_{1563247}+80 {X}^*_{1563274}+95 {X}^*_{1563427}\\
+141 {X}^*_{1563472}+132 {X}^*_{1563724}-22 {X}^*_{1563742}+53 {X}^*_{1564237}+27 {X}^*_{1564273}+27 {X}^*_{1564327}\\
-9 {X}^*_{1564372}+115 {X}^*_{1564723}-35 {X}^*_{1564732}+176 {X}^*_{1567234}+176 {X}^*_{1567243}+176 {X}^*_{1567324}\\
-176 {X}^*_{1567432}+28 {X}^*_{1572346}+106{X}^*_{1572364}+84 {X}^*_{1572436}+95 {X}^*_{1572463}+132 {X}^*_{1572634}\\
+20 {X}^*_{1572643}+42 {X}^*_{1573246}-3 {X}^*_{1573264}+42 {X}^*_{1573426}+45 {X}^*_{1573462}+85 {X}^*_{1573624}\\
-{X}^*_{1573642}+84 {X}^*_{1574236}-11 {X}^*_{1574263}+28 {X}^*_{1574326}-78 {X}^*_{1574362}+15 {X}^*_{1574623}\\
-153 {X}^*_{1574632}+176 {X}^*_{1576234}+61 {X}^*_{1576243}+70 {X}^*_{1576324}-26 {X}^*_{1576342}-35 {X}^*_{1576423}\\
-176 {X}^*_{1576432}-149 {X}^*_{1623457}+27 {X}^*_{1623475}+27 {X}^*_{1623547}+176 {X}^*_{1623574}+176 {X}^*_{1623745}\\
+176 {X}^*_{1623754}+27 {X}^*_{1624357}+53 {X}^*_{1624375}+114 {X}^*_{1624537}+176 {X}^*_{1624573}+135 {X}^*_{1624735}\\
+46{X}^*_{1624753}+135 {X}^*_{1625347}+154 {X}^*_{1625374}+149 {X}^*_{1625437}+44 {X}^*_{1625473}+143 {X}^*_{1625734}\\
-19 {X}^*_{1625743}+176 {X}^*_{1627345}+124 {X}^*_{1627354}+162 {X}^*_{1627435}+18 {X}^*_{1627453}+27 {X}^*_{1627534}\\
-149 {X}^*_{1627543}+27 {X}^*_{1632457}+74 {X}^*_{1632475}+53 {X}^*_{1632547}+35 {X}^*_{1632574}+162 {X}^*_{1632745}\\
+12 {X}^*_{1632754}+162 {X}^*_{1634257}+115 {X}^*_{1634275}+176 {X}^*_{1634527}+176 {X}^*_{1634572}+176 {X}^*_{1634725}\\
+58 {X}^*_{1634752}+147 {X}^*_{1635247}+101{X}^*_{1635274}+141 {X}^*_{1635427}+47 {X}^*_{1635472}+95 {X}^*_{1635724}\\
-70 {X}^*_{1635742}+135 {X}^*_{1637245}+6 {X}^*_{1637254}+88 {X}^*_{1637425}-12 {X}^*_{1637452}-11 {X}^*_{1637524}\\
-176 {X}^*_{1637542}+141 {X}^*_{1642357}+94 {X}^*_{1642375}+147 {X}^*_{1642537}+46 {X}^*_{1642573}+88 {X}^*_{1642735}\\
-71 {X}^*_{1642753}+176 {X}^*_{1643257}+162 {X}^*_{1643527}+47 {X}^*_{1643572}+61 {X}^*_{1643725}-71 {X}^*_{1643752}\\
+53 {X}^*_{1645237}+18 {X}^*_{1645273}+27 {X}^*_{1645327}-47 {X}^*_{1645372}+106 {X}^*_{1645723}-70 {X}^*_{1645732}\\
+176 {X}^*_{1647235}+105 {X}^*_{1647253}+61 {X}^*_{1647325}+12 {X}^*_{1647352}+59 {X}^*_{1647523}-82 {X}^*_{1647532}\\
+149 {X}^*_{1652347}+105 {X}^*_{1652374}+135 {X}^*_{1652437}-19 {X}^*_{1652473}+94 {X}^*_{1652734}-82 {X}^*_{1652743}\\
+114 {X}^*_{1653247}-62 {X}^*_{1653274}+27 {X}^*_{1653427}-26 {X}^*_{1653472}+20 {X}^*_{1653724}-156 {X}^*_{1653742}\\
+27 {X}^*_{1654237}-149 {X}^*_{1654273}-149 {X}^*_{1654327}-176 {X}^*_{1654372}-176 {X}^*_{1654732}+176 {X}^*_{1657234}\\
+40 {X}^*_{1657243}+61 {X}^*_{1657324}-40 {X}^*_{1657342}-35 {X}^*_{1657423}-176 {X}^*_{1657432}+27 {X}^*_{1672345}\\
+36 {X}^*_{1672354}+74 {X}^*_{1672435}+106 {X}^*_{1672453}+115 {X}^*_{1672534}+53 {X}^*_{1673245}-46 {X}^*_{1673254}\\
+94 {X}^*_{1673425}+11 {X}^*_{1673452}+15 {X}^*_{1673524}-50 {X}^*_{1673542}+115 {X}^*_{1674235}+59 {X}^*_{1674253}\\
-3 {X}^*_{1674352}+56 {X}^*_{1674523}-12 {X}^*_{1674532}-35 {X}^*_{1675243}-35 {X}^*_{1675324}-12 {X}^*_{1675342}\\
-12 {X}^*_{1675423}-12 {X}^*_{1675432}-164 {X}^*_{1723456}+12 {X}^*_{1723465}+12 {X}^*_{1723546}+176 {X}^*_{1723564}\\
+176 {X}^*_{1723645}+176 {X}^*_{1723654}+12 {X}^*_{1724356}+12 {X}^*_{1724365}+153 {X}^*_{1724536}+176 {X}^*_{1724563}\\
+82 {X}^*_{1724635}+47 {X}^*_{1724653}+165 {X}^*_{1725346}+153 {X}^*_{1725364}+165 {X}^*_{1725436}+47 {X}^*_{1725463}\\
+141 {X}^*_{1725634}-26 {X}^*_{1725643}+176 {X}^*_{1726345}+35 {X}^*_{1726354}+70 {X}^*_{1726435}-47 {X}^*_{1726453}\\
-9 {X}^*_{1726534}-176 {X}^*_{1726543}+12 {X}^*_{1732456}+12 {X}^*_{1732465}+33 {X}^*_{1732546}+26 {X}^*_{1732564}\\
+40 {X}^*_{1732645}+165 {X}^*_{1734256}+50 {X}^*_{1734265}+176 {X}^*_{1734526}+148 {X}^*_{1734562}+176 {X}^*_{1734625}\\
+11 {X}^*_{1734652}+132 {X}^*_{1735246}+{X}^*_{1735264}+109 {X}^*_{1735426}+67 {X}^*_{1735462}+45 {X}^*_{1735624}\\
-109 {X}^*_{1735642}+156 {X}^*_{1736245}+22 {X}^*_{1736254}+70 {X}^*_{1736425}-95 {X}^*_{1736452}-78 {X}^*_{1736524}\\
-176 {X}^*_{1736542}+153 {X}^*_{1742356}+3 {X}^*_{1742365}+129 {X}^*_{1742536}+58 {X}^*_{1742563}-12 {X}^*_{1742635}\\
-71 {X}^*_{1742653}+165 {X}^*_{1743256}-11 {X}^*_{1743265}+95 {X}^*_{1743526}+11 {X}^*_{1743562}+12 {X}^*_{1743625}\\
-147 {X}^*_{1743652}+147 {X}^*_{1745236}-12 {X}^*_{1745263}-11 {X}^*_{1745326}-95 {X}^*_{1745362}+11 {X}^*_{1745623}\\
-165 {X}^*_{1745632}+71 {X}^*_{1746235}+12 {X}^*_{1746253}-58 {X}^*_{1746325}-129 {X}^*_{1746352}-3 {X}^*_{1746523}\\
-153 {X}^*_{1746532}+176 {X}^*_{1752346}+78 {X}^*_{1752364}+95 {X}^*_{1752436}-70 {X}^*_{1752463}-22 {X}^*_{1752634}\\
-156 {X}^*_{1752643}+109 {X}^*_{1753246}-45 {X}^*_{1753264}-67 {X}^*_{1753426}-109 {X}^*_{1753462}-{X}^*_{1753624}\\
-132 {X}^*_{1753642}-11 {X}^*_{1754236}-176 {X}^*_{1754263}-148 {X}^*_{1754326}-176 {X}^*_{1754362}-50 {X}^*_{1754623}\\
-165 {X}^*_{1754632}-40 {X}^*_{1756243}-26 {X}^*_{1756324}-33 {X}^*_{1756342}-12 {X}^*_{1756423}-12 {X}^*_{1756432}\\
+176 {X}^*_{1762345}+9{X}^*_{1762354}+47 {X}^*_{1762435}-70 {X}^*_{1762453}-35 {X}^*_{1762534}-176 {X}^*_{1762543}\\
+26 {X}^*_{1763245}-141 {X}^*_{1763254}-47 {X}^*_{1763425}-165 {X}^*_{1763452}-153 {X}^*_{1763524}-165 {X}^*_{1763542}\\
-47{X}^*_{1764235}-82 {X}^*_{1764253}-176 {X}^*_{1764325}-153 {X}^*_{1764352}-12 {X}^*_{1764523}-12 {X}^*_{1764532}\\
-176 {X}^*_{1765234}-176 {X}^*_{1765243}-176 {X}^*_{1765324}-12 {X}^*_{1765342}-12 {X}^*_{1765423}+164 {X}^*_{1765432}\Bigr)
$.\\

Up to $\pm$ sign, this takes values on $\mathbf M_7^{\Lie}$ from the set
\[\Bigl\{0,{\frac {1}{176}},{\frac {3}{176}},{\frac {3}{88}},{\frac {7}{176}},
{\frac {9}{176}},{\frac {5}{88}},\frac{1}{16},{\frac {3}{44}},{\frac {13}{176}
},{\frac {7}{88}},{\frac {15}{176}},{\frac {17}{176}},{\frac {9}{88}},
{\frac {19}{176}},{\frac {5}{44}},\]\[{\frac {21}{176}},\frac18,{\frac {23}{
176}},{\frac {3}{22}},{\frac {25}{176}},{\frac {13}{88}},{\frac {27}{
176}},{\frac {7}{44}},{\frac {15}{88}},\frac2{11},\frac3{16},{\frac {35}{176}},{
\frac {9}{44}},{\frac {19}{88}},{\frac {5}{22}},{\frac {41}{176}},\]\[{
\frac {21}{88}},{\frac {43}{176}},\frac14,{\frac {45}{176}},{\frac {23}{88
}},{\frac {47}{176}},\frac3{11},{\frac {49}{176}},{\frac {25}{88}},{\frac {
13}{44}},{\frac {53}{176}},{\frac {5}{16}},{\frac {7}{22}},{\frac {29}
{88}},{\frac {59}{176}},{\frac {15}{44}},\]\[{\frac {61}{176}},{\frac {31}
{88}},{\frac {63}{176}},{\frac {65}{176}},\frac38,{\frac {67}{176}},{
\frac {17}{44}},{\frac {35}{88}},{\frac {71}{176}},{\frac {73}{176}},{
\frac {37}{88}},{\frac {75}{176}},{\frac {39}{88}},{\frac {79}{176}},{
\frac {5}{11}},{\frac {81}{176}},\]\[{\frac {41}{88}},{\frac {83}{176}},{
\frac {21}{44}},{\frac {85}{176}},{\frac {43}{88}},{\frac {87}{176}},\frac12,{\frac {89}{176}},{\frac {47}{88}},{\frac {95}{176}},{\frac {6}{11}
},{\frac {97}{176}},{\frac {49}{88}},{\frac {9}{16}},{\frac {25}{44}},
{\frac {101}{176}},\]\[{\frac {103}{176}},{\frac {105}{176}},{\frac {53}{
88}},{\frac {107}{176}},{\frac {27}{44}},{\frac {109}{176}},{\frac {7}
{11}},{\frac {57}{88}},{\frac {115}{176}},{\frac {29}{44}},{\frac {117
}{176}},{\frac {59}{88}},{\frac {11}{16}},{\frac {61}{88}},{\frac {123
}{176}},{\frac {31}{44}},\]\[{\frac {63}{88}},{\frac {129}{176}},{\frac {
65}{88}},{\frac {131}{176}},\frac34,{\frac {67}{88}},{\frac {135}{176}},{
\frac {17}{22}},{\frac {137}{176}},{\frac {141}{176}},{\frac {71}{88}}
,{\frac {13}{16}},{\frac {9}{11}},{\frac {147}{176}},{\frac {37}{44}},
{\frac {149}{176}},\]\[{\frac {75}{88}},{\frac {153}{176}},{\frac {7}{8}},
{\frac {155}{176}},{\frac {39}{44}},{\frac {79}{88}},{\frac {159}{176}
},{\frac {81}{88}},{\frac {163}{176}},{\frac {41}{44}},{\frac {15}{16}
},{\frac {167}{176}},{\frac {21}{22}},{\frac {85}{88}},{\frac {171}{176}},1\Bigr\}.
\]
\\
Yet, $\alpha_7$ takes the value $\frac{1621}{1540}$ on $\mu^{\Lie}_7(X_1\ldots,X_7)$.
This proves $\Theta_7^{\Lie}\geq\frac1{7!}\cdot\frac{1621}{1540}$.
On the other hand, one can prove that $\alpha_7=1$ defines the critical face
of $\mathcal I_7^{\Lie}$, such that $\frac{\mu^{\Lie}_7(X_1\ldots,X_7)}{\alpha_7(\mu^{\Lie}_7(X_1\ldots,X_7))}$ is
an interior point of it.
In fact, a minimal presentation, which uses all vertices of the critical face, is given by\\
\begin{align*}
\mu_7^{\Lie}&(X_1,X_2,X_3,X_4,X_5,X_6,X_7)=\dfrac1{207900}\Bigl(
-60\,X_{{[[[1,7],[3,5]],[2,[4,6]]]}}^{{\dyh16}}\\&
+80\,X_{{[1,[2,[[3,4],[5,[6,7]]]]]}}^{{\dyh8}}
-80\,X_{{[1,[2,[[[3,[6,7]],4],5]]]}}^{{\dyh16}}
-80\,X_{{[1,[2,[[[3,[6,7]],5],4]]]}}^{{\dyh16}}\\&
+80\,X_{{[1,[[[2,[3,[4,7]]],5],6]]}}^{{\dyh16}}
-80\,X_{{[1,[[[[2,[4,7]],3],5],6]]}}^{{\dyh16}}
-80\,X_{{[[1,4],[[[[2,5],6],3],7]]}}^{{\dyh16}}\\&
-210\,X_{{[[[1,4],[[2,[5,7]],6]],3]}}^{{\dyh8}}
+240\,X_{{[1,[2,[[3,5],[[4,6],7]]]]}}^{{\dyh16}}
-240\,X_{{[[1,3],[2,[[[4,7],5],6]]]}}^{{\dyh16}}\\&
-240\,X_{{[[1,3],[2,[[[4,7],6],5]]]}}^{{\dyh16}}
-240\,X_{{[[1,6],[[[[2,5],3],4],7]]}}^{{\dyh8}}
-240\,X_{{[[1,6],[[[[2,5],4],3],7]]}}^{{\dyh8}}\\&
+240\,X_{{[[1,6],[[[2,4],[3,5]],7]]}}^{{\dyh16}}
-240\,X_{{[[[1,[4,7]],[3,6]],[2,5]]}}^{{\dyh8}}
-330\,X_{{[[[1,4],[[2,6],[5,7]]],3]}}^{{\dyh8}}\\&
-336\,X_{{[[1,4],[2,[[[3,7],6],5]]]}}^{{\dyh16}}
+336\,X_{{[[1,5],[[[2,4],[3,6]],7]]}}^{{\dyh16}}
-336\,X_{{[[1,5],[[[[2,4],6],3],7]]}}^{{\dyh16}}\\&
-336\,X_{{[[1,5],[[[[2,6],3],4],7]]}}^{{\dyh8}}
-336\,X_{{[[1,5],[[[[2,6],4],3],7]]}}^{{\dyh16}}
+555\,X_{{[[[[1,5],[3,7]],[2,6]],4]}}^{{\dyh1}}\\&
+580\,X_{{[1,[2,[[3,[4,5]],[6,7]]]]}}^{{\dyh8}}
+616\,X_{{[1,[[2,[3,[4,[5,7]]]],6]]}}^{{\dyh16}}
-616\,X_{{[[1,3],[[[[2,5],6],4],7]]}}^{{\dyh16}}\\&
-616\,X_{{[[1,3],[[[[2,6],5],4],7]]}}^{{\dyh16}}
+770\,X_{{[1,[2,[3,[[4,5],[6,7]]]]]}}^{{\dyh16}}
+770\,X_{{[1,[[2,3],[[4,6],[5,7]]]]}}^{{\dyh8}}\\&
+770\,X_{{[1,[[2,3],[[4,[5,7]],6]]]}}^{{\dyh16}}
-780\,X_{{[[1,5],[[[[2,6],3],7],4]]}}^{{\dyh8}}
-890\,X_{{[[1,4],[[[[2,6],3],5],7]]}}^{{\dyh16}}\\&
-890\,X_{{[[1,4],[[[[2,6],5],3],7]]}}^{{\dyh16}}
-1140\,X_{{[1,[[[2,[[3,7],6]],4],5]]}}^{{\dyh16}}
+1464\,X_{{[1,[2,[[[3,5],[4,6]],7]]]}}^{{\dyh16}}\\&
+1464\,X_{{[1,[2,[[[3,[4,6]],5],7]]]}}^{{\dyh16}}
-1464\,X_{{[1,[2,[[[[3,5],6],4],7]]]}}^{{\dyh16}}
-1464\,X_{{[1,[2,[[[[3,6],4],5],7]]]}}^{{\dyh16}}\\&
-1464\,X_{{[1,[2,[[[[3,6],5],4],7]]]}}^{{\dyh16}}
-1500\,X_{{[1,[[[2,[[3,6],7]],5],4]]}}^{{\dyh16}}
-1576\,X_{{[[1,3],[[[[2,6],4],5],7]]}}^{{\dyh16}}\\&
-1580\,X_{{[1,[[2,3],[[[4,7],5],6]]]}}^{{\dyh16}}
-1604\,X_{{[1,[2,[[[[3,5],6],7],4]]]}}^{{\dyh16}}
-1740\,X_{{[1,[[2,[[[3,6],7],4]],5]]}}^{{\dyh16}}\\&
+1740\,X_{{[[1,[3,6]],[[2,[4,7]],5]]}}^{{\dyh16}}
+1800\,X_{{[1,[[2,[3,5]],[4,[6,7]]]]}}^{{\dyh16}}
+1800\,X_{{[1,[[2,4],[[3,[5,6]],7]]]}}^{{\dyh16}}\\&
+1900\,X_{{[1,[[2,3],[[4,5],[6,7]]]]}}^{{\dyh4}}
+1920\,X_{{[[1,3],[[2,5],[[4,6],7]]]}}^{{\dyh8}}
+1980\,X_{{[[1,3],[[2,[4,6]],[5,7]]]}}^{{\dyh16}}\\&
+2100\,X_{{[1,[[2,[3,[[4,5],7]]],6]]}}^{{\dyh8}}
+2166\,X_{{[1,[2,[3,[[4,[5,6]],7]]]]}}^{{\dyh16}}
-2166\,X_{{[1,[2,[3,[[[4,7],6],5]]]]}}^{{\dyh16}}\\&
-2754\,X_{{[1,[2,[3,[[[4,7],5],6]]]]}}^{{\dyh16}}
-2760\,X_{{[[1,6],[[[2,5],4],[3,7]]]}}^{{\dyh8}}
+2832\,X_{{[1,[2,[[3,[[4,5],6]],7]]]}}^{{\dyh8}}\\&
+3480\,X_{{[1,[2,[[3,[4,[6,7]]],5]]]}}^{{\dyh16}}
+4430\,X_{{[[1,3],[[[[2,7],6],5],4]]}}^{{\dyh16}}
+4520\,X_{{[1,[[2,[3,[4,[6,7]]]],5]]}}^{{\dyh16}}\\&
+4632\,X_{{[1,[2,[[[3,[4,5]],6],7]]]}}^{{\dyh8}}
+4680\,X_{{[1,[[2,[3,6]],[[4,5],7]]]}}^{{\dyh16}}
+4740\,X_{{[1,[[2,5],[[[3,4],6],7]]]}}^{{\dyh16}}\\&
+4800\,X_{{[[1,3],[[2,[4,[5,7]]],6]]}}^{{\dyh8}}
-4848\,X_{{[1,[[2,4],[[[3,7],6],5]]]}}^{{\dyh16}}
-4860\,X_{{[1,[[2,[[[3,7],4],6]],5]]}}^{{\dyh16}}\\&
+4862\,X_{{[1,[2,[[[3,4],[5,6]],7]]]}}^{{\dyh16}}
-4980\,X_{{[[1,4],[[[[2,7],3],5],6]]}}^{{\dyh8}}
+5058\,X_{{[1,[2,[3,[4,[[5,6],7]]]]]}}^{{\dyh16}}\\&
-5148\,X_{{[1,[[2,[[[3,7],4],5]],6]]}}^{{\dyh16}}
+5220\,X_{{[1,[[2,[3,5]],[[4,6],7]]]}}^{{\dyh8}}
-5280\,X_{{[1,[[[[2,[5,7]],3],4],6]]}}^{{\dyh16}}\\&
-5760\,X_{{[[1,3],[[[[2,7],4],5],6]]}}^{{\dyh16}}
+5880\,X_{{[[1,[3,5]],[[2,[4,7]],6]]}}^{{\dyh16}}
+6540\,X_{{[[1,4],[[2,[3,[5,7]]],6]]}}^{{\dyh16}}\\&
+6820\,X_{{[1,[2,[[3,4],[[5,6],7]]]]}}^{{\dyh16}}
+6900\,X_{{[1,[[[[[2,7],4],5],6],3]]}}^{{\dyh16}}
+7280\,X_{{[1,[[2,[3,[4,7]]],[5,6]]]}}^{{\dyh16}}\\&
-8268\,X_{{[1,[2,[[[[3,7],4],5],6]]]}}^{{\dyh16}}
+8384\,X_{{[1,[2,[[3,[4,[5,7]]],6]]]}}^{{\dyh16}}
-8950\,X_{{[[1,3],[[[[2,5],6],7],4]]}}^{{\dyh16}}\\&
+9016\,X_{{[1,[[2,4],[[[3,5],6],7]]]}}^{{\dyh16}}
+9048\,X_{{[1,[2,[[[[3,7],6],5],4]]]}}^{{\dyh16}}
+9780\,X_{{[[1,[2,[4,6]]],[3,[5,7]]]}}^{{\dyh16}}
\Bigr).
\end{align*}
Here the critical face is not a $(7-1)!-1$ dimensional hyperface but a $(7-1)!-1-3$ dimensional one.
Thus, restricted to $\mathfrak M_7^{\Lie}$, there are several linear functionals which define the critical face
(forming the interior of a $3$ dimensional polytope with $20$ vertices).
However, it is true that $\alpha_7$ is the only one which is symmetrical to the dihedral action.
(The corresponding $20$ hyperface-defining linear functionals are also not particulary nice either.
This is expected, as their dihedral symmetrization is the $\alpha_7$ above.
There are no good integral face-defining integral linear functionals because those would lead
to the symmetrized $\alpha_7$ with coefficients with multiples of $1/16$.)

[There are $1141$ critical vertices spanning a $(7-1)!-3-1$ dimensional face.
Thus, formally, the minimal presentations form a $1141-(7-1)!+3=424$ dimensional abstract polytope.
The dihedrally invariant (dual) space is $60$ dimensional.
In particular, there is a minimal presentation with at most $60$ dihedrally symmetrized Lie monomials.]
\qedexer
\end{example}
\begin{example}
\[\Theta_8^{\Lie}=\frac1{8!}\cdot
\frac{5242130799984621832318}{2419342933460499216625}
=\frac1{8!}\cdot2.166758060\ldots\quad.\]

This can be done as in the previous cases.
There is little point in printing out the actual formulas, which are rather long;
although the minimal presentations can be still fairly short.

There are $5934$ critical vertices spanning a $(8-1)!-1$ dimensional hyperface.
Thus, formally, the minimal presentations form a $5934-(8-1)!=894$ dimensional abstract polytope.
The dihedrally invariant (dual) space is $306$ dimensional.
In particular, there is a minimal presentation with at most $306$ dihedrally symmetrized Lie monomials.
\qedexer
\end{example}

\begin{example}
\begin{align*}
\Theta_9^{\Lie}=&\frac1{9!}\cdot
\frac{
\left(\begin{array}{l}
\phantom{+\,00000000000000000000000000000000000000}\!2588638005\cdot10^{48\cdot5}\\
+\,137011373976287809546942042433176351096359261476\cdot10^{48\cdot4}\\
+\,123833498094907230678264389179973397107285604160\cdot10^{48\cdot3}\\
+\,536705848906821186613961139113214863359234861372\cdot10^{48\cdot2}\\
+\,558922154447684390889961674834042514266130495835\cdot10^{48}\\
+\,178471091815540406692440312028332551019654336199
\end{array}\right)
}{
\left(\begin{array}{l}
\phantom{+\,000000000000000000000000000000000000000}512183840\cdot10^{48\cdot5}\\
+\,494197597994580290019457156210276499193075866836\cdot10^{48\cdot4}\\
+\,789345554793363030217885865082978090374815737781\cdot10^{48\cdot3}\\
+\,676980678683847851476408188922894174486569439031\cdot10^{48\cdot2}\\
+\,333056158992821946623532028245051548010609826960\cdot10^{48}\\
+\,543053601386302542647543865480289939077893915552
\end{array}\right)
}
\\
=&\frac1{9!}\cdot5.054118854\ldots\qquad.
\end{align*}

There are $47715$ critical vertices spanning a $(9-1)!-1$ dimensional hyperface.
Thus, formally, the minimal presentations form a $47715-(9-1)!=7395$ dimensional abstract polytope.
The dihedrally invariant (dual) space is  2132 dimensional.
In particular, there is a minimal presentation with at most $2132$ dihedrally symmetrized Lie monomials.
\qedexer
\end{example}

\snewpage
\section{The BCH expansion and $\Gamma^{\Lie}_{n,m}$ in low degrees}\plabel{sec:someGamma}

In the literature, one can find several accounts for the concrete description of the terms of the BCH expansion.
The early investigations aim for a relatively straightforward algebraic-algorithmic presentation,
see Dynkin \cite{Dy},  Goldberg \cite{G};
they transform $\BCH^{\mathrm{alg}}(X_1,X_2)$ into $\BCH^{\Lie}(X_1,X_2)$
(but see \cite{L0} for minor variants of this procedure).
See Newman, Thompson \cite{NT} for a longer list, and references therein and thereto.
Later, lists concentrate on presenting the BCH terms in specific bases; see Michel \cite{Mi}, Macdonald \cite{Mac}, Casas, Murua \cite{CM}.
Thirdly, in order to obtain even more economical presentations, there are searches for  presentations
with given generating sets but subject to certain minimality conditions; see Oteo \cite{O}, Kolsrud \cite{K}.
For our purposes, however, the relevant presentations are those where the sum of the absolute coefficients of the Lie monomials is minimal,
i. e. they are minimal presentations in our  terminology.

One can apply to the same methods from linear programming as in the case of the Magnus expansion
in order to obtain minimal presentations.
We will consider (nonzero) Lie monomials to be the same if they give same free Lie algebra element, which is
the same as that they evaluate to the same free noncommutative algebra element by the commutator expansion.
If $M$ is a Lie monomial, then we also consider $-M$ as a Lie monomial.
Two Lie monomials are essentially the same if they differ only by sign.
We will consider two minimal representations to be the same if they are essentially the same, i. e. they can be
brought to the same form by changing the signs of the constituting monomials.
A small difference to the multiplicity-free Magnus case is that the accounting of the Lie monomials is not very predictable.
In the multiplicity-free case,  Lie monomials are essentially represented by  rooted unordered binary trees on the variables.
With multiplicities (as in the case of the BCH expansion)  the situation is more complicated, as the identity
\[X_{[1,[2,[1,2]]]}=-X_{[2,[1,[1,2]]]}\]
shows.
(In particular, the expressions $X_{[1,[2,[1,2]]]}, -X_{[2,[1,[1,2]]]},  \frac13 X_{[1,[2,[1,2]]]} -\frac23  X_{[2,[1,[1,2]]]}$
are essentially the same according to our terminology.)
It may also happen that a monomial is the fractional part of another, e. g.
\[X_{[[1, [1, [2, [2, [1, 2]]]]], [[1, 2], [1, [2, [1, 2]]]]]}=-\frac13X_{[[1, [1, [2, [2, [1, 2]]]]], [2, [2, [1, [1, [1, 2]]]]]]}.\]
This behaviour is expected, though; as any almost independent constellation of monomials, like
\[-{X}_{{[1,[1,[2,[2,[1,2]]]]]}}+{X}_{{[2,[2,[1,[1,[1,2]]]]]}}+3\,{X}_{{[[1,2],[1,[2,[1,2]]]]}}=0,\]
can be reduced further by taking a commutator with one of the monomials.
The identity above also excludes the expression ${X}_{{[[1,2],[1,[2,[1,2]]]]}}$ from any minimal presentation.
Practically, for any Lie monomial up to sign we fix, basically arbitrarily,   a single monomial representative to be used.
In our case, we will prefer $X_{[1,[2,[1,2]]]}$ to  $X_{[2,[1,[1,2]]]}$ or to, say, $X_{[2,[1,[2,1]]]}$,  just by lexicographical principles.

It is easy to see that
\[\BCH^{\Lie}(X_1,X_2)=X_1+X_2+\sum_{n_1=1}^\infty \sum_{n_2=1}^\infty \BCH^{\Lie}_{n_1,n_2}(X_1,X_2),\]
where $\BCH^{\Lie}_{n_1,n_2}(X_1,X_2)$ is homogeneous of degree $n_1$ in $X_1$ and  homogeneous of degree $n_2$ in $X_2$.
If $\BCH^{\Lie}_{n_1,n_2}(X_1,X_2)\neq0$, then $\Gamma_{n_1,n_2}^{\Lie}$ is the smallest nonnegative real such that
$\BCH^{\Lie}_{n_1,n_2}(X_1,X_2)/\Gamma_{n_1,n_2}^{\Lie}$ is in the convex hull of the Lie monomials.
In fact, $\BCH^{\Lie}_{n_1,n_2}(X_1,X_2)/\Gamma_{n_1,n_2}^{\Lie}$ is the best to be thought of as the place
where the positive ray of $\BCH^{\Lie}_{n_1,n_2}(X_1,X_2)$ intersects the boundary the convex hull of the Lie monomials.
(But, as have seen, not every Lie monomial is a vertex of the convex hull.)
Then, $\BCH^{\Lie}_{n_1,n_2}(X_1,X_2)/\Gamma_{n_1,n_2}^{\Lie}$ is in a face of the convex hull of the Lie monomials.
The smallest dimensional such face, for what $\BCH^{\Lie}_{n_1,n_2}(X_1,X_2)/\Gamma_{n_1,n_2}^{\Lie}$ is in its interior, is the critical face.
Then the minimal presentations are convex combinations of the Lie monomials lying on the critical face, i. e. of the  critical Lie monomials,
yielding  $\BCH^{\Lie}_{n_1,n_2}(X_1,X_2)$.
(The vertices of the critical face are also sufficient but they do not necessarily account for all critical Lie monomials.)
Relative to those critical Lie monomials, the minimal presentations form an abstract polytope,
whose vertices are the extremal minimal presentations.
It is not hard to see that a minimal presentation is extremal if and only if the set of its  Lie monomials cannot be reduced further.

The shape of a minimal presentation in itself is of little importance; one either prefers an extremal one which is
supposed to be relatively simple, or one of barycentric type which uses every critical Lie monomial.
In the next example, we look for extremal ones as we try to give a more comprehensive view of the (lack of) patterns
of possible shapes.
Due to the  easy-to-see identity
\[\BCH^{\Lie}_{n_1,n_2}(X_1,X_2)=(-1)^{n_1+n_2+1}\BCH^{\Lie}_{n_2,n_1}(X_2,X_1),\]
it is sufficient to consider the cases $n_1\leq n_2$.
\begin{example}\plabel{ex:data1}
Some accounts of minimal (extremal) representations for BCH terms are as follows.

Case $(1,1)$:
\[\BCH^{\Lie}_{1,1}(X_1,X_2)=\frac12\,{X}_{{[1,2]}}.\]

Case $(1,2)$:
\[\BCH^{\Lie}_{1,2}(X_1,X_2)=-\frac{1}{12}\,{X}_{{[2,[1,2]]}}.\]

Case $(1,3)$:
\[\BCH^{\Lie}_{1,3}(X_1,X_2)=0.\]

Case $(2,2)$:
\[\BCH^{\Lie}_{2,2}(X_1,X_2)=-\frac{1}{24}\,{X}_{{[1,[2,[1,2]]]}}.\]

Case $(1,4)$:
\[\BCH^{\Lie}_{1,4}(X_1,X_2)={\frac {{X}_{{[2,[2,[2,[1,2]]]]}}}{720}}.\]

Case $(2,3)$:
\[\BCH^{\Lie}_{2,3}(X_1,X_2)={\frac {{X}_{{[2,[1,[2,[1,2]]]]}}}{180}}
-{\frac {{X}_{{[[1,2],[2,[1,2]]]}}}{360}}
.\]

Case $(1,5)$:
\[\BCH^{\Lie}_{1,5}(X_1,X_2)=0.\]

Case $(2,4)$:
\[\BCH^{\Lie}_{2,4}(X_1,X_2)={\frac {{X}_{{[1,[2,[2,[2,[1,2]]]]]}}}{1440}}.\]

Case $(3,3)$:
\begin{align*}
\BCH_{3,3}^{\Lie}(X_1,X_2)&={\frac {{X}_{{[1,[2,[1,[2,[1,2]]]]]}}}{360}}
-{\frac {{X}_{{[1,[[1,2],[2,[1,2]]]]}}}{720}}
\\&={\frac {{X}_{{[2,[1,[1,[2,[1,2]]]]]}}}{360}}
+{\frac {{X}_{{[2,[[1,2],[1,[1,2]]]]}}}{720}}
\\&={\frac {{X}_{{[1,[2,[1,[2,[1,2]]]]]}}}{720}}
+{\frac {{X}_{{[2,[1,[1,[2,[1,2]]]]]}}}{720}}
-{\frac {{X}_{{[[1,[1,2]],[2,[1,2]]]}}}{720}}
\end{align*}
(i. e. 3 extremal minimal presentations altogether) spanning a $2$ dimensional abstract polytope, which is a triangle.
So, $6$ is the lowest degree where the minimal presentation is not unique.

\textit{Remark.} Here the last presentation is not the shortest but the most symmetrical.

Case $(1,6)$:
\[\BCH_{1,6}^{\Lie}(X_1,X_2)=-{\frac {{X}_{{[2,[2,[2,[2,[2,[1,2]]]]]]}}}{30240}}.\]

Case $(2,5)$:
\begin{align*}
\BCH_{2,5}^{\Lie}(X_1,X_2)&=-{\frac {{X}_{{[2,[2,[1,[2,[2,[1,2]]]]]]}}}{5040}}
+{\frac {{X}_{{[[1,2],[2,[2,[2,[1,2]]]]]}}}{10080}}
\\&=-{\frac {{X}_{{[2,[1,[2,[2,[2,[1,2]]]]]]}}}{10080}}
-{\frac {{X}_{{[2,[2,[2,[1,[2,[1,2]]]]]]}}}{10080}}
+{\frac {{X}_{{[[1,2],[2,[2,[2,[1,2]]]]]}}}{10080}}
\\&=-{\frac {{X}_{{[2,[1,[2,[2,[2,[1,2]]]]]]}}}{20160}}
-{\frac {{X}_{{[2,[2,[2,[1,[2,[1,2]]]]]]}}}{6720}}
-{\frac {{X}_{{[[2,[1,2]],[2,[2,[1,2]]]]}}}{10080}}
\\&=-{\frac {{X}_{{[2,[2,[1,[2,[2,[1,2]]]]]]}}}{10080}}
-{\frac {{X}_{{[2,[2,[2,[1,[2,[1,2]]]]]]}}}{10080}}
-{\frac {{X}_{{[[2,[1,2]],[2,[2,[1,2]]]]}}}{10080}}
\end{align*}
(i. e. 4 extremal minimal presentations altogether) spanning a $2$ dimensional abstract polytope, which is a trapezoid.

\textit{Remark.} Here the identity
\[{X}_{{[2,[2,[1,[2,[2,[1,2]]]]]]}}=\frac12\,{X}_{{[2,[1,[2,[2,[2,[1,2]]]]]]}}+\frac12\,{X}_{{[2,[2,[2,[1,[2,[1,2]]]]]]}}\]
exhibits a critical Lie monomial ``for the first time'', which is not an extremal point (vertex) of the critical face.
Notice, it contributes to the shortest  extremal minimal presentation of length $2$.

Case $(3,4)$:
\begin{align*}
\BCH_{3,4}^{\Lie}(X_1,X_2)
&
={\frac {{X}_{{[1,[1,[2,[2,[2,[1,2]]]]]]}}}{10080}}
-{\frac {{X}_{{[2,[2,[[1,2],[1,[1,2]]]]]}}}{7560}}
-{\frac {{X}_{{[2,[1,[2,[1,[2,[1,2]]]]]]}}}{5040}}\\&\qquad
+{\frac {{X}_{{[[1,2],[2,[1,[2,[1,2]]]]]}}}{2520}}
-{\frac {{X}_{{[2,[2,[1,[1,[2,[1,2]]]]]]}}}{6048}}
\\
&
={\frac {{X}_{{[1,[1,[2,[2,[2,[1,2]]]]]]}}}{10080}}
-{\frac {{X}_{{[2,[2,[[1,2],[1,[1,2]]]]]}}}{7560}}
-{\frac {{X}_{{[[2,[1,2]],[1,[2,[1,2]]]]}}}{5040}}\\&\qquad
+{\frac {{X}_{{[[1,2],[2,[1,[2,[1,2]]]]]}}}{5040}}
-{\frac {11\,{X}_{{[2,[2,[1,[1,[2,[1,2]]]]]]}}}{30240}}
\\
&
={\frac {{X}_{{[1,[1,[2,[2,[2,[1,2]]]]]]}}}{10080}}
-{\frac {{X}_{{[2,[2,[[1,2],[1,[1,2]]]]]}}}{20160}}
-{\frac {11\,{X}_{{[2,[1,[2,[1,[2,[1,2]]]]]]}}}{30240}}\\&\qquad
+{\frac {{X}_{{[[1,2],[2,[1,[2,[1,2]]]]]}}}{2520}}
+{\frac {{X}_{{[2,[1,[[1,2],[2,[1,2]]]]]}}}{12096}}
\\
&
={\frac {{X}_{{[1,[1,[2,[2,[2,[1,2]]]]]]}}}{10080}}
-{\frac {{X}_{{[2,[2,[[1,2],[1,[1,2]]]]]}}}{30240}}
-{\frac {{X}_{{[[2,[1,2]],[1,[2,[1,2]]]]}}}{2520}}\\&\qquad
+{\frac {{X}_{{[2,[1,[[1,2],[2,[1,2]]]]]}}}{10080}}
-{\frac {11\,{X}_{{[2,[2,[1,[1,[2,[1,2]]]]]]}}}{30240}}
\\
&
={\frac {{X}_{{[1,[1,[2,[2,[2,[1,2]]]]]]}}}{10080}}
+{\frac {{X}_{{[2,[[1,[1,2]],[2,[1,2]]]]}}}{7560}}
-{\frac {{X}_{{[2,[1,[2,[1,[2,[1,2]]]]]]}}}{3024}}\\&\qquad
+{\frac {{X}_{{[[1,2],[2,[1,[2,[1,2]]]]]}}}{2520}}
-{\frac {{X}_{{[2,[2,[1,[1,[2,[1,2]]]]]]}}}{30240}}
\\
&
={\frac {{X}_{{[1,[1,[2,[2,[2,[1,2]]]]]]}}}{10080}}
+{\frac {{X}_{{[2,[[1,[1,2]],[2,[1,2]]]]}}}{7560}}
-{\frac {{X}_{{[[2,[1,2]],[1,[2,[1,2]]]]}}}{3024}}\\&\qquad
+{\frac {{X}_{{[[1,2],[2,[1,[2,[1,2]]]]]}}}{15120}}
-{\frac {11\,{X}_{{[2,[2,[1,[1,[2,[1,2]]]]]]}}}{30240}}
\\
&
={\frac {{X}_{{[1,[1,[2,[2,[2,[1,2]]]]]]}}}{10080}}
+{\frac {{X}_{{[2,[[1,[1,2]],[2,[1,2]]]]}}}{10080}}
-{\frac {11\,{X}_{{[2,[1,[2,[1,[2,[1,2]]]]]]}}}{30240}}\\&\qquad
+{\frac {{X}_{{[[1,2],[2,[1,[2,[1,2]]]]]}}}{2520}}
+{\frac {{X}_{{[2,[1,[[1,2],[2,[1,2]]]]]}}}{30240}}
\\
&
={\frac {{X}_{{[1,[1,[2,[2,[2,[1,2]]]]]]}}}{10080}}
+{\frac {{X}_{{[2,[[1,[1,2]],[2,[1,2]]]]}}}{15120}}
-{\frac {{X}_{{[[2,[1,2]],[1,[2,[1,2]]]]}}}{2520}}\\&\qquad
+{\frac {{X}_{{[2,[1,[[1,2],[2,[1,2]]]]]}}}{15120}}
-{\frac {11\,{X}_{{[2,[2,[1,[1,[2,[1,2]]]]]]}}}{30240}}
\\
&
={\frac {{X}_{{[1,[1,[2,[2,[2,[1,2]]]]]]}}}{10080}}
-{\frac {{X}_{{[[2,[1,2]],[1,[2,[1,2]]]]}}}{10080}}
-{\frac {11\,{X}_{{[2,[1,[2,[1,[2,[1,2]]]]]]}}}{30240}}\\&\qquad
+{\frac {{X}_{{[[1,2],[2,[1,[2,[1,2]]]]]}}}{3360}}
+{\frac {{X}_{{[2,[1,[[1,2],[2,[1,2]]]]]}}}{7560}}
\\
&
={\frac {{X}_{{[1,[1,[2,[2,[2,[1,2]]]]]]}}}{10080}}
-{\frac {{X}_{{[[2,[1,2]],[1,[2,[1,2]]]]}}}{2520}}
-{\frac {{X}_{{[2,[1,[2,[1,[2,[1,2]]]]]]}}}{15120}}\\&\qquad
+{\frac {{X}_{{[2,[1,[[1,2],[2,[1,2]]]]]}}}{7560}}
-{\frac {{X}_{{[2,[2,[1,[1,[2,[1,2]]]]]]}}}{3360}}
\end{align*}
(i. e. 10 extremal minimal presentations of length 5 altogether) spanning a $3$ dimensional abstract polytope.

Case $(1,7)$:
\[\BCH_{1,7}^{\Lie}(X_1,X_2)=0.\]

Case $(2,6)$:
\[\BCH_{2,6}^{\Lie}(X_1,X_2)=-{\frac {{X}_{{[1,[2,[2,[2,[2,[2,[1,2]]]]]]]}}}{60480}}.\]

Case $(3,5)$:
\begin{align*}
\BCH_{3,5}^{\Lie}(X_1,X_2)&=-{\frac {{X}_{{[1,[2,[2,[1,[2,[2,[1,2]]]]]]]}}}{10080}}
+{\frac {{X}_{{[1,[[1,2],[2,[2,[2,[1,2]]]]]]}}}{20160}}
\\&=-{\frac {{X}_{{[1,[2,[1,[2,[2,[2,[1,2]]]]]]]}}}{20160}}
-{\frac {{X}_{{[1,[2,[2,[2,[1,[2,[1,2]]]]]]]}}}{20160}}
+{\frac {{X}_{{[1,[[1,2],[2,[2,[2,[1,2]]]]]]}}}{20160}}
\\&=-{\frac {{X}_{{[1,[2,[1,[2,[2,[2,[1,2]]]]]]]}}}{40320}}
-{\frac {{X}_{{[1,[2,[2,[2,[1,[2,[1,2]]]]]]]}}}{13440}}
-{\frac {{X}_{{[1,[[2,[1,2]],[2,[2,[1,2]]]]]}}}{20160}}
\\&=-{\frac {{X}_{{[1,[2,[2,[1,[2,[2,[1,2]]]]]]]}}}{20160}}
-{\frac {{X}_{{[1,[2,[2,[2,[1,[2,[1,2]]]]]]]}}}{20160}}
-{\frac {{X}_{{[1,[[2,[1,2]],[2,[2,[1,2]]]]]}}}{20160}}
\\&=-{\frac {{X}_{{[2,[1,[1,[2,[2,[2,[1,2]]]]]]]}}}{20160}}
-{\frac {{X}_{{[1,[2,[2,[2,[1,[2,[1,2]]]]]]]}}}{20160}}
+{\frac {{X}_{{[[1,[1,2]],[2,[2,[2,[1,2]]]]]}}}{20160}}
\end{align*}
(i. e. 5 extremal minimal presentations altogether) spanning a $3$ dimensional abstract polytope, which is a trapezoid based pyramid.

Case $(4,4)$:
\begin{align*}
\BCH_{4,4}^{\Lie}(X_1,X_2)
&={\frac {{X}_{{[1,[2,[[1,[1,2]],[2,[1,2]]]]]}}}{20160}}
-{\frac {{X}_{{[1,[[1,2],[[1,2],[2,[1,2]]]]]}}}{60480}}\\&\qquad\qquad\qquad
-{\frac {23\,{X}_{{[1,[2,[1,[2,[1,[2,[1,2]]]]]]]}}}{120960}}
+{\frac {23\,{X}_{{[1,[[1,2],[2,[1,[2,[1,2]]]]]]}}}{120960}}
\\&=
{\frac {{X}_{{[1,[2,[[1,[1,2]],[2,[1,2]]]]]}}}{20160}}
-{\frac {{X}_{{[1,[[1,2],[[1,2],[2,[1,2]]]]]}}}{60480}}\\&\qquad\qquad\qquad
-{\frac {23\,{X}_{{[1,[2,[2,[1,[1,[2,[1,2]]]]]]]}}}{120960}}
-{\frac {23\,{X}_{{[1,[[2,[1,2]],[1,[2,[1,2]]]]]}}}{120960}}
\\&=
{\frac {{X}_{{[1,[2,[[1,[1,2]],[2,[1,2]]]]]}}}{20160}}
-{\frac {{X}_{{[1,[[1,2],[[1,2],[2,[1,2]]]]]}}}{60480}}\\&\qquad\qquad\qquad
-{\frac {23\,{X}_{{[2,[1,[1,[2,[1,[2,[1,2]]]]]]]}}}{120960}}
+{\frac {23\,{X}_{{[2,[[1,[1,2]],[1,[2,[1,2]]]]]}}}{120960}}
\\&=
{\frac {{X}_{{[1,[2,[[1,[1,2]],[2,[1,2]]]]]}}}{20160}}
-{\frac {{X}_{{[1,[[1,2],[[1,2],[2,[1,2]]]]]}}}{60480}}\\&\qquad\qquad\qquad
-{\frac {23\,{X}_{{[2,[1,[2,[1,[1,[2,[1,2]]]]]]]}}}{120960}}
-{\frac {23\,{X}_{{[2,[[1,2],[1,[1,[2,[1,2]]]]]]}}}{120960}}
\\&=
{\frac {{X}_{{[2,[1,[[1,[1,2]],[2,[1,2]]]]]}}}{20160}}
-{\frac {{X}_{{[2,[[1,2],[[1,2],[1,[1,2]]]]]}}}{60480}}\\&\qquad\qquad\qquad
-{\frac {23\,{X}_{{[1,[2,[1,[2,[1,[2,[1,2]]]]]]]}}}{120960}}
+{\frac {23\,{X}_{{[1,[[1,2],[2,[1,[2,[1,2]]]]]]}}}{120960}}
\\&=
{\frac {{X}_{{[2,[1,[[1,[1,2]],[2,[1,2]]]]]}}}{20160}}
-{\frac {{X}_{{[2,[[1,2],[[1,2],[1,[1,2]]]]]}}}{60480}}\\&\qquad\qquad\qquad
-{\frac {23\,{X}_{{[1,[2,[2,[1,[1,[2,[1,2]]]]]]]}}}{120960}}
-{\frac {23\,{X}_{{[1,[[2,[1,2]],[1,[2,[1,2]]]]]}}}{120960}}
\\&=
{\frac {{X}_{{[2,[1,[[1,[1,2]],[2,[1,2]]]]]}}}{20160}}
-{\frac {{X}_{{[2,[[1,2],[[1,2],[1,[1,2]]]]]}}}{60480}}\\&\qquad\qquad\qquad
-{\frac {23\,{X}_{{[2,[1,[1,[2,[1,[2,[1,2]]]]]]]}}}{120960}}
+{\frac {23\,{X}_{{[2,[[1,[1,2]],[1,[2,[1,2]]]]]}}}{120960}}
\\&=
{\frac {{X}_{{[2,[1,[[1,[1,2]],[2,[1,2]]]]]}}}{20160}}
-{\frac {{X}_{{[2,[[1,2],[[1,2],[1,[1,2]]]]]}}}{60480}}\\&\qquad\qquad\qquad
-{\frac {23\,{X}_{{[2,[1,[2,[1,[1,[2,[1,2]]]]]]]}}}{120960}}
-{\frac {23\,{X}_{{[2,[[1,2],[1,[1,[2,[1,2]]]]]]}}}{120960}}
\\&\qquad\text{(which are 8 extremal minimal presentations of length 4)}
\\&=16\text{ extremal minimal presentations of length }5
\\&=36\text{ extremal minimal presentations of length }6
\\&=28\text{ extremal minimal presentations of length }7
\\&=64\text{ extremal minimal presentations of length }8
\end{align*}
(i. e. 152 extremal minimal presentations altogether) spanning a $10$ dimensional abstract polytope.

Case $(1,8)$:
\[\BCH_{1,8}^{\Lie}(X_1,X_2)={\frac {{X}_{{[2,[2,[2,[2,[2,[2,[2,[1,2]]]]]]]]}}}{1209600}}.\]

Case $(2,7)$:
\begin{align*}
&\BCH_{2,7}^{\Lie}(X_1,X_2)=
\\&={\frac {{X}_{{[2,[1,[2,[2,[2,[2,[2,[1,2]]]]]]]]}}}{302400}}
+{\frac {{X}_{{[2,[2,[2,[2,[2,[1,[2,[1,2]]]]]]]]}}}{302400}}-{\frac {{X}_{{[[1,2],[2,[2,[2,[2,[2,[1,2]]]]]]]}}}{302400}}
\\&={\frac {{X}_{{[2,[1,[2,[2,[2,[2,[2,[1,2]]]]]]]]}}}{907200}}
+{\frac {{X}_{{[2,[2,[2,[2,[1,[2,[2,[1,2]]]]]]]]}}}{181440}}+{\frac {{X}_{{[[2,[1,2]],[2,[2,[2,[2,[1,2]]]]]]}}}{302400}}
\\&={\frac {{X}_{{[2,[2,[1,[2,[2,[2,[2,[1,2]]]]]]]]}}}{302400}}
+{\frac {{X}_{{[2,[2,[2,[2,[2,[1,[2,[1,2]]]]]]]]}}}{302400}}+{\frac {{X}_{{[[2,[1,2]],[2,[2,[2,[2,[1,2]]]]]]}}}{302400}}
\\&={\frac {{X}_{{[2,[2,[2,[1,[2,[2,[2,[1,2]]]]]]]]}}}{201600}}
+{\frac {{X}_{{[2,[2,[2,[2,[2,[1,[2,[1,2]]]]]]]]}}}{604800}}-{\frac {{X}_{{[[2,[2,[1,2]]],[2,[2,[2,[1,2]]]]]}}}{302400}}
\\&={\frac {{X}_{{[2,[2,[2,[1,[2,[2,[2,[1,2]]]]]]]]}}}{302400}}
+{\frac {{X}_{{[2,[2,[2,[2,[1,[2,[2,[1,2]]]]]]]]}}}{302400}}-{\frac {{X}_{{[[2,[2,[1,2]]],[2,[2,[2,[1,2]]]]]}}}{302400}}
\\&\qquad\text{(which are 5 extremal minimal presentations of length 3)}
\\&=5\text{ extremal minimal presentations of length }4
\end{align*}
(i. e. 10 extremal minimal presentations altogether) spanning a $4$ dimensional abstract polytope.

Case $(3,6)$:
\begin{align*}
&\BCH_{3,6}^{\Lie}(X_1,X_2)=
\\&\quad-{\frac {{X}_{{[[1,2],[2,[2,[1,[2,[2,[1,2]]]]]]]}}}{75600}}
+{\frac {{X}_{{[2,[1,[2,[2,[2,[1,[2,[1,2]]]]]]]]}}}{100800}}
-{\frac {{X}_{{[2,[1,[[1,2],[2,[2,[2,[1,2]]]]]]]}}}{151200}}
\\&\quad+{\frac {{X}_{{[2,[2,[2,[2,[[1,2],[1,[1,2]]]]]]]}}}{453600}}
+{\frac {{X}_{{[2,[2,[1,[1,[2,[2,[2,[1,2]]]]]]]]}}}{604800}}
+{\frac {{X}_{{[[2,[1,2]],[1,[2,[2,[2,[1,2]]]]]]}}}{100800}}
\\&\quad+{\frac {11\,{X}_{{[2,[2,[2,[2,[1,[1,[2,[1,2]]]]]]]]}}}{1814400}}
\qquad\text{(which is of length 7)}
\\&=129\text{ extremal minimal presentations of length }8
\\&=4943\text{ extremal minimal presentations of length }9
\end{align*}
(i. e. 5073 extremal minimal presentations altogether) spanning a $21$ dimensional abstract polytope.

Case $(4,5)$:
\begin{align*}
&\BCH_{4,5}^{\Lie}(X_1,X_2)=
\\&=28\text{ extremal minimal presentations of length 13; one of them is}
\\&\quad-{\frac {{X}_{{[[1,[2,[1,2]]],[2,[1,[2,[1,2]]]]]}}}{40320}}
-{\frac {{X}_{{[1,[1,[2,[1,[2,[2,[2,[1,2]]]]]]]]}}}{1814400}}
+{\frac {{X}_{{[2,[1,[2,[2,[1,[1,[2,[1,2]]]]]]]]}}}{72576}}
\\&\quad+{\frac {{X}_{{[2,[2,[1,[1,[2,[1,[2,[1,2]]]]]]]]}}}{120960}}
-{\frac {{X}_{{[1,[1,[2,[2,[2,[1,[2,[1,2]]]]]]]]}}}{72576}}
+{\frac {{X}_{{[1,[1,[[1,2],[2,[2,[2,[1,2]]]]]]]}}}{453600}}
\\&\quad-{\frac {{X}_{{[1,[2,[1,[1,[2,[2,[2,[1,2]]]]]]]]}}}{181440}}
+{\frac {{X}_{{[2,[1,[2,[1,[2,[1,[2,[1,2]]]]]]]]}}}{90720}}
+{\frac {{X}_{{[2,[1,[[2,[1,2]],[1,[2,[1,2]]]]]]}}}{120960}}
\\&\quad+{\frac {{X}_{{[2,[2,[[1,2],[[1,2],[1,[1,2]]]]]]}}}{181440}}
-{\frac {{X}_{{[[1,2],[2,[1,[2,[1,[2,[1,2]]]]]]]}}}{181440}}
-{\frac {{X}_{{[[1,2],[2,[2,[1,[1,[2,[1,2]]]]]]]}}}{90720}}
\\&\quad-{\frac {{X}_{{[[1,2],[[1,2],[[1,2],[2,[1,2]]]]]}}}{90720}}
\\&=216\text{ extremal minimal presentations of length }14
\end{align*}
(i. e. 244 extremal minimal presentations altogether) spanning a $11$ dimensional abstract polytope.

Case $(1,9)$:
\[\BCH_{1,9}^{\Lie}(X_1,X_2)=0.\]

Case $(2,8)$:
\[\BCH_{2,8}^{\Lie}(X_1,X_2)={\frac {{X}_{{[1,[2,[2,[2,[2,[2,[2,[2,[1,2]]]]]]]]]}}}{2419200}}.\]

Case $(3,7)$:
\begin{align*}
&\BCH_{3,7}^{\Lie}(X_1,X_2)=
\\&={\frac {{X}_{{[1,[2,[1,[2,[2,[2,[2,[2,[1,2]]]]]]]]]}}}{604800}}+
{\frac {{X}_{{[1,[2,[2,[2,[2,[2,[1,[2,[1,2]]]]]]]]]}}}{604800}}-
{\frac {{X}_{{[1,[[1,2],[2,[2,[2,[2,[2,[1,2]]]]]]]]}}}{604800}}
\\&={\frac {{X}_{{[1,[2,[1,[2,[2,[2,[2,[2,[1,2]]]]]]]]]}}}{1814400}}
+{\frac {{X}_{{[1,[2,[2,[2,[2,[1,[2,[2,[1,2]]]]]]]]]}}}{362880}}
+{\frac {{X}_{{[1,[[2,[1,2]],[2,[2,[2,[2,[1,2]]]]]]]}}}{604800}}
\\&={\frac {{X}_{{[1,[2,[2,[1,[2,[2,[2,[2,[1,2]]]]]]]]]}}}{604800}}+
{\frac {{X}_{{[1,[2,[2,[2,[2,[2,[1,[2,[1,2]]]]]]]]]}}}{604800}}+
{\frac {{X}_{{[1,[[2,[1,2]],[2,[2,[2,[2,[1,2]]]]]]]}}}{604800}}
\\&={\frac {{X}_{{[1,[2,[2,[2,[1,[2,[2,[2,[1,2]]]]]]]]]}}}{403200}}+
{\frac {{X}_{{[1,[2,[2,[2,[2,[2,[1,[2,[1,2]]]]]]]]]}}}{1209600}}
-{\frac {{X}_{{[1,[[2,[2,[1,2]]],[2,[2,[2,[1,2]]]]]]}}}{604800}}
\\&={\frac {{X}_{{[1,[2,[2,[2,[1,[2,[2,[2,[1,2]]]]]]]]]}}}{604800}}+
{\frac {{X}_{{[1,[2,[2,[2,[2,[1,[2,[2,[1,2]]]]]]]]]}}}{604800}}-
{\frac {{X}_{{[1,[[2,[2,[1,2]]],[2,[2,[2,[1,2]]]]]]}}}{604800}}
\\&={\frac {{X}_{{[2,[1,[1,[2,[2,[2,[2,[2,[1,2]]]]]]]]]}}}{604800}}
+{\frac {{X}_{{[1,[2,[2,[2,[2,[2,[1,[2,[1,2]]]]]]]]]}}}{604800}}
-{\frac {{X}_{{[[1,[1,2]],[2,[2,[2,[2,[2,[1,2]]]]]]]}}}{604800}}
\\&\quad\text{(which are 6 extremal minimal presentations of length 3)}
\\&=6\text{ extremal minimal presentations of length }4
\end{align*}
(i. e. 12 extremal minimal presentations altogether) spanning a $5$ dimensional abstract polytope.

Case $(4,6)$:
\begin{align*}
&\BCH_{4,6}^{\Lie}(X_1,X_2)=
\\&=25\text{ extremal minimal presentations of length 7, one of those is}
\\&\qquad{\frac {{X}_{{[1,[2,[[2,[1,2]],[[1,2],[2,[1,2]]]]]]}}}{302400}}
-{\frac {{X}_{{[1,[[1,2],[2,[1,[2,[2,[2,[1,2]]]]]]]]}}}{725760}}
-{\frac {{X}_{{[1,[[1,2],[2,[2,[2,[1,[2,[1,2]]]]]]]]}}}{302400}}
\\&\qquad-{\frac {{X}_{{[1,[[2,[2,[1,2]]],[2,[1,[2,[1,2]]]]]]}}}{181440}}
+{\frac {{X}_{{[1,[2,[1,[2,[1,[2,[2,[2,[1,2]]]]]]]]]}}}{518400}}
+{\frac {{X}_{{[1,[2,[2,[1,[2,[2,[1,[2,[1,2]]]]]]]]]}}}{120960}}
\\&\qquad-{\frac {{X}_{{[1,[2,[2,[2,[1,[[1,2],[2,[1,2]]]]]]]]}}}{907200}}
\\&=\text{ other extremal minimal presentations of length }\geq8
\end{align*}
spanning a 33 dimensional abstract polytope.

Case $(5,5)$:
\begin{align*}
&\BCH_{5,5}^{\Lie}(X_1,X_2)=
\\&={\frac {{X}_{{[1,[2,[1,[2,[1,[2,[1,[2,[1,2]]]]]]]]]}}}{1209600}}
+{\frac {{X}_{{[1,[2,[1,[[1,2],[[1,2],[2,[1,2]]]]]]]}}}{604800}}
+{\frac {{X}_{{[1,[2,[1,[[2,[1,2]],[1,[2,[1,2]]]]]]]}}}{604800}}
\\&\quad+{\frac {{X}_{{[1,[2,[2,[1,[1,[2,[1,[2,[1,2]]]]]]]]]}}}{151200}}
+{\frac {{X}_{{[1,[2,[2,[1,[2,[1,[1,[2,[1,2]]]]]]]]]}}}{134400}}
+{\frac {{X}_{{[1,[2,[[2,[1,2]],[[1,2],[1,[1,2]]]]]]}}}{201600}}
\\&\quad+{\frac {{X}_{{[1,[[1,2],[[1,2],[2,[1,[2,[1,2]]]]]]]}}}{1209600}}
+{\frac {{X}_{{[1,[[2,[1,2]],[1,[2,[1,[2,[1,2]]]]]]]}}}{134400}}
+{\frac {{X}_{{[1,[[2,[1,2]],[2,[1,[1,[2,[1,2]]]]]]]}}}{134400}}
\\&\quad+{\frac {{X}_{{[1,[[2,[1,2]],[2,[[1,2],[1,[1,2]]]]]]}}}{1209600}}
-{\frac {{X}_{{[1,[2,[2,[1,[1,[[1,2],[2,[1,2]]]]]]]]}}}{1209600}}
-{\frac {{X}_{{[1,[[1,2],[[1,2],[[1,2],[2,[1,2]]]]]]}}}{604800}}
\\&\quad-{\frac {{X}_{{[1,[[1,[2,[1,2]]],[2,[1,[2,[1,2]]]]]]}}}{134400}}
\\&=\text{ other extremal minimal presentations. }
\end{align*}
spanning a 61 dimensional abstract polytope.
\qedexer
\end{example}
Verifying the previous example requires much computation.
However, the part that we see minimal presentations is easy to demonstrate.
First, in order show that we have a presentation for $\BCH_{n_1,n_2}^{\Lie}(X_1,X_2)$,
it is sufficient to expand it in the noncommutative polynomial algebra
in order to see that it gives the algebraic $\BCH^{\mathrm{alg}}_{n_1,n_2}(X_1,X_2)$.
In order to see minimality, it is sufficient to exhibit a linear functional $\alpha_{n_1,n_2}$
which takes values on the Lie monomials from $[-1,1]$ (i.~e.~it is a face-defining linear functional,
at this point allowing to define the $-1$ dimensional face),
but which takes $1$ on the constituting Lie monomials of the given presentation with the sign of the coefficients incorporated.
Or, equivalently,  $\alpha_{n_1,n_2}(\BCH^{\Lie}_{n_1,n_2}(X_1,X_2))$ is equal to the
sum of the absolute values of the coefficients in the given presentation; or,
the defined face contains $\BCH_{n_1,n_2}^{\Lie}(X_1,X_2)/\Gamma_{n_1,n_2}^{\Lie}$; or, the face-defining linear
functional is maximal with respect to $\BCH_{n_1,n_2}^{\Lie}(X_1,X_2)$.
(This latter one is the dual LP problem.)

We can generate linear functionals on a  free Lie algebra as follows.
Let $X^*_{i_1,\ldots,i_n}$ be the linear functional which gives the coefficient of $X_{i_1}\ldots X_{i_n}$.
In what follows we use these linear functionals as restricted the Lie algebra.
Further linear functionals (on the $(n_1,n_2)$-homogeneous part of the free Lie algebra) can be constructed by linear combinations.
In fact, this makes the description of the linear functionals highly non-unique.
(For example, $X^*_{i_1,\ldots,i_n}=(-1)^{n-1}X^*_{i_n,\ldots,i_1}$ always holds;
or,  the $(k,1)$-homogeneous part of the free Lie algebra is $1$ dimensional,
but we can use several monomials, etc.)
On the other hand, there is little benefit in using any prearranged basis for the linear functionals.

The face-defining linear functionals which take $1$ on $\BCH^{\Lie}_{n_1,n_2}(X_1,X_2)/\Gamma_{n_1,n_2}^{\Lie}$
(that is the dual solutions)
form a polytope itself. Its vertices (extremal points) are the relevant hyperface-defining linear functionals,
and its interior points are linear functionals defining the critical face.
Again, one may prefer one type to the other, but we will seek out hyperface-defining linear functionals;
some of those may be quite simple.
Such linear functionals can be unique (if the critical face is hyperface itself) or not.
The sort of simplest linear functionals are integral, i. e. they take values on  the Lie monomials
from the set $\{-1,0,1\}=[-1,1]\cap\mathbb Z$.

\begin{example}\plabel{ex:data2}
 Here we exhibit some hyperface-defining linear functionals which take
$1$ on $\BCH_{n_1,n_2}(X_1,X_2)/\Gamma_{n_1,n_2}^{\Lie}$.
We will give comments regarding the properties of the given hyperface-defining linear functional.
We may give alternative presentations for the chosen linear function; nevertheless in each case
we give only one single linear function.

Case $(1,1)$: (unique, integral)
\[\alpha_{1,1}=X^*_{{12}}.\]

Case $(1,2)$: (unique, integral)
\[\alpha_{1,2}=X^*_{{122}}.\]

Case $(2,2)$: (unique, integral)
\[\alpha_{2,2}= X^*_{{1122}} = -\frac12 X^*_{{1212}}  .\]

Case $(1,4)$: (unique, integral)
\[\alpha_{1,4}=  -X^*_{{12222}}  .\]

Case $(2,3)$: (unique, integral)
\[\alpha_{2,3}= -X^*_{{21122}}   .\]

Case $(2,4)$: (non-unique, necessarily integral)
\[\alpha_{2,4}= -X^*_{{112222}} .\]

Case $(3,3)$: (unique, integral)
\[\alpha_{3,3}= -X^*_{{122112}}=\frac14 X^*_{{121212}} .\]

Case $(1,6)$: (unique, integral)
\[\alpha_{1,6}= X^*_{{1222222}}  .\]

Case $(2,5)$: (unique, integral)
\[\alpha_{2,5}= X^*_{{2112222}}+X^*_{{2211222}}  =X^*_{{2221122}}+X^*_{{2222112}}.\]

Case $(3,4)$: (unique, non-integral)
\[\alpha_{3,4}=  -X^*_{{1221122}}+\frac{X^*_{{1221221}}}{12} =-\frac{X^*_{{1122122}}}2-\frac{X^*_{{1122212}}}2-\frac{X^*_{{1221122}}}2
 .\]

Case $(2,6)$: (non-unique, unique integral)
\[\alpha_{2,6}=  X^*_{{11222222}}  .\]

Case $(3,5)$: (non-unique, unique integral)
\[\alpha_{3,5}= X^*_{{12211222}}+X^*_{{12222112}}   .\]

Case $(4,4)$: (unique, integral)
\[\alpha_{4,4}= -X^*_{{11221122}} =-\frac18 X^*_{{12121212}}  .\]

Case $(1,8)$: (unique, integral)
\[\alpha_{1,8}= -X^*_{{122222222}} .\]

Case $(2,7)$: (unique, integral)
\[\alpha_{2,7}=  X^*_{{112222222}}+\frac{X^*_{{122222221}}}2
=-X^*_{{222211222}}-X^*_{{222221122}}-X^*_{{222222112}}.\]

Case $(3,6)$: (unique, integral)
\[\alpha_{3,6}=   X^*_{{122112222}}+X^*_{{122221122}} .\]

Case $(4,5)$: (unique, non-integral)
\[\alpha_{4,5}=
\frac{X^*_{{112222112}}}{3}+\frac{X^*_{{121222121}}}{24}-\frac{X^*_{{211122212}}}{12}-{\frac {5\,X^*_{{211212212}}}{108}}
+{\frac {53\,X^*_{{211221122}}}{54}}+{\frac {X^*_{{221121122}}}{108}}\]
\[=-\frac{X^*_{{121221221}}}{12}-\frac{X^*_{{122121221}}}{24}
-\frac{X^*_{{211122212}}}{12}-{\frac {5\,X^*_{{211212212}}}{108}}
+{\frac {53\,X^*_{{211221122}}}{54}}+{\frac {X^*_{{221121122}}}{108}}\]
\[=\frac{X^*_{{112222112}}}{3}+\frac{X^*_{{121222121}}}{24}-\frac{X^*_{{211112222}}}{3}-\frac{X^*_{{211121222}}}{6}
-\frac{X^*_{{211211222}}}{12}+X^*_{{211221122}}-\frac{X^*_{{211222112}}}{24}.\]

Case $(2,8)$: (non-unique, ``the simplest'' integral one)
\[\alpha_{2,8}=  -X^*_{{1122222222}}  .\]

Case $(3,7)$: (non-unique, unique integral)
\[\alpha_{3,7}=   -X^*_{{1221122222}}-X^*_{{1222211222}}-X^*_{{1222222112}}
 .\]

Case $(4,6)$: (non-unique, unique integral)
\[\alpha_{4,6}=  X^*_{{1122112222}}+X^*_{{1122221122}}  .\]

Case $(5,5)$: (unique, integral)
\[\alpha_{5,5}= X^*_{{1221122112}}   =\frac1{16} X^*_{{1212121212}} .\]
\qedexer
\end{example}
For $n=n_1+n_2\geq11$, expressions for minimal presentations
and good face-defining linear functionals became more complicated.
Nevertheless, one can obtain minimal presenations for a couple of other degrees.
This yields, for $n=n_1+n_2\leq16$,

\begin{example}\plabel{ex:data3}
\[\Gamma^{\Lie}_1(x,y)=x+y;\]

\[\Gamma^{\Lie}_2(x,y)=\frac12\,xy;\]

\[\Gamma^{\Lie}_3(x,y)=\frac1{12}\,({x}^{2}y+x{y}^{2});\]

\[\Gamma^{\Lie}_4(x,y)=\frac1{24}\,{x}^{2}{y}^{2};\]

\[\Gamma^{\Lie}_5(x,y)={\frac {1}{720}}({x}^{4}y+x{y}^{4})
+{\frac {1}{120}}({x}^{3} {y}^{2}+{x}^{2}{y}^{3});\]

\[\Gamma^{\Lie}_6(x,y)={\frac {1}{1440}}({x}^{4}{y}^{2}+{x}^{2}{y}^{4})
+{\frac {1}{240}}{x}^{3}{y}^{3};\]

\[\Gamma^{\Lie}_7(x,y)={\frac {1}{30240}}({x}^{6}y+x{y}^{6})
+{\frac {1}{3360}({x}^{5}{y}^{2}}+{x}^{2}{y}^{5})
+{\frac {1}{1008}}({x}^{4}{y}^{3}+{x}^{3}{y}^{4});\]

\[\Gamma^{\Lie}_8(x,y)={\frac {1}{60480}}({x}^{6}{y}^{2}+{x}^{2}{y}^{6})+{\frac {1}{6720}}({x}^{5}{y}^{3}+{x}^{3}{y}^{5})
+{\frac{1}{2240}}{x}^{4}{y}^{4};\]

\begin{align*}
\Gamma^{\Lie}_9(x,y)=\,&{\frac {1}{1209600}}({x}^{8}y+x{y}^{8})
+{\frac {1}{100800}}({x}^{7}{y}^{2}+{x}^{2}{y}^{7})
+{\frac {1}{20160}}({x}^{6}{y}^{3}+{x}^{3}{y}^{6})
\\[3mm]&
+{\frac {11}{90720}}({x}^{5}{y}^{4}+{x}^{4}{y}^{5});
\end{align*}

\begin{align*}
\Gamma^{\Lie}_{10}(x,y)=\,&{\frac {1}{2419200}}({x}^{8}{y}^{2}+{x}^{2}{y}^{8})
+{\frac {1}{201600}}({x}^{7}{y}^{3}+{x}^{3}{y}^{7})
+{\frac {1}{40320}}({x}^{6}{y}^{4}+{x}^{4}{y}^{6})
\\[3mm]&
+{\frac {1}{20160}}{x}^{5}{y}^{5};
\end{align*}

\begin{align*}
\Gamma^{\Lie}_{11}(x,y)=\,&{\frac {1}{47900160}}({x}^{10}y+x{y}^{10})
+{\frac {1}{3193344}}({x}^{9}{y}^{2}+{x}^{2}{y}^{9})
+{\frac {19}{8870400}}({x}^{8}{y}^{3}+{x}^{3}{y}^{8})
\\[3mm]&+{\frac {4573}{585446400}}({x}^{7}{y}^{4}+{x}^{4}{y}^{7})
+{\frac {12176029}{782632287360}}({x}^{6}{y}^{5}+{x}^{5}{y}^{6});
\end{align*}

\begin{align*}
\Gamma^{\Lie}_{12}(x,y)=\,&{\frac {1}{95800320}}({x}^{10}{y}^{2}+{x}^{2}{y}^{10})
+{\frac {1}{6386688}}({x}^{9}{y}^{3}+{x}^{3}{y}^{9})
+{\frac {19}{17740800}}({x}^{8}{y}^{4}+{x}^{4}{y}^{8})
\\[3mm]&+{\frac {1}{266112}}({x}^{7}{y}^{5}+{x}^{5}{y}^{7})
+{\frac {4759889}{810949708800}}{x}^{6}{y}^{6};
\end{align*}

\begin{align*}
\Gamma^{\Lie}_{13}(x,y)=\,&{\frac {691}{1307674368000}}({x}^{12}y+x{y}^{12})
\\[3mm]&+{\frac {691}{72648576000}}({x}^{11}{y}^{2}+{x}^{2}{y}^{11})
\\[3mm]&+{\frac {397}{4843238400}}({x}^{10}{y}^{3}+{x}^{3}{y}^{10})
\\[3mm]&+{\frac {2394627036176344242739}{5942332300516463608428672000}}({x}^{9}{y}^{4}+{x}^{4}{y}^{9})
\\[3mm]&+{\frac {50716867585528101162848154608067859}{
42820360683492185111833447903002621388800}}({x}^{8}{y}^{5}+{x}^{5}{y}^{8})
\\[3mm]&+{\frac {6351369008653149273148491526516654975706237}{
3088991399411560045433904125019163555602524928000}}({x}^{7}{y}^{6}+{x}^{6}{y}^{7});
\end{align*}

\begin{align*}
\Gamma^{\Lie}_{14}(x,y)=\,&
{\frac {691}{2615348736000}}({x}^{12}{y}^{2}+{x}^{2}{y}^{12})
\\[3mm]&+{\frac {691}{145297152000}}({x}^{11}{y}^{3}+{x}^{3}{y}^{11})
\\[3mm]&+{\frac {397}{9686476800}}({x}^{10}{y}^{4}+{x}^{4}{y}^{10})
\\[3mm]&+{\frac {83}{415134720}}({x}^{9}{y}^{5}+{x}^{5}{y}^{9})
\\[3mm]&+{\frac {7566301}{13948526592000}}({x}^{8}{y}^{6}+{x}^{6}{y}^{8})
\\[3mm]&+{\frac {212007212778877515989524910307665642456}{
286576587590801020714633877109661457031377625}}
{x}^{7}{y}^{7};
\end{align*}

\begin{align*}
&\Gamma^{\Lie}_{15}(x,y)=\\
\,&=
{\frac {1}{74724249600}}({x}^{14}{y}^{1}+{x}^{1}{y}^{14})
\\[3mm]&+{\frac {1}{3558297600}}({x}^{13}{y}^{2}+{x}^{2}{y}^{13})
\\[3mm]&+{\frac {2539}{871782912000}}({x}^{12}{y}^{3}+{x}^{3}{y}^{12})
\\[3mm]&+{\frac {
\left(\begin{array}{l}
\phantom{+\,00000000000000000000000000000000000000000000000}6\cdot10^{48}\\
+\,649786271967048955772044693373139641994744816961
\end{array}\right)
}{
\left(\begin{array}{l}
\phantom{+\,000000000000000000000000000000000000000}368722311\cdot10^{48}\\
+\,334408788686773393216110793289900907143561164800
\end{array}\right)
}}({x}^{11}{y}^{4}+{x}^{4}{y}^{11})
\\[3mm]&+{\frac {
\left(\begin{array}{l}
\phantom{+\,0000000000000000000000000000000000000000}63539675\cdot10^{48\cdot2}\\
+\,177049942043736407508334515735755770609745356893\cdot10^{48}\\
+\,980819471090768039178249054110571285451802903379
\end{array}\right)
}{
\left(\begin{array}{l}
\phantom{+\,000000000000000000000000000000000}901179990029545\cdot10^{48\cdot2}\\
+\,428718183918894693730342545775180196979531661104\cdot10^{48}\\
+\,845166086556190579801394158637374455917897216000
\end{array}\right)
}}({x}^{10}{y}^{5}+{x}^{5}{y}^{10})
\\[3mm]&+{\frac {
\left(\begin{array}{l}
\phantom{+\,000000000000000000}171547210173971722827445817547\cdot10^{48\cdot3}\\
+\,158731057405060500907292864526365317351811500143\cdot10^{48\cdot2}\\
+\,234296372605725177399144212353368954720896224907\cdot10^{48}\\
+\,598509491845326719170918708821316040741501970831
\end{array}\right)
}{
\left(\begin{array}{l}
\phantom{+\,000000000000}971295787097604635818160263849761551\cdot10^{48\cdot3}\\
+\,154019031183387500535145962154730827642729015649\cdot10^{48\cdot2}\\
+\,740280703987385235238923509068650210122205356558\cdot10^{48}\\
+\,302555224140950596575082693099267943804364864000
\end{array}\right)
}}({x}^{9}{y}^{6}+{x}^{6}{y}^{9})
\\[3mm]&+{\frac {
\left(\begin{array}{l}
\phantom{+\,00000000000000000000000000000000}5671493461617548\cdot10^{48\cdot4}\\
+\,853621987120367629212557003017717845581747445783\cdot10^{48\cdot3}\\
+\,441881655787020451150146106662972262482462298769\cdot10^{48\cdot2}\\
+\,251164579404266971191426668019910548820190417317\cdot10^{48}\\
+\,188721284699183953507808922928091909537502803051
\end{array}\right)
}{
\left(\begin{array}{l}
\phantom{+\,0000000000000000000000000}20364015450827782518605\cdot10^{48\cdot4}\\
+\,630081136937808190741805219343692215234730604277\cdot10^{48\cdot3}\\
+\,194653471919913806009542300049722574486542668065\cdot10^{48\cdot2}\\
+\,756386569878027310355088616768038436283503499148\cdot10^{48}\\
+\,947339629768820473677940940597942068844188672000
\end{array}\right)
}}
({x}^{8}{y}^{7}+{x}^{7}{y}^{8});
\end{align*}

\begin{align*}
&\Gamma^{\Lie}_{16}(x,y)=\\
\,&=
{\frac {1}{149448499200}}({x}^{14}{y}^{2}+{x}^{2}{y}^{14})
\\[3mm]&+{\frac {1}{7116595200}}({x}^{13}{y}^{3}+{x}^{3}{y}^{13})
\\[3mm]&+{\frac {2539}{1743565824000}}({x}^{12}{y}^{4}+{x}^{4}{y}^{12})
\\[3mm]&+{\frac {157}{17435658240}}({x}^{11}{y}^{5}+{x}^{5}{y}^{11})
\\[3mm]&+{\frac{19}{553512960}}({x}^{10}{y}^{6}+{x}^{6}{y}^{10})
\\[3mm]&+{\frac {
\left(\begin{array}{l}
\phantom{+\,00000000000000000000000000000000000}1039344362106\cdot10^{48\cdot4}\\
+\,720482120590271345733001556438706166521066488208\cdot10^{48\cdot3}\\
+\,726625544052437884549419900300665368847830165478\cdot10^{48\cdot2}\\
+\,198635405960119950700109064795116626269153588831\cdot10^{48}\\
+\,240655249022713140251858157212550524022240480531
\end{array}\right)
}{
\left(\begin{array}{l}
\phantom{+\,0000000000000000000000000000}13584389095357793053\cdot10^{48\cdot4}\\
+\,631780482183200751759399356789111798277774545854\cdot10^{48\cdot3}\\
+\,636488620407855713853143925839382798933802859905\cdot10^{48\cdot2}\\
+\,162717394221351011042265156501583938830121202814\cdot10^{48}\\
+\,851034123647812897279838265566424510922940416000
\end{array}\right)
}}
({x}^{9}{y}^{7}+{x}^{7}{y}^{9})
\\[3mm]&+{\frac {
\left(\begin{array}{l}
\phantom{+\,0000}21312541008580826861428594864343359866245654\cdot10^{48\cdot3}\\
+\,133759893174786848640898706040925391107362008393\cdot10^{48\cdot2}\\
+\,516262196621988387021856747472037509392076751209\cdot10^{48}\\
+\,273452853756327705607893223372184874826949407441
\end{array}\right)
}{
\left(\begin{array}{l}
\phantom{+\,000000000000000000000000000000000000000000000}217\cdot10^{48\cdot4}\\
+\,963975520630356755781078497827921973491487950724\cdot10^{48\cdot3}\\
+\,469334556024266318816004856827827803373892270751\cdot10^{48\cdot2}\\
+\,183288576930080464558924755845753928235607454138\cdot10^{48}\\
+\,885087283472978603661406075692737456207793920000
\end{array}\right)
}}
{x}^{8}{y}^{8}.
\end{align*}

\qedexer
\end{example}
\begin{remark}
Almost all computations in this paper were done with Maple 2015 in an x86-64, 6GB RAM environment.
Even $\Gamma^{\Lie}_{n_1,n_2}\,:\, n_1+n_2=14$ and  $\Theta^{\Lie}_9$ (!) were obtained in this way,
however, the latter case already required special assistance.
For $\Gamma^{\Lie}_{n_1,n_2}\,:\, 15\leq n_1+n_2\leq16$, the QSopt\_ex rational solver
(developed by David Applegate, William Cook, Sanjeeb Dash, Daniel Espinoza, see \cite{ACDE};
development code published by Daniel Espinoza,
cf.~also \cite{Esp}; install version  by Jon Lund Steffensen,
\verb"https://github.com/jonls")
was invoked in an x86-64, 22GB RAM environment, making the optimization more straightforward.
\qedremark
\end{remark}

In Examples \ref{ex:data1} and \ref{ex:data2},
the  data for $n=n_1+n_2\leq10$ shows some interesting patterns.
Yet, we are left with very few clues regarding the shape of minimal
presentations for $n_1\sim n_2$ in general.
Some of the obvious guesses turn out to be wrong:
For example, Example \ref{ex:data1} suggests that the minimal presentations for $\BCH^{\Lie}_{s,s}(X_1,X_2)$ can be
reduced to monomials of shape $[X_1,[\ldots]]$.
However, this  fails for $s=6$.
Similarly, Example \ref{ex:data2} suggests that the face-defining linear functional
\[ X^*_{{112211221122}}   =-\frac1{32} X^*_{{1212121212}} \]
could be a good dual solution with respect to $\BCH^{\Lie}_{6,6}(X_1,X_2)$.
This also fails.
(But in Section \ref{sec:bchupper}, essentially such linear functionals were utilized to give the upper estimate for the convergence radius.)
Some other patterns seem to continue
(like that the critical face is a hyperface if $n=n_1+n_2$ is odd)
but there is no obvious guarantee to extend to all degrees.
In general, $n=n_1+n_2\leq10$ means still relatively low degrees;
but in the case of higher degrees, $n=n_1+n_2\geq11$, results start to get difficult to print out.

\snewpage

\end{document}